\documentclass{article}
\usepackage[export]{adjustbox}
\usepackage{graphicx}%
\usepackage{amsmath,amssymb,amsfonts}%
\usepackage{amsthm}%
\usepackage{mathrsfs}%
\usepackage[title]{appendix}%
\usepackage{xcolor}%
\usepackage{algorithm}%
\usepackage{algorithmicx}%
\usepackage{algpseudocode}%
\usepackage{amsmath,bm}
\usepackage{physics}
\usepackage{mathtools}
\DeclarePairedDelimiterX{\inp}[2]{\langle}{\rangle}{#1, #2} 
\usepackage{acronym}
\newtheorem{theorem}{Theorem}

\newtheorem{problem}{Problem}

\newtheorem{remark}{Remark}%
\newtheorem{corollary}{Corollary}
\newtheorem{definition}{Definition}%
\usepackage{hyperref}
\usepackage{cleveref}
\usepackage{caption}
\usepackage{subcaption}
\usepackage{orcidlink}
\def\niceparagraph#1{\vspace{5pt} \noindent \textbf{#1}}

\newcommand{\OFF}[1]{}

\algnewcommand{\Initialize}[1]{%
  \State \textbf{Initialize:}
  \Statex \hspace*{\algorithmicindent}\parbox[t]{.8\linewidth}{\raggedright #1}
} 

\algnewcommand{\Setting}[1]{%
  \State \textbf{Setting:}
  \Statex \hspace*{\algorithmicindent}\parbox[t]{.8\linewidth}{\raggedright #1}
} 

\title{Koopman Regularization and Independent Koopman Eigenfunction Approximation}
\author{Ido Cohen\thanks{Faculty of Electrical and Computer Engineering, Technion - Israel Institute of Technology, 3200003, Haifa, Israel \orcidlink{0000-0001-8142-9092} orcid 0000-0001-8142-9092
       {\tt\small ido.coh@gmail.com, idoc@technion.ac.il}}
}
\date{\today}

\begin{document}
\maketitle
\begin{abstract}
The geometry induced from a dynamical system dictates its latent structures, conservation laws, and intrinsic coordinates. This paper provides two procedures to extract the governing equations of a dynamical system from its samples, \acf{KR} and \acf{IKEA}. These procedures leverage the prior knowledge from the Koopman Eigenfunction Space. Both techniques extract a functionally independent set of Koopman Eigenfunctions, which are the intrinsic coordinates. Consequently, despite its finite cardinality, the governing equations are immaculately restored. Thus, the principle of parsimony is implemented.

Samples can be taken either from a vector field or from orbits of the system. \acs{KR} extracts the governing law from samples of the vector field, and \acs{IKEA} extracts it from samples of trajectories of the system. Both algorithms are based on constrained optimization  methods for learning the governing equations from sparse and corrupted samples of the vector field or trajectories. These procedures yield promising results in denoising, generalization, and dimensionality reduction and achieve better results than \acs{DMD}, \acs{EDMD}, \acs{PINN}, and \acs{SINDy} in dynamical system restoration.
\end{abstract}


\section{Introduction}
Extracting the main coherent structures of a system is a technique for analyzing and reconstructing dynamical systems. For example, in linear dynamical systems, the stable, unstable, and center modes are the main geometric structures. These modes are the eigenvectors spanning the corresponding subspaces. They describe the system's behavior accurately, qualitatively, and quantitatively  \cite{lawrence1991differential}. Consequently, finding these eigenvectors and their corresponding eigenvalues is tantamount to revealing the governing equations in linear systems.

This concept was applied to nonlinear systems. The seminal work by Schmid  \cite{schmid2010dynamic} presents a spatio-temporal decomposition of orbits into coherent modes that evolve with time, known as \ac{DMD}. Due to its simplicity and the lack of prior knowledge requirements, this decomposition has been one of the most important cornerstones of many works in dynamical system analysis for the last two decades. This method was applied to many disciplines \cite{kutz2016dynamic}, and was extended to different variations e.g. \ac{EDMD} \cite{tu2013dynamic}, tls\ac{DMD} \cite{hemati2017biasing}, fb\ac{DMD} \cite{dawson2016characterizing}, Symmetric \ac{DMD} \cite{cohen2021modes}, optimized \ac{DMD}  \cite{askham2018variable}, Constrained \ac{DMD} \cite{azencot2019consistent}. Revealing the geometric structures of a nonlinear dynamical system is not trivial, and therefore, \ac{DMD} must be adapted for almost every single system. Even then, the system representation is not necessarily concise. \ac{EDMD} tries to generalize the \ac{DMD} representation by adding nonlinear functions of the state-space coordinates, but results in a redundant representation \cite{li2017extended}. The main variations of \ac{DMD} are summarized in \cite{schmid2022dynamic}.

The authors in \cite{rowley2009spectral} established that \ac{DMD} approximates the modes and the eigenvalues of the Koopman operator. Koopman theory has become popular in recent decades in the dynamical systems community as a system identification tool. Many algorithms strive to find the dynamic structure of a system from samples via Koopman theory, from linear approximation approaches (\ac{DMD} family \cite{schmid2022dynamic}) through dictionary-based algorithms (such as \ac{SINDy} \cite{brunton2016discovering}) to advantaging the computational power of NN (\ac{KAE} \cite{lusch2018deep}).

A sparse representation, such as \ac{SINDy}, circumvents this redundancy \cite{brunton2016discovering,9454411}. However, even these methods do not take into account the geometry of the dynamics since they are dictionary-based. The accuracy of these algorithms highly depends on the optional functions from the dictionary. Theoretically, a dictionary can recover any dynamical system flawlessly if it includes a subalgebra in $C^1$ which is essential to restore any functions in this space. This demand is not realistic. Restoring systems from samples requires a high level of flexibility and adaptability. These attributes, which are expected from a subalgebra, inherently exist in neural networks.

An autoencoder is one of the most common NN structures in manifold learning in data science. Leveraging this structure for dynamical system identification results in \acf{KAE} \cite{lusch2018deep}. This procedure can be described as a version of dictionary-based algorithm, where the encoder extracts the functions adaptively from the samples. The main \ac{DMD} cost function is the Frobenius norm of the approximated dynamical system error accumulating in one time step forward. In the cost function suggested in \cite{lusch2018deep}, the number of time steps whose error is accumulated is a parameter. Azencot et al. constrain consistency by minimizing the error one time step forward and one step backward \cite{azencot2020forecasting}, as suggested in \cite{azencot2019consistent}. The considerations of initializing this algorithm are discussed in \cite{miller2023eigenvalue}, which focuses on the reconstruction of stable or bounded systems. The authors in \cite{nayak2025temporally} presented better robustness to noise by minimizing the Hilbert Transform of the Frobenius norm of the temporal error. Unfortunately, even though \ac{KAE} achieves better results than the \ac{DMD} family, it inherits its flaws. The adaptivity of the NN, which is necessary for system identification, is also exploited to recognize the noise as part of the system, in particular, when the dimensionality of the latent space in the NN is higher than that of the system. Furthermore, including more data in the cost function of \ac{DMD} may result in the ``\ac{DMD} Paradox''\footnote{\ac{DMD} Paradox is when adding more data samples or increasing the number of the dominant modes reduces the \ac{DMD} cost function but does not improve the reconstruction's accuracy.} \cite{cohen2021modes} (Theorem 3.4). 

Another approach to inferring the governing law in nonlinear dynamics, based on deep learning,  is to embed the physical laws in the trained model. An algorithm from this approach is \ac{PINN}, which restores unknown parameters in dynamical systems \cite{raissi2017physics}. In practice, this method solves a \ac{PDE} with unknown parameters based on samples of the solution. Therefore, to use this kind of algorithm, the dynamical systems should be known in advance up to some parameters. Although this method yields better results compared to those of \ac{DMD} and its family, the physical laws behind the scenes are not always known. Therefore, dynamical system restoration relying on inherent and common attributes may lead to more robust results.

None of these methods leverage the functionally independent structure of the Koopman eigenfunction space.

This work extracts the main structures of the dynamics and reveals the physical laws behind the scenes. The theory of Bernard Osgood Koopman combines these two inherent attributes. Koopman argued that, given a nonlinear dynamical system, there exists a set of measurements of the state space that behaves linearly under the governing law \cite{koopman1931hamiltonian}. This measurement is termed \acf{KEF}, and its limitations are thoroughly discussed in \cite{cohen_gilboa_2023}.

For decades, researchers mistakenly considered the \ac{KEF} space functionally infinite-dimensional. Recently, however, it was proven that a finite set of KEFs generates the entire Koopman eigenfunction space \cite{cohen2024functional}. The cardinality of this generating set is at most the dimension of the system. Moreover, from this set, one can restore the governing equations and reveal the conservation (physical) laws of the system.

To recapitulate, the conclusions from \cite{cohen2024functional} are as follows. Let us consider the following nonlinear dynamical system
\begin{equation}\label{eq:introDynamics}
    \dot{\bm{x}}=p(\bm{x}), \,\, \bm{x}\in\mathcal{D}\subset\mathbb{R}^N, \,\, t\in I=[0,T]
\end{equation}
where $p:\mathbb{R}^N\to \mathbb{R}^N$ in $C^1$ (and therefore $\bm{x}(t)\in C^2$). Given this dynamical system, we are looking for the general form solution of \ac{KPDE}, defined by
\begin{equation}\label{eq:introKPDE}
    \nabla^T \phi(\bm{x})p(\bm{x})=\lambda \phi(\bm{x})
\end{equation}
where $\phi:\mathbb{R}^N\to \mathbb{R}$ and $\lambda$ is some constant in $\mathbb{C}$. From the general form given in \cite{cohen2024functional}, one can conclude the following.

\begin{enumerate}
    \item There are at most $N$ functionally independent \acp{KEF}.
    \item This set precisely restores the governing equations of the system.
    \item For any $N$-dimension dynamical system, there are at most $N-1$ functionally independent conservation laws.
    \item The conservation (physical) laws in 3 are derived from the set of \acp{KEF} in 1.
\end{enumerate}
Thus, the theory of Koopman not only reveals the inherent geometry of the system but also extracts the physical knowledge from this geometry.

Trying to find a Koopman eigenfunction by solving \eqref{eq:introKPDE} with a neural network may entail two main problems. The first is that the NN converges to the trivial eigenfunction, meaning a constant Koopman Eigenfunction with the corresponding eigenvalue ($\lambda=0$). The second is that even if it converges to non-trivial solutions, the NN results are not necessarily the most informative, as these eigenfunctions may be functionally dependent. As a result, the recovered governing law may be ill-posed or not robust enough to additive noise or measurement errors.

The former can be solved by finding an equivalent set of functions instead of \acp{KEF}. The suggested function set is the \ac{UVM} as discussed in \cite{cohen2023minimal}. A solution to the latter necessitates a thorough consideration of the functional dimensionality of the Koopman Eigenfunction Space and how to form such a set with the same dimensionality. In this work, we discuss the conditions on a set of functions to be functionally independent. From these conditions, the process of finding solutions to \eqref{eq:introKPDE} becomes a constrained optimization problem, where the objective functional refers to the PDE in \eqref{eq:introKPDE} and the feasibility functional refers to the functionally independent constraint.



This work provides two procedures for system identification from samples that take the geometry of the dynamical system into consideration. These algorithms, \acf{KR} and \acf{IKEA}, recover the governing equations from samples of the vector field and orbits, respectively. Both algorithms are based on the hidden geometry of the dynamical system, thus revealing a functionally independent set of Koopman Eigenfunctions; namely, the minimum number of independent variables required to describe the state space at any given time.

\subsection{Main Contributions}
\begin{enumerate}
    \item {\bf Koopman Regularization (\acs{KR})} -- is an algorithm that restores and generalizes a vector field from noisy or sparse vector field samples. By finding a set of at most $N$ Koopman eigenfunctions, it finds the parsimonious representation of the dynamical system. 

    \item {\bf Independent Koopman Eigenfunction Approximation (\acs{IKEA})} --  restores the governing law from samples of orbits. In the same vein as \ac{KR},  it finds the parsimonious representation of the dynamical system.

    \item {\bf{Condition for Linear Independence}} -- Given a set of vectors, $\{v_j\}_{j=1}^K$, we formulate algebraic and geometric 
    criteria derived from the Gershgorin circle theorem to determine whether the set is linearly independent.
    
    \item \label{item:CondFI} {\bf{Condition for Functional Independence}} -- Given a set of functions, $\{f_j(\bm{x})\}_{j=1}^K$ where $f_j:\mathbb{R}^N\to \mathbb{R}$, we formulate algebraic and geometric criteria derived from the Gershgorin circle theorem to determine whether the set is functionally independent.
    
    \item {\bf{Finding Functionally Independent Set of \acp{KEF}}} -- We solve the \ac{KPDE} (Eq. \eqref{eq:KEFPDE}) under the condition established in \Cref{item:CondFI}.

    \item {\bf Applications} -- \ac{KR} and \ac{IKEA} yield promising results for denoising, generalization, and dimensionality reduction.
    
\end{enumerate}

The paper is organized as follows. \Cref{sec:Preparatory} summarizes the necessary theoretical background and \Cref{sec:towards} defines the scope of this work. In \Cref{sec:sysReConstOpt}, we formulate the problem of finding \acp{KEF} as a constrained optimization problem. In \Cref{sec:penMethBack}, an algorithm is suggested to solve this optimization problem. In \Cref{sec:res}, we demonstrate \ac{KR} and \ac{IKEA} on samples of linear and nonlinear systems. Finally, in \Cref{sec:conclusions} we conclude our work.

\section{Preparatory Theory}\label{sec:Preparatory}
Essential mathematical notation, definitions, and theorems for this work are listed below.
\subsection{Notations}
\Cref{tab:Notations} summarizes the notations of differential, geometric, and algebraic operators relevant to this work.
\begin{table}[htbp]\caption{Necessary notations for this work.}
\begingroup
\renewcommand{\arraystretch}{1.5} 
    \begin{tabular}{|l| l|}
        \hline
        \multicolumn{1}{|c|}{\textbf{Notation}} & 
        \multicolumn{1}{c|}{\textbf{Meaning}} \\ 
        \hline
        \renewcommand{\arraystretch}{5}
        $i$& Index of the data points goes from $1$ to $L$\\
        \hline
        $j,k$& Indices of vectors in a set or of functions in a set go from $1$ to $K\leq N$\\
        \hline
        $o$&Index of orbits goes from $1$ to $O$\\
        \hline
        $s$& Index of samples from an orbit goes from $1$ to $S$\\
        \hline
        $\{\bm{x}_i\}_{i=1}^L$& Point set in $\mathcal{D}\subset \mathbb{R}^N$\\
        \hline
        $\{\bm{p}_i\}_{i=1}^L$& The corresponding vector field at the points $\{\bm{x}_i\}_{i=1}^L$, i.e. $\bm{p}_i=p(\bm{x}_i)$\\
        \hline
        $\bm{x}_s^o$& The $s$-th sample from the $o$-th orbit\\
        \hline
        $t_s^o$& The time point associated to the sample $\bm{x}_s^o$\\
        \hline
        $\nabla f$ & Gradient of a differentiable function $f:\mathbb{R}^N\to \mathbb{R}$ \\
        \hline
        $\bm{f}$ & Bold letter refers to a vector (of variables or functions)\\
        \hline
       $\bm{1}$& $K(\leq N)$ dimensional vector where each entry takes the value $1$\\
    &$\quad\quad\quad\quad\quad\quad\quad\quad\quad\quad\quad\quad \bm{1} = \begin{bmatrix}
        1&1&\ldots&1
    \end{bmatrix}^T$\\
    \hline
       $J(\bm{f}(\bm{x}))$ & Jacobian matrix defined as\\
       & 
    {\quad\quad\quad\quad\quad\quad\quad\quad$J(\bm{f})=
    \begin{bmatrix}
            \nabla^T f_1(\bm{x})\\
            \vdots\\
            \nabla^T f_K(\bm{x})
        \end{bmatrix}
    $}
    \\
    \hline
    $\inp{\cdot}{\cdot}$ & Standard euclidean inner product\\
    \hline
    $\norm{\cdot}$ & Standard euclidean norm\\
    \hline
    $\dot{}$ & Time derivative\\
    \hline
    $^T$& (superscript $T$) Transpose\\
    \hline
    \end{tabular}
    \endgroup
    \label{tab:Notations}
\end{table}

\subsection{Definitions and Theorems}

\begin{definition} [Functionally Independent Set (FIS)] \label{def:funcInd} Let $\{f_j(\bm{x})\}_{j=1}^K$ be a set of differentiable functions from $E\subset \mathbb{R}^N$ to $\mathbb{R}$. This set is functionally independent in the neighborhood of a point $\bm{x}\in E$ if there is no differentiable function $g$ such that
    \begin{equation}
            f_j(\bm{x}) = g(f_1(\bm{x}),\cdots,f_{j-1}(\bm{x}),f_{j+1}(\bm{x}),\cdots,f_K(\bm{x}))
    \end{equation}
    in the neighborhood of $\bm{x}$ and for all $f_j(\cdot),\,j=1,\,\cdots\,K$. Equivalently, this set\\
    $\{f_j(\bm{x})\}_{j=1}^K$ is functionally independent in the neighborhood of $\bm{x}$ if the corresponding Jacobian matrix is full rank (\cite{zorich2016mathematical} Section 8.6) in this neighborhood.
    
\end{definition}

\begin{definition}[Nonlinear Dynamical System]
\label{def:NdynamicalSystem}
A (nonlinear) dynamical system is specified by a domain
$\mathcal{D} \subseteq \mathbb{R}^N$, a time interval $I = [0,T]$, and a
vector field $p : \mathbb{R}^N \to \mathbb{R}^N$ of class $C^1$, governing
the evolution
\begin{equation}\label{eq:dynamics}
    \dot{\bm{x}}(t) = p(\bm{x}(t)), \quad \bm{x}(t=0) \in \mathcal{D}, \ t \in I.
\end{equation}
Since $p \in C^1$, the Picard--Lindel\"of theorem guarantees, for every
initial condition $\bm{x}(0) = \bm{x}_0 \in \mathcal{D}$, a unique local
solution $\bm{x}(\cdot)$; we assume $\mathcal{D}$ and $I$ are such that
this solution exists on all of $I$. Because
$\dot{\bm{x}}(t) = p(\bm{x}(t))$ is then a composition of $C^1$
functions, $\bm{x}(\cdot) \in C^2(I; \mathbb{R}^N)$.
\end{definition}


\begin{definition}[Configuration Space \cite{Arnold1989-1}]
A configuration space is a $K$ dimensional manifold embedded in the state space. This space defines the reachable area of the dynamical system. The admissible initial condition set or the physical constraints shape this manifold and its dimension.
\end{definition}

\begin{definition}[Dimensionality of Dynamical Systems (Degree of Freedom)] \label{def:DOF} The dimensionality of a dynamical system is the minimal number of independent variables that are required to define its state uniquely at any given moment. Equivalently, the dimensionality of a dynamical system is the configuration space dimensionality. This number is also termed the \textbf{degree of freedom} of the system.
\end{definition}

\begin{definition}[Orbit of an initial condition]\label{def:Orbit} Given an initial condition, $\bm{x}(t=0)= \bm{x}_0$, the unique solution of \eqref{eq:dynamics} can be seen as a curve in $\mathbb{R}^N$. This trajectory is denoted $\mathcal{X}(\bm{x}_0)$, and termed the \emph{orbit of $\bm{x}_0$}.
    
\end{definition}

\begin{definition}[Measurement] A measurement is a function from $\mathcal{D}$ to $\mathbb{C}$. 
    
\end{definition}

\begin{definition}[Koopman Operator]\label{def:KoopmanOperator}
    Let $g(\bm{x})$ be a measurement. The Koopman operator $K_\tau$ acts on the infinite-dimensional vector space of measurements and admits the following,
\begin{equation}\label{eq:koopDef}\tag{\bf{KO}}
    K_\tau(g(\bm{x}(t)))=g(\bm{x}(t+\tau)),\quad t,t+\tau\in I,
\end{equation}
where $\tau>0$ and $\bm{x}(t)$ is the solution of Eq. \eqref{eq:dynamics}. This operator is linear \cite{koopman1931hamiltonian}.
\end{definition}

\begin{definition}[\acf{KEF}]\label{def:KEF}
    Assuming the initial condition $\bm{x}_0$, a measurement $\varphi(\bm{x})$, satisfying the following relation along the orbit $\mathcal{X}(\bm{x}_0)$
\begin{equation}\label{eq:KEFdiff}
    \dfrac{d\varphi (\bm{x})}{dt}=\lambda \varphi(\bm{x}), \quad \forall \bm{x}\in \mathcal{X}(\bm{x}_0)
\end{equation}
for some value $\lambda\in \mathbb{C}$, is a \acf{KEF}.
    
\end{definition}
A \ac{KEF} therefore behaves linearly under the dynamics $p$.

\begin{definition}[\acf{KPDE}] The \acf{KPDE} is formulated as follows,
\begin{equation}\label{eq:KEFPDE}
    \nabla^T \Phi (\bm{x}) p(\bm{x}) = \lambda \Phi(\bm{x}), \quad \forall \bm{x}\in \mathcal{D}. 
\end{equation}
where $\Phi(\bm{x})$ is a differentiable measurement, $\nabla$ denotes the gradient of $\Phi$ with respect to the state vector $\bm{x}$, and $\lambda\in \mathbb{C}$.
\end{definition}

\ac{KPDE} is the ODE in Eq. \eqref{eq:KEFdiff} after applying the chain rule. A solution to Eq. \eqref{eq:KEFPDE} is a solution to Eq. \eqref{eq:KEFdiff}, but not vice versa. 

There are local solutions to this equation (e.g., \cite{Folland1995Introduction}, Chapter 1, or \cite{Debnath2012}, Theorem 3.5.3). To write the general-form solution to \ac{KPDE}, two more definitions should be noted.

\begin{definition}[Conservation Laws]\label{def:ConLaw} Conservation laws are the solution of Eq. \eqref{eq:KEFPDE} associated with $\lambda = 0$. This type of solutions is denoted as $h(\bm{x})$, (namely, admitting $\nabla^T h(\bm{x})p(\bm{x})=0$).
    
\end{definition}

\begin{definition}[Unit Velocity Measurement (UVM)]\label{def:UVM} A unit velocity measurement is a differentiable measurement satisfying the following \ac{PDE}
\begin{equation}\label{eq:unitPDE}
    \nabla^T m(\bm{x}) p(\bm{x})=1,\quad \forall \bm{x}\in \mathcal{D}.
\end{equation}

\end{definition}
Following \Cref{def:ConLaw,def:UVM}, we can formulate the solution to \ac{KPDE}.
\begin{theorem}[General Form Solution of \acs{KPDE} \cite{cohen2024functional}]\label{theo:Gery}   
The solution of Eq. \eqref{eq:KEFPDE} is of the form 
\begin{equation}\label{eq:generalSolKPDE}
    \Phi(\bm{x})=f(h_1(\bm{x}),\ldots,h_{N-1}(\bm{x}))\cdot e^{\lambda m(\bm{x})}
\end{equation}
where $m$ is a unit velocity measurement, $\{h_i\}_{i=1}^{N-1}$ is a set of functionally independent conservation laws, and $f$ is a differentiable function.
\end{theorem}

\Cref{theo:Gery} holds as long as $\mathcal{D}$ does not contain equilibria, singular points, closed, or recurrent orbits. For a more detailed discussion, see \cite{cohen2024functional}.

\begin{remark}[\ac{UVM} and \ac{KPDE} Relation]
    Let $\Phi(\bm{x})$ be a solution of Eq. \eqref{eq:KEFPDE} associated with eigenvalue $\lambda\ne0$. If $\Phi(\bm{x})\ne 0$ for all $\bm{x}\in\mathcal{D}$, a UVM derived from $\Phi(\bm{x})$ is given by \cite{cohen2023minimal}
\begin{equation}\label{eq:Relation}
    m(\bm{x})=\frac{1}{\lambda}\ln{\Phi(\bm{x})}.
\end{equation}
In the same vein, one can formulate a Koopman Eigenfunction with UVM as
\begin{equation}
    \Phi(\bm{x})=e^{\lambda m(\bm{x})}
\end{equation}
There are at most $N$ functionally independent unit velocity measurements \cite{cohen2023minimal, cohen2024functional}. In addition, Eq. \eqref{eq:unitPDE} is a chain rule applied to the ODE 
\begin{equation}
\frac{dm(\bm{x})}{dt}=1.    
\end{equation}
Therefore, a the graph of the solution to Eq. \eqref{eq:unitPDE} is a surface on which the velocity of the dynamics is one\footnote{This function is denoted by $m$ as a shorthand for $\mu o v \acute\alpha \delta \alpha$ ``unit'' in Greek.}.

\end{remark}

\subsection{Parsimonious Representation via \texorpdfstring{\acp{KEF}}{}}
The parsimonious representation of dynamical systems is finding the most informative variable set that describes the system with negligible\footnote{The definition of ``negligible'' depends on the application and the requirements at hand.} error. The most informative set is a set of functionally independent variables with maximal cardinality.  This statement coincides with the dimensionality of dynamical systems (\Cref{def:DOF}). Here, we are laying the theoretical foundation for restoring a dynamical system with independent functionally \acp{KEF} or \acp{UVM}.

\begin{definition}[Dynamic Reconstruction via \acp{KEF}]\label{def:reconMiniSet}
Let $\{\Phi_i(\bm{x})\}_{i=1}^N$ be a functionally independent set of \acp{KEF}. The dynamical system reconstruction is given by \cite{cohen2021latent}
\begin{equation}\label{eq:synReconKEF}
    \hat{p}(\bm{x}) =J(\bm{\Phi})^{-1}\Lambda \bm{\Phi}
\end{equation}
where $\Lambda=\text{diag}[\lambda_1,\ldots,\lambda_N]$, and $\bm{\Phi}(\bm{x})=[\Phi_1(\bm{x}),\ldots,\Phi_N(\bm{x})]^T$. This reconstruction exists as long as the Jacobian is full-rank.     
\end{definition}

\begin{definition}[Dynamic Reconstruction via UVM] Given a functionally independent set of $N$ unit velocity measurements, one can restore the dynamics $\hat{p}$ as
\begin{equation}\label{eq:dynReconUnit}
    \hat{p}(\bm{x}) = J(\bm{m})^{-1}\bm{1}.
\end{equation}
\end{definition}

\begin{remark}[Dimensionality Reduction via Koopman Theory]\label{remark:dim}  Dynamical system reconstruction from \acp{KEF} or from \acp{UVM} (Eqs. \eqref{eq:synReconKEF} and \eqref{eq:dynReconUnit}) are entailed from the general form of \ac{KEF}, proved in \Cref{theo:Gery}. This form is valid when the configuration space contains an $N-1$-dimension hypersurface embedded in $\mathbb{R}^N$ as a boundary condition for \eqref{eq:KEFPDE}. However, if the configuration space is a $K(<N)$-dimension manifold, then the dimensionality of the system is practically $K$ \cite{Arnold1989-1}. In that case, the general form of \ac{KEF} is
\begin{equation}
    \phi(\bm{x})=f(h_1(\bm{x}),\ldots,h_{K-1}(\bm{x}))e^{\lambda m(\bm{x})}
\end{equation}
where $f$, $\{h_i\}_{i=1}^{K-1}$, and $m$ are the same as in \eqref{eq:generalSolKPDE}.

As a result
\begin{enumerate}
    \item The cardinality of the largest functionally independent set of \acp{KEF} or \acp{UVM} is $K$.
    \item These sets of \acp{KEF} or \acp{UVM} can be interpreted as generalized coordinate systems. Therefore, the vector field belongs to the respective spans as follows
    \begin{equation}
        p(\bm{x})\in \text{Span}\{\nabla \phi_1(\bm{x}),\ldots,\nabla\phi_{K}(\bm{x})\}
    \end{equation}
    or equivalently
    \begin{equation}
        p(\bm{x})\in \text{Span}\{\nabla m_1(\bm{x}),\ldots,\nabla m_{K}(\bm{x})\}.
    \end{equation}
    \item The dynamical system reconstructions Eqs. \eqref{eq:synReconKEF} and \eqref{eq:dynReconUnit} become
    \begin{equation}\label{eq:UVMDR}
        \hat{p}(\bm{x}) =J^T(\bm{\Phi})\left(J(\bm{\Phi})J^T(\bm{\Phi})\right)^{-1}\Lambda \bm{\Phi}
    \end{equation}
    or alternatively
    \begin{equation}\label{eq:dynReconUnitK}
        \hat{p}(\bm{x}) = J^T(\bm{m})\left(J(\bm{m})J^T(\bm{m})\right)^{-1}\bm{1}.
    \end{equation}
\end{enumerate}    
\end{remark}
In \Cref{remark:dim}, the importance of functional independence is emphasized. The closer $JJ^T$ is to singularity, the less informative the set is. The remainder of this section discusses sufficient conditions for linear and functional independence.

\subsection{Conditions for Linear and Functional Independence}\label{subsec:cond}
Functional independence in the domain $E$ (following \Cref{def:funcInd}) is equivalent to linear independence of gradients in this domain. Here, a sufficient condition for linear and functional independence is established.
\subsubsection{Sufficient Conditions for Linear Independence}
Recall the Gershgorin theorem.
\begin{theorem}[Gershgorin circle theorem \cite{gershgorin1931uber}]\label{theo:gersh} Given a $K\times K$ matrix $A$, where $[A]_{j,k}=a_{j,k}$. The eigenvalues of $A$ are in the following domain in $\mathbb{C}$
\begin{equation}\label{eq:gershCir}
    \bigcup_{j=1}^K B\left(a_{j,j},\sum_{k=1,k\ne j}^K\abs{a_{j,k}}\right)
\end{equation}
where $B(a,r)$ is a ball in $\mathbb{C}$ centered in $a$ with radius $r$.  
\end{theorem}

\begin{corollary}[Sufficient Conditions for Linear Independence]\label{coro:sufConLin}
Let $\{v_j\}_{j=1}^K$ be a set of vectors in $\mathbb{R}^N$ where $K\leq N$. This set is linearly independent if 
\begin{equation}
    \sum_{j=1,j\ne k}^K \abs{\inp{v_j}{v_k}}<\norm{v_k}^2,\quad k=1,\ldots,K
\end{equation}
or 
\begin{equation}
    \sum_{j=1,j\ne k}^K\frac{\abs{\inp{v_j}{v_k}}}{\norm{v_j}\norm{v_k}}<1,\quad k=1,\ldots,K.
\end{equation}

\end{corollary}
\begin{proof} The set $\{v_j\}_{j=1}^K$ is linearly independent if the rank of the matrix 
\begin{equation}
    V = \begin{bmatrix}
        v_1&\cdots &v_K
    \end{bmatrix}
\end{equation}
is full. Equivalently, the rank of $V^TV$ is $K$. The matrix $V^TV$ is a Gram matrix, meaning it is symmetric, hermitian, and positive semi-definite. Excluding the origin from Gershgorin circles guarantees $V^TV$ does not have zero eigenvalues. To check whether the origin includes in the circles or not, one has to compare the radii and the centers in Eq. \eqref{eq:gershCir}. In our case, the $j$'th circle's center is at $\norm{v_j}^2$ and its radius is $\sum_{k=1,k\ne j}^K \abs{\inp{v_j}{v_k}}$, respectively. Therefore, if  
\begin{equation}\label{eq:gershFirstCon}
    \sum_{j=1,j\ne k}^K \abs{\inp{v_j}{v_k}}<\norm{v_k}^2,\quad k=1,\ldots,K
\end{equation}
the matrix $V^TV$ is not singular, therefore, the set $\{v_j\}_{j=1}^K$ is linearly independent.

In the same vein, we could start this proof by defining the matrix
\begin{equation}
    \hat{V}=\begin{bmatrix}
        \hat{v}_1&\cdots &\hat{v}_K
    \end{bmatrix}
\end{equation}
where $\hat{v}_j = \dfrac{v_j}{\norm{v_j}}$. The corresponding condition to \eqref{eq:gershFirstCon}
is 
\begin{equation}
    \sum_{j=1,j\ne k}^K\frac{\abs{\inp{v_j}{v_k}}}{\norm{v_j}\norm{v_k}}<1,\quad k=1,\ldots,K.
\end{equation}

\end{proof}

\subsubsection{Sufficient Conditions for Functional Independence}

\begin{theorem}[Sufficient Conditions for Functional Independence]\label{theo:sufConFun}
    Let $\{f_j\}_{j=1}^K$ ($K\leq N$) be a set of functions from $\mathcal{D}\subset \mathbb{R}^N$ to $\mathbb{R}$ in $C^1$, and let the gradient of these functions do not vanish in $\mathcal{D}$. This set is functionally independent in $\mathcal{D}$ if 
    \begin{equation}
        \sum_{j=1,j\ne k}^K \abs{\inp{\nabla f_j(x)}{\nabla f_k(x)}}<\norm{\nabla f_k(x)}^2, \quad k=1,\ldots,K
    \end{equation}
    or
    \begin{equation}\label{eq:geoCond}
        \sum_{j=1,j\ne k}^K \frac{\abs{\inp{\nabla f_j(x)}{\nabla f_k(x)}}}{\norm{\nabla f_j(x)}\norm{\nabla f_k(x)}}<1, \quad k=1,\ldots,K
    \end{equation}
    $\forall x\in \mathcal{D}$.
\end{theorem}
\begin{proof}
    According to \Cref{def:funcInd} the Jacobian $J$ must be a full rank matrix at any point $x$ in $\mathcal{D}$. If conditions of \Cref{coro:sufConLin} hold at any point in $\mathcal{D}$, the set $\{f_k(x)\}_{k=1}^K$ is functionally independent.
\end{proof}

\section{Towards \texorpdfstring{\ac{KR}}{} and \texorpdfstring{\ac{IKEA}}{}}\label{sec:towards}
Let us consider a dynamical system in the form given in \Cref{def:NdynamicalSystem}. This work presents two methods for reconstructing a dynamical system from samples. The samples can be obtained either from the vector field or from orbits corresponding to different initial conditions. In both methods, the reconstruction is on a domain that contains the input data and excludes equilibria, singular points, closed orbits, and recurrent orbits. These conditions ensure that \Cref{theo:Gery} holds. Consequently, a set of \acp{KEF} exists and can be reconstructed from the given data.

\begin{figure}[phtb!]
    \centering 
    \includegraphics[trim=100 100 90 75, clip,width=0.85\textwidth,valign = t]{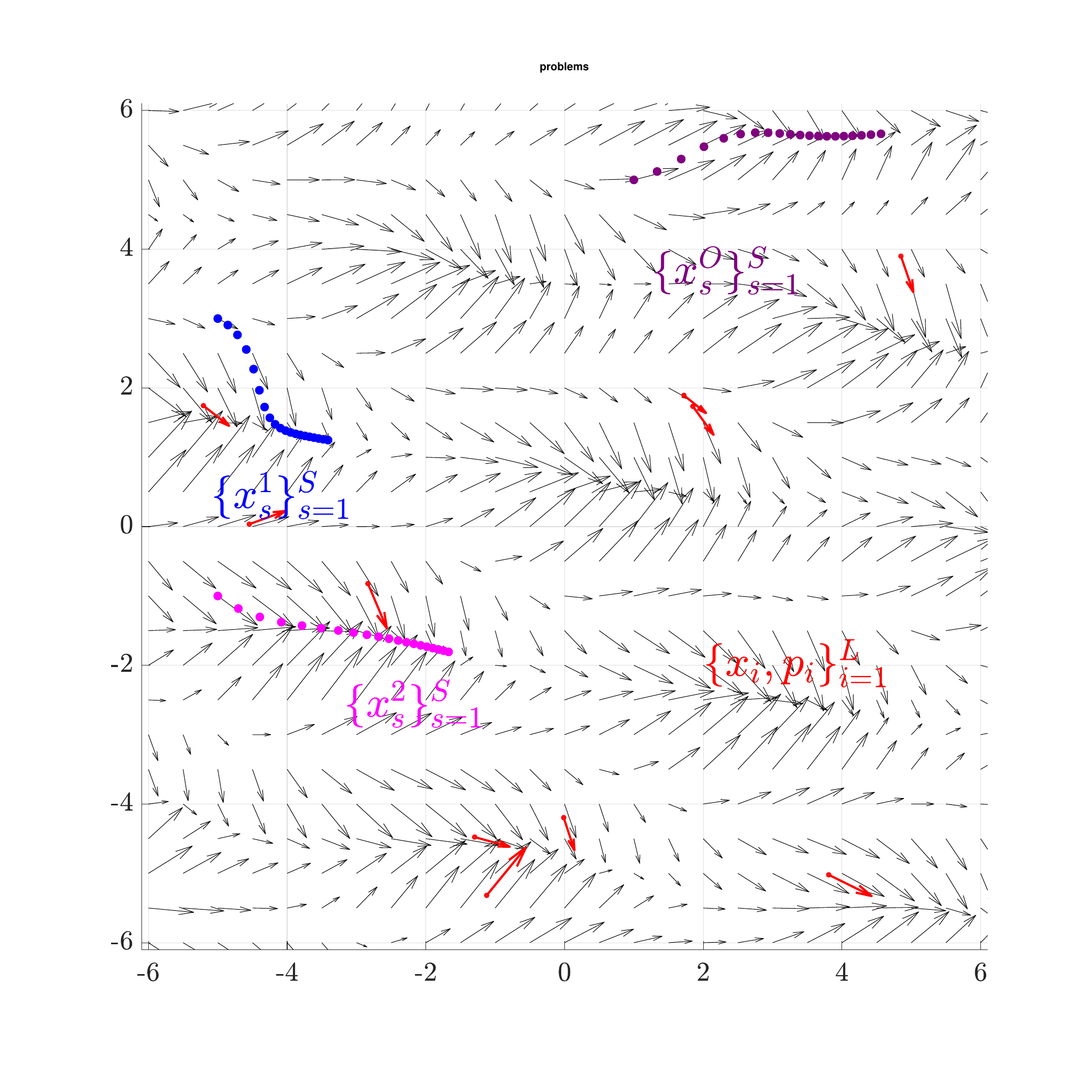}
    \caption{{\bf{Input Data Configurations}} -- System restoration either from $O$ sets of sampled orbits $\{\bm{x}_s^o\}_{s=1,o=1}^{S,O}$ of a system differently initiated, or from a set of samples of a vector field $\{\bm{x}_i,\,\bm{p}_i\}_{i=1}^L$.}
    \label{fig:problems} 
\end{figure}

An example illustrating these methods is shown in \Cref{fig:problems}. In this example, the objective is to reconstruct a two-dimensional dynamical system on the domain $[-6,6]\times[-6,6]$. The black arrows represent the dynamical system. To achieve this objective, the given data consist of either samples of the vector field, i.e., $\{\bm{x}_i,p_i\}_{i=1}^L$ (represented by the red arrows and dots), or samples from orbits corresponding to different initial conditions together with their corresponding times, i.e., $\{\bm{x}_s^o,t_s^o\}_{s=1,o=1}^{S,O}$ (with each orbit depicted in a different color). Next, we define the problems addressed by \ac{KR} and \ac{IKEA}, and their respective settings.

\subsection{Problem Settings for \texorpdfstring{\ac{KR}}{}}\label{subsec:PSKR}
\acf{KR} extracts the governing law from samples of the vector field. A formal formulation of the problem is given below.

\begin{problem}[\emph{Koopman Regularization}]\label{Pro:Main}
Let us consider a dynamical system as defined in \Cref{def:NdynamicalSystem}. Let $X = \{\bm{x}_i\}_{i=1}^L$ be a set of points in $\mathcal{D}\subset \mathbb{R}^N$, and let $P=\{\bm{p}_i\}_{i=1}^L$ denote the corresponding samples of the vector field $p$ at these points. The set of samples $P$ may be sparse, noisy, or embedded in a lower-dimensional manifold. Assuming that $\mathcal{D}$ does not contain singular points, closed or recurrent orbits of the dynamics, we seek to reconstruct the vector field $p$ at any point in this domain.
\end{problem}

We consider Problem \ref{Pro:Main} in three different settings. The first is the denoising of corrupted vector field samples. The second is the generalization of sparse samples. The third is the identification of the lower-dimensional dynamical system.

\paragraph{Noise Reduction} \label{sec:settingNoiseRe}
The goal is to reconstruct the vector field from corrupted samples. The given data consist of a set of points $\{\bm{x}_i\}_{i=1}^L$, where $\bm{x}_i\in \mathcal{D}$, and the corresponding corrupted samples of the vector field, $\{\tilde{\bm{p}}_i\}_{i=1}^L$. The samples are corrupted by additive noise, i.e.,
\begin{equation}
\tilde{\bm{p}}_i = \bm{p}_i+\bm{n}_i,
\end{equation}
where $\bm{p}_i=p(\bm{x}_i)$ is the noise-free sample and $\bm{n}_i$ denotes the additive noise.

\paragraph{Generalization} \label{sec:settingGe}
The generalization process is common in image processing (inpainting) \cite{shi2016weighted} or machine learning when the labeled data is sparse \cite{calder2018game}. In both cases, diffusion provides the intuitive solution. Borrowing this approach to generalize sparse samples in a vector field demands rethinking the functional to minimize. If the vector field is conservative, restoring the potential function is enough to ``inpaint'' the domain. On the other hand, when it is not conservative, there is not only one real function that describes the whole domain but $N$ functions. One can use a minimal set of Koopman eigenfunctions from which the dynamics can be restored. The geometry induced by these eigenfunctions dictates the vector field at any point in the domain. 

\paragraph{Dimensionality Reduction} \label{sec:settingDimRe}
The demand for parsimonious representations is not limited to data analysis. The growing size of neural networks raises the need for an accurate and compact dynamic representation in order to accelerate the training process. The authors in \cite{gilboa2024system} showed a high (temporal) correlation between the weights in the training process. They used this correlation for prediction in the training process. Additionally, the temporal correlation is essential for pruning. Here, the functional dependence plays the role of the temporal correlation. It is shown that a small set of functionally independent functions can restore a high-dimensional dynamic.

\subsection{Problem Settings for \texorpdfstring{\ac{IKEA}}{}}
\acf{IKEA} extracts the governing law from sampled orbits of different initial conditions. A formal formulation is given below.
\begin{problem}[\emph{Independent Koopman Eigenfunction Approximation}]\label{prob:IKEA}
Let us consider a dynamical system (\Cref{def:NdynamicalSystem}), and let $\{\mathcal{X}^o\}_{o=1}^{O}$ be a set of curves corresponding to different initial conditions $\{\bm{x}_0^o\}_{o=1}^{O}$. Let $\{\bm{x}_s^o\}_{s=1}^{S}$ denote the $S$ samples of the $o^{\text{th}}$ orbit, corresponding to the time points $\{t_s^o\}_{s=1}^{S}$, respectively. Assume that there exists a simply connected domain $\mathcal{D}$ containing the points $\{\bm{x}_s^o\}_{s=1,o=1}^{S,O}$ that does not contain singular points, closed or recurrent orbits of the dynamics. Given the sampled orbits, we seek to reconstruct the governing law of the dynamical system on the domain $\mathcal{D}$.
\end{problem}

In this context, we focus on extracting the governing law from noisy samples of the orbits. In practice, the input data is $\{\tilde{\bm{x}}_s^o,t_s^o\}_{s=1,o=1}^{S,O}$ where $\tilde{\bm{x}}_s^o$ is the clear state vector at time $t_s^o$ along the $o^{th}$ orbit with additive noise, meaning 
\begin{equation}\label{eq:ikeaNoise}
    \tilde{\bm{x}}_s^o = \bm{x}_s^o+n_s^o
\end{equation}
where $\bm{x}_s^o$ is the clear sample and $n_s^o$ is a noise drawn from some distribution. Please note that the data should not be sampled at equal intervals.

The formulations of Problems \ref{Pro:Main} and \ref{prob:IKEA} and their settings fit the form of ``inverse problem'' from the domain of variational calculus. Namely, the solutions to these problems are equivalent to the solutions of optimization problems under geometric constraints. In the next section, the objective functional and the feasible region are defined.

\section{System Reconstruction as a Constrained Optimization Problem}\label{sec:sysReConstOpt}
A functionally independent set of \acp{UVM} should satisfy Eq. \eqref{eq:unitPDE} and the condition in \Cref{theo:sufConFun}. To this end, the input data are incorporated into Eq. \eqref{eq:unitPDE} through an objective functional and into Eq. \eqref{eq:geoCond} through a feasible functional. The general form of the optimization problem for \ac{KR} and \ac{IKEA} is
\begin{equation}
\begin{split}
\min_{\bm{m}}\mathcal{O}(\bm{m})\quad \text{s.t.}\,\, \mathcal{F}(\bm{m})=0,
\end{split}
\end{equation}
where $\mathcal{O}(\bm{m})$ and $\mathcal{F}(\bm{m})$ are the objective and feasible functionals, respectively. These functionals are defined differently depending on the type of input data. In both cases, however, the functional $\mathcal{O}$ enforces Eq. \eqref{eq:unitPDE}, while $\mathcal{F}$ enforces the geometric condition in Eq. \eqref{eq:geoCond}.

\subsection{Objective Functional} 
Let $\{m_j\}_{j=1}^K$ be the set of \acp{UVM}. The objective functionals of \ac{KR} and \ac{IKEA} enforce the satisfaction of Eq. \eqref{eq:unitPDE} by this set. The specific form of the objective functional depends on the type of input data.

\subsubsection{\acf{KR}}
The input data for \ac{KR} are the pairs $\{\bm{x}_i,\bm{p}_i\}_{i=1}^L$. Each element of the \ac{UVM} set is required to satisfy Eq. \eqref{eq:unitPDE} at every point in the input data, i.e.,
\begin{equation}
\nabla^T m_j(\bm{x}_i)\bm{p}_i=1, \quad \forall i=1,\ldots,L,\,\, \forall j=1,\ldots,K.
\end{equation}
Consequently, the set of \acp{UVM} minimizes the following functional,
\begin{equation}
\mathcal{O}(\bm{m})=\sum_{i=1}^L\sum_{j=1}^{K}\norm{\nabla^T m_j(\bm{x}_i)\bm{p}_i-1}^2.
\end{equation}
In matrix form, this functional can be written as
\begin{equation}\label{eq:objFunctionalKR}
\mathcal{O}(\bm{m})=\sum_{i=1}^L\norm{J(\bm{m}(\bm{x}_i))\bm{p}_i-\bm{1}}^2.
\end{equation}

\subsubsection{\acf{IKEA}}
The input data for \ac{IKEA} are the pairs $\{\bm{x}_s^o,t_s^o\}_{s=1,o=1}^{S,O}$, where $\bm{x}_s^o$ is the $s^{\text{th}}$ sample of the $o^{\text{th}}$ orbit and $t_s^o$ is the corresponding time. As discussed in \cite{cohen2023minimal}, a \ac{UVM} is a mapping from the state space to the time axis, up to an additive constant. Therefore, for any two samples $\bm{x}_s^o$ and $\bm{x}_{s'}^o$ from the same orbit $o$, a \ac{UVM} satisfies
\begin{equation}
m_j(\bm{x}_{s'}^o) - m_j(\bm{x}_s^o)=t_{s'}^o-t_s^o.
\end{equation}
By incorporating all pairs of samples from all orbits and all \acp{UVM} into a single functional, we obtain
\begin{equation}\label{eq:UVMfunctionalIKEA}
\begin{split}
\mathcal{O}(\bm{m}(\bm{x}))=
\sum_{j=1}^{K}\sum_{o=1}^O\sum_{s=1}^{S}\sum_{s'=1}^{S}
\norm{m_j(\bm{x}_{s'}^o)-m_j(\bm{x}_{s}^o)-(t_{s'}^o-t_{s}^o)}^2.
\end{split}
\end{equation}

\subsection{Feasible Functional} 
The feasible functional determines whether the set of \acp{UVM} is functionally independent. As a result from \Cref{theo:sufConFun},  if the \ac{UVM} set, $\{m_j(\bm{x})\}_{j=1}^K$, is functionally independent in $\mathcal{D}$, the following inequation must be satisfied
\begin{equation}\label{eq:geoCondKRIKEA}
\sum_{j=1,j\ne k}^K
\frac{\abs{\inp{\nabla m_j(\bm{x})}{\nabla m_k(\bm{x})}}}
{\norm{\nabla m_j(\bm{x})}\norm{\nabla m_k(\bm{x})}}
<1, \quad k=1,\ldots,K
\end{equation}
at every point $\bm{x}\in\mathcal{D}$. However, unfortunately, we have to limit this demand only to sampled points in $\mathcal{D}$. Next, this condition is applied to the input data of \ac{KR} and \ac{IKEA}.

\subsubsection{\acf{KR}}
The input data for \ac{KR} are the pairs $\{\bm{x}_i,\bm{p}_i\}_{i=1}^L$. Applying the condition in Eq. \eqref{eq:geoCondKRIKEA} to the input data, we get
\begin{equation}\label{eq:geoCondKR}
\sum_{j=1,j\ne k}^K
\frac{\abs{\inp{\nabla m_k(\bm{x}_i)}{\nabla m_j(\bm{x}_i)}}}
{\norm{\nabla m_k(\bm{x}_i)}\norm{\nabla m_j(\bm{x}_i)}}
<1, \quad k=1,\ldots,K,
\end{equation}
at every point $\bm{x}_i$ in the input data.

We define the maximal bias from functional independence as
\begin{equation}\label{eq:maxBaisKR}
\max_{k,i}\left\{\sum_{j=1,j\ne k}^K \frac{\abs{\inp{\nabla m_k(\bm{x}_i)}{\nabla m_j(\bm{x}_i)}}}
{\norm{\nabla m_k(\bm{x}_i)}\norm{\nabla m_j(\bm{x}_i)}}
\right\}.
\end{equation}
By constraining this maximal bias, the minimum of Eq. \eqref{eq:geoCondKR} can be restricted to the feasible region. Hence, the chosen feasible functional ensures that the maximal bias is below a nonnegative threshold $\delta\,(<1)$, thereby ensuring that the \ac{UVM} set is functionally independent. Formally, the feasible functional is given by
\begin{equation}\label{eq:feaFunKR}
\mathcal{F}(\bm{m})=
\max\left\{\max_{k,i}\left\{ \sum_{j=1,j\ne k}^K 
\frac{\abs{\inp{\nabla m_k(\bm{x}_i)}{\nabla m_j(\bm{x}_i)}}}
{\norm{\nabla m_k(\bm{x}_i)}\norm{\nabla m_j(\bm{x}_i)}}
\right\}-\delta,0
\right\},
\end{equation}
where $\delta\in[0,1)$.

\subsubsection{\acf{IKEA}}
The input data for \ac{IKEA} are the samples of the orbits and their corresponding times, $\{\bm{x}_s^o,t_s^o\}_{s=1,o=1}^{S,O}$. In this case, the geometric condition is
\begin{equation}\label{eq:geoCondIKEA}
\sum_{j=1,j\ne k}^K
\frac{\abs{\inp{\nabla m_k(\bm{x}_s^o)}{\nabla m_j(\bm{x}_s^o)}}}
{\norm{\nabla m_k(\bm{x}_s^o)}\norm{\nabla m_j(\bm{x}_s^o)}}
<1, \quad k=1,\ldots,K,
\end{equation}
at every sampled point $\bm{x}_s^o$. The corresponding feasible functional is defined as
\begin{equation}\label{eq:feaFunIKEA}
\mathcal{F}(\bm{m})=
\max\left\{\max_{k,s,o}\left\{ \sum_{j=1,j\ne k}^K \frac{\abs{\inp{\nabla m_k(\bm{x}_s^o)}{\nabla m_j(\bm{x}_s^o)}}}
{\norm{\nabla m_k(\bm{x}_s^o)}\norm{\nabla m_j(\bm{x}_s^o)}}
\right\}-\delta,0
\right\},
\end{equation}
where $\delta\in[0,1)$.

\begin{remark}[Recapitulation of \ac{KR} and \ac{IKEA}]
\ac{KR} and \ac{IKEA} are both constrained optimization problems for extracting functionally independent \acp{UVM} from different types of data. As discussed above, the objective and feasible functionals are defined according to the type of input data and enforce the constraints and conditions given by Eq. \eqref{eq:unitPDE} and \Cref{theo:sufConFun}. \Cref{table:KRIKEA} summarizes the input data, output, objective functional, and feasible functional for both methods.

\begin{table}[pthb!]
        \caption{Main Components of \ac{KR} and \ac{IKEA}}
    \label{table:KRIKEA}

\hspace{-10em}
\begin{tabular}{|c ||c| c|}
 \hline
 &\acs{KR}& \acs{IKEA}\\
 \hline\hline
Input & $\{\bm{x}_i,\bm{p}_i\}_{i=1}^L$ & $\{\bm{x}_s^o,t_s^o\}_{s=1,o=1}^{s=S,o=O}$\\
\hline
Output & $\{m_j(\bm{x})\}_{j=1}^K$ & $\{m_j(\bm{x})\}_{j=1}^K$\\
\hline
$\mathcal{O}(\bm{m})$ & $\sum_{i=1}^L\norm{J(\bm{m}(\bm{x}_i))\bm{p}_i-\bm{1}}^2$ & $\sum_{j,o,s,s}\norm{m_j(\bm{x}_{s'}^o)-m_j(\bm{x}_{s}^o)-(t_{s'}^o-t_{s}^o)}^2$\\
\hline
$\mathcal{F}(\bm{m})$ & $\max\left\{\max_{k,i}\left\{\sum_{j=1,j\ne k}^K\frac{\abs{\inp{\nabla m_k(\bm{x}_i)}{\nabla m_j(\bm{x}_i)}}}{\norm{\nabla m_k(\bm{x}_i)}\norm{\nabla m_j(\bm{x}_i)}}\right\}-\delta,0\right\}$ & $\max\left\{\max_{k,s,o}\left\{\sum_{j=1,j\ne k}^K\frac{\abs{\inp{\nabla m_k(\bm{x}_s^o)}{\nabla m_j(\bm{x}_s^o)}}}{\norm{\nabla m_k(\bm{x}_s^o)}\norm{\nabla m_j(\bm{x}_s^o)}}\right\}-\delta,0\right\}$\\ 
\hline

\end{tabular}
 
\end{table}

\end{remark}

\section{Penalty Method with Backtracking}\label{sec:penMethBack}
In this section, the optimization process of \ac{KR} and \ac{IKEA} is explained. Although \ac{KR} and \ac{IKEA} employ different objective and feasibility functionals, they share the same optimization procedure. 

The aim is to minimize the objective functional $\mathcal{O}(\bm{m})$ subject to the constraint of functional independence. Constrained optimization problems in variational calculus are well studied \cite{OPTbao2004computing, OPTbao2003ground, OPTcaliari2009minimisation, OPTCOHEN20181138, OPTdemyanov2011exact, OPTdem2004exact, OPTekeland1999convex, OPTgarcia2001optimizing}. 
To unify the objective and feasible functionals into one loss function, the following formulation is suggested
\begin{equation}\label{eq:Lalpha}
    \mathcal{L}_\alpha (\bm{m})=\alpha \mathcal{O}(\bm{m})+\mathcal{F}(\bm{m}).
\end{equation}
where $\alpha$ is a nonnegative real number. To find a minimal set $\{m_i\}_{i=1}^K$, we find the minimizer of $\mathcal{L}_\alpha (\cdot)$ for a sequence of values of $\alpha$ increasing from $0$ to $1$. The minimization of $\mathcal{L}\alpha$ is performed subject to the feasibility condition $\mathcal{F}(\bm{m})<1$. Our method applies the gradient descent flow as long as it remains in the feasible region. Thus, the feasibility functional is invoked to ensure that the iterative optimization process remains within the feasible region.

\subsection{Approaching the Feasible Region}\label{subsubsec:FR}
Given an initial condition, we first have to set it in the feasible region by reducing the feasibility functional $\mathcal{F}$ below a prescribed threshold $\theta_{\mathcal{F}}$. Thus, if $\mathcal{F}(\bm{m})<\theta_{\mathcal{F}}$ then the geometric conditions in Eqs. \eqref{eq:geoCondKR}  and \eqref{eq:geoCondIKEA} are less than $\theta_{\mathcal{F}}+\delta$

\begin{subequations}
\begin{equation}
    \sum_{j=1,j\ne k}^K\frac{\abs{\inp{\nabla m_k(\bm{x}_i)}{\nabla m_j(\bm{x}_i)}}}{\norm{\nabla m_k(\bm{x}_i)}\norm{\nabla m_j(\bm{x}_i)}}<\theta_{\mathcal{F}}+\delta, \quad \forall k=1,\ldots,K,
\end{equation}
\begin{equation}
    \sum_{j=1,j\ne k}^K\frac{\abs{\inp{\nabla m_k(\bm{x}_s^o)}{\nabla m_j(\bm{x}_s^o)}}}{\norm{\nabla m_k(\bm{x}_s^o)}\norm{\nabla m_j(\bm{x}_s^o)}}<\theta_{\mathcal{F}}+\delta, \quad \forall k=1,\ldots,K,
\end{equation}
\end{subequations}
To fulfill the conditions, the upper bound must be less than $1$, or equivalently,
\begin{equation}\label{eq:thetaF}
    \theta_{\mathcal{F}}<1 - \delta.
\end{equation}
Since the feasibility functional admits a trivial minimizer for which the expressions in Eqs. \eqref{eq:feaFunKR} and \eqref{eq:feaFunIKEA} vanish, minimizing $\mathcal{F}$ provides a procedure for moving the initial condition toward the feasible region. 
This procedure is summarized in \Cref{alg:FR}. This optimization problem is the first step of the constrained optimization in Eq. \eqref{eq:Lalpha}, meaning where $\alpha=0$. The subsequent step is to minimize $\mathcal{L}_\alpha$ where $\alpha$ is positive.

\begin{algorithm}[phtb!]
\caption{Feasible Set Initialization}\label{alg:FR}
\begin{algorithmic}[1]

\Setting{Set $\eta_{th}>0$, $\eta>0$, $\beta\in(0,1)$, $\delta\in[0,1)$, $\theta_{\mathcal{F}}\,(<1-\delta)$, and $\bm{m}$ arbitrarily }
\While{$\eta > \eta_{th}$}
    \If{$\mathcal{F}(\bm{m})<\theta_{\mathcal{F}}$}
        \State{Return}
    \EndIf
    \State{$\bm{m} = \bm{m}-\partial \mathcal{F}(\bm{m})\eta$}
    \If{$\mathcal{F}(\bm{m})$ decreased}
        \State{Continue}
    \Else
        \State{Restore the previous point $\{\bm{m}\}$}
        \State{$\eta = \beta\cdot \eta$}
    \EndIf
\EndWhile
\end{algorithmic}
\end{algorithm}

\subsection{Minimization process of \texorpdfstring{$\mathcal{L}_\alpha $}{L} with a fixed \texorpdfstring{$\alpha$}{a}} \label{subsubsec:La}
The aim is now to perform gradient descent on $\mathcal L_\alpha$ within the feasible region. The suggested method here is to apply one step forward in the flow. If the resulting iterate remains feasible, proceed to another step. If the result is not feasible, restore the previous iterate and apply one smaller step in the flow. The stop condition is when the step size is smaller than a threshold. 
This process is summarized in \Cref{alg:Mini}. 
\begin{algorithm}[phtb!]
\caption{Minimizing $\mathcal{L}_\alpha$ under the constraint of $\mathcal{F}$}\label{alg:Mini}
\begin{algorithmic}[1]

\Setting{Set $\eta_{th}>0$, $\eta>0$, $\beta\in(0,1)$, and $\theta_{\mathcal{F}}$ according to Eq.\eqref{eq:thetaF} }

\While{$\eta>\eta_{th}$}
    \State $\bm{m} = \bm{m}-\partial \mathcal{L}_\alpha(\bm{m})\eta$
    \If{$\mathcal{L}_{\alpha}(\bm{m})$ decreased $\&\, \mathcal{F}(\bm{m})<\theta_{\mathcal{F}}$}
        \State{Continue}
    \Else
        \State{Restore the previous point $\{\bm{m}\}$}
        \State{$\eta = \beta\cdot \eta$}
    \EndIf
\EndWhile
\end{algorithmic}
\end{algorithm}

It should be noted that when the quantities in Eqs. \eqref{eq:geoCondKR} and \eqref{eq:geoCondIKEA} are below $\delta$, the feasibility functional $\mathcal{F}$ is inactive and therefore does not contribute to the gradient of $\mathcal{L}_\alpha$. Consequently, for any $\delta\in[0,1)$, the resulting set remains functionally independent, while $\delta$ controls the tolerance with which the gradient approaches orthogonality.

\subsection{\texorpdfstring{\ac{KR}}{} and \texorpdfstring{\ac{IKEA}}{} -- Optimization Process}
Now, the \ac{KR} and \ac{IKEA} optimization procedures combine the two stages described above: initialization in the feasible region (\Cref{subsubsec:FR}) followed by minimization of $\mathcal L_\alpha$ for increasing values of $\alpha$ (\Cref{subsubsec:La}). The initialization step is to locate the initial condition in the feasible region by invoking \ref{alg:FR} ($\alpha=0$). The next step is to repeatedly invoke \Cref{alg:Mini} for increasing values of $\alpha$ from $0$ to $\alpha_{th} \, (>0)$ by increments of $\Delta \alpha$.
This approach is summarized in \Cref{alg:KR}.

\begin{algorithm}[phtb!]
\caption{\acl{KR} and \acl{IKEA}}\label{alg:KR}
\begin{algorithmic}[1]

\Setting{Set $\eta_{th}>0$, $\eta>0$, $\beta\in(0,1)$, $\alpha=\Delta\alpha$, $\alpha_{th} $, $\delta\in[0,1)$, and $\theta_{\mathcal{F}}$ according to Eq.\eqref{eq:thetaF} }

\Initialize{Minimize $\mathcal{F}$ \Comment{Invoke \Cref{alg:FR} using $\eta_{th}$, $\eta$, $\beta$, $\delta$, and $\theta_{\mathcal{F}}$}}
\While{$\alpha \leq  \alpha_{th}$}
    \State{Minimize $\mathcal{L}_\alpha(\bm{m})$ \Comment{Invoke \Cref{alg:Mini} using $\eta_{th}$, $\eta$, $\beta$, $\alpha$, $\delta$, and $\theta_{\mathcal{F}}$}}
    \State{$\alpha\leftarrow\alpha+\Delta\alpha$}
\EndWhile
\end{algorithmic}
\end{algorithm}



\section{Results}\label{sec:res}
The goal of this section is to illustrate the theory discussed above. In addition, it elaborates the framework of the experiments and the settings such that one can repeat them and restore the results achieved here. We start with the NN configuration, then we discuss the results of \acl{KR} in different settings and statistics; the last are the results of \acl{IKEA}. A comparison between \ac{KR} and \ac{PINN} is performed. \ac{DMD} is not considered in this first stage, as it restores dynamical systems from samples of the vector field and not samples of orbits. Instead, \ac{IKEA} is compared to Extended \ac{DMD}. This section targets only small-scale synthetic systems, as it serves only to illustrate the algorithms' behaviour in practice.

\niceparagraph{Error Measuring} The error measurement between the clean vector field $p(\bm{x})$ and the restored one $\hat{p}(\bm{x})$ is formulated as
\begin{equation}\label{eq:errorVF}
E = \frac{\norm{p(\bm{x})-\hat{p}(\bm{x})}}{\norm{p(\bm{x})}},
\end{equation}
and the SNR is inversely related to the error in $[dB]$.



\subsection{Neural Network Configuration}
\subsubsection{Noise Reduction and Generalization}
\begin{figure}[phtb!]
    \centering 
    \includegraphics[trim=100 130 130 140, clip,width=0.8\textwidth,valign = t]{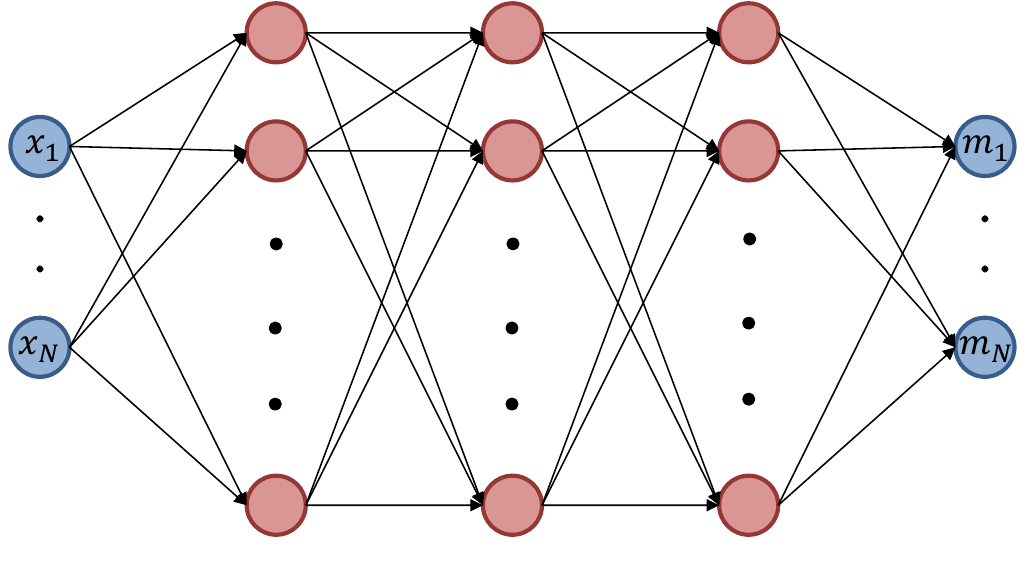}
    \caption{{\bf{NN Configuration}} -- $N$ inputs, $N$ outputs, three hidden layers with $100$ nodes each. Fully connected with $tanh$ as an activation function.}
    \label{fig:NN} 
\end{figure}

The configuration is very basic. For \ac{IKEA} and \ac{KR} with the settings of denoising and generalization, $N$ inputs and $N$ outputs, and three hidden layers with $100$ nodes of a fully connected net with activation function $tanh$, depicted in \Cref{fig:NN}.

\subsubsection{Dimensionality Reduction}
Dimensionality reduction in the context of dynamical systems means reconstructing a vector field with as few \acp{KEF} as possible. Generally, the vector field is a linear combination of the unit velocity measurements’ gradients since Eq. \eqref{eq:dynReconUnitK} holds (point-wise). However, the Jacobian matrix is invertible if there are fewer than $N$ unit functions. Let us assume that there exists a functionally independent set of $K$ \acp{UVM} where $K<N$. Therefore, the objective functional gets the form

\begin{equation}
    \mathcal{O}(\bm{m},\bm{\gamma}) = \frac{1}{2}\sum_{j=1}^L\norm{J(\bm{m}(\bm{x}_j))\bm{p}_j-1}^2+\frac{1}{2}\sum_{j=1}^L\norm{\bm{p}_j-J^T(\bm{m}(\bm{x}_j))\bm{\gamma}(\bm{x}_j)}^2,
\end{equation}
where $\bm{\gamma}(\bm{x}) = [\gamma_1(\bm{x}),\ldots,\gamma_K(\bm{x})]^T$ is a vector of real functions. The first addend forces the set $\{m_i(\bm{x})\}_{i=1}^K$ to solve Eq. \eqref{eq:unitPDE} and the second constrains it to reconstruct the vector field $P(\bm{x})$. These addends yield a parsimonious representation of the vector field. 

As discussed in \Cref{remark:dim}, it is expected that $\bm{\gamma}(\bm{x}) = \left(J(\bm{m})J^T(\bm{m}\right)^{-1} \bm{1}$. However, we exclude this term from the cost function to avoid singularity. Therefore, the coefficients $(\{\gamma_i\}_{i=1}^K)$ is learned independently from the \acp{UVM}. Usually, we assume $N \gg K$; therefore, the benefit of this representation is clear, since training $2K$ functions ($\bm{m}$ and $\bm{\gamma}$) is faster than training $N$ functions.

For the dimensionality reduction use case, the setup is: $N$ inputs, $2K$ outputs, and the same configuration for the hidden layers. $N$ inputs for $N$ coordinates $K$ outputs for $K$ unit velocity measurements, and $K$ outputs for $K$ coefficients $\bm{\gamma}(\bm{x})$.

\subsection{Koopman Regularization}
The following experiments examine robustness to noise, generalization capability, and dimensionality reduction. In addition, a comparison with \ac{PINN} is presented here. The noise reduction and generalization settings are tested on 2D linear and nonlinear dynamical systems, whereas dimensionality reduction is tested on the Lorenz butterfly.

A 2D linear dynamical system either has real, complex, or imaginary eigenvalues. Here, all the options are discussed with the following systems. The linear dynamics are given by
\begin{equation}\label{eq:linA}
    \dot{\bm{x}}=A\bm{x}
\end{equation}
where $A$ in each experiment gets the following values
\begin{equation}\label{eq:lin}
    \frac{1}{200}\begin{bmatrix}
        11&-5\\-5&11
    \end{bmatrix},\,\, \frac{1}{10}\begin{bmatrix}
        -0.4&0.1\\-0.4&-0.5
    \end{bmatrix}, \,\, \frac{1}{10}\begin{bmatrix}
        0&1\\-1&0
    \end{bmatrix}
\end{equation}
where the eigenvalues of these dynamical systems are real, complex, and imaginary, respectively.

The nonlinear system, discussed here, has an unstable equilibrium at $(0,0)$ and a stable limit cycle in the unit circle. The nonlinear dynamical system is given by,
\begin{equation}\label{eq:nonlin}
    \begin{split}
        \dot{x}_1 &= \frac{1}{1000}(-x_2 + x_1(1-x_1^2-x_2^2))\\
        \dot{x}_2 &= \frac{1}{1000}(x_1 + x_2(1-x_1^2-x_2^2))
    \end{split}.
\end{equation}
The factors in Eqs.\eqref{eq:lin} and \eqref{eq:nonlin} are to keep the vector field in the same order of magnitude. The domain of interest is $\mathcal{D}=[6,12]\times [-3,3]$ and the vector field is sampled at a resolution of $dx=0.1$.

Dimensionality reduction is tested on the Lorenz Butterfly. We show that the local behaviour in this setting can be approximated by a 2D dynamical system embedded in 3D.

\subsubsection{Noise Reduction}
\begin{figure}[phtb!]
    \centering 
    \captionsetup[subfigure]{justification=centering}
    \begin{subfigure}[t]{0.495\textwidth} 
\includegraphics[trim=30 30 30 20, clip,width=1\textwidth,valign = t]{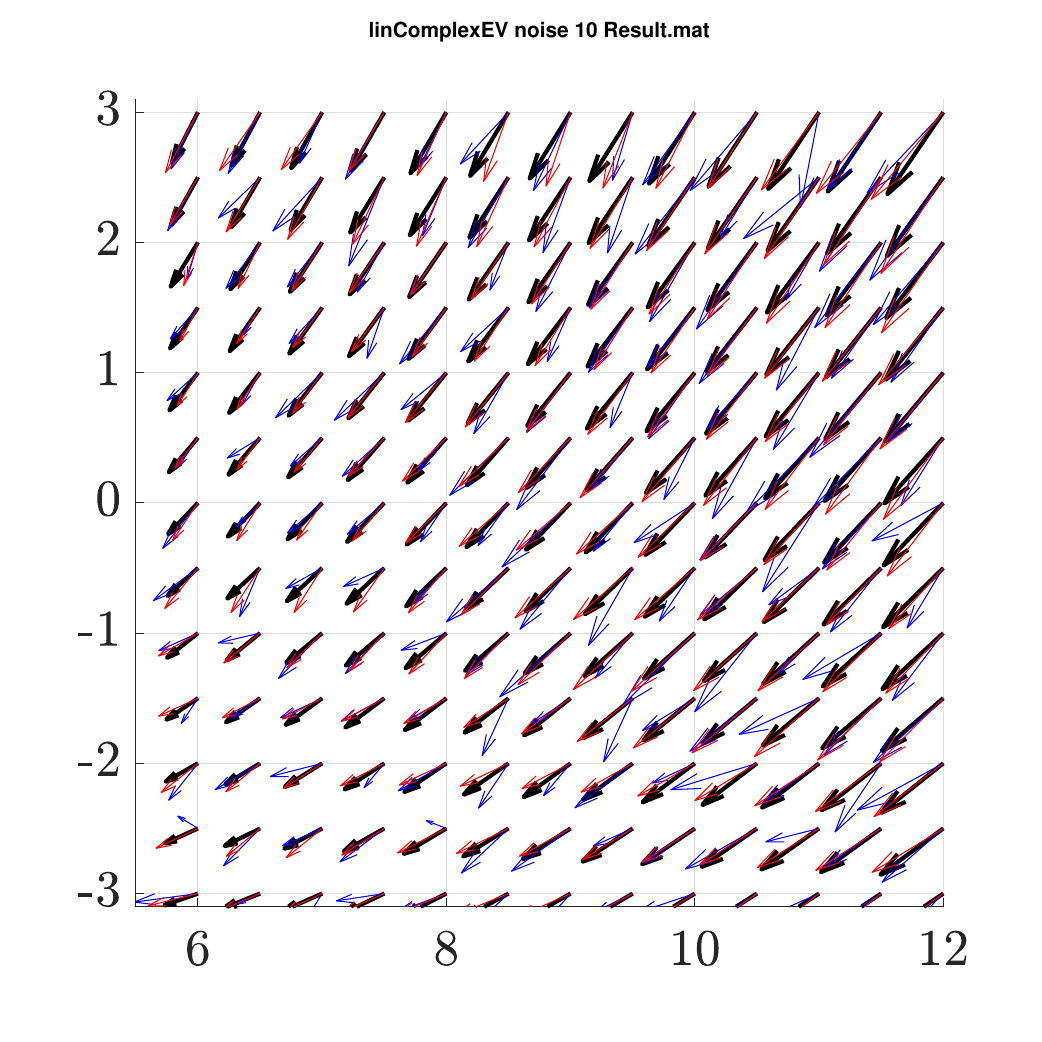}
        \caption{Noisy samples of vector field in $\mathcal{D}$ of Eq. \eqref{eq:linA} where $A=\dfrac{1}{10}\begin{bmatrix}
            -0.4&0.1\\-0.4&-0.5
        \end{bmatrix}$ and their restorations}
        \label{subfig:linComplexEV_noise_10_ResultVF}
    \end{subfigure}
    \begin{subfigure}[t]{0.495\textwidth} 
\includegraphics[trim=30 30 30 20, clip,width=1\textwidth,valign = t]{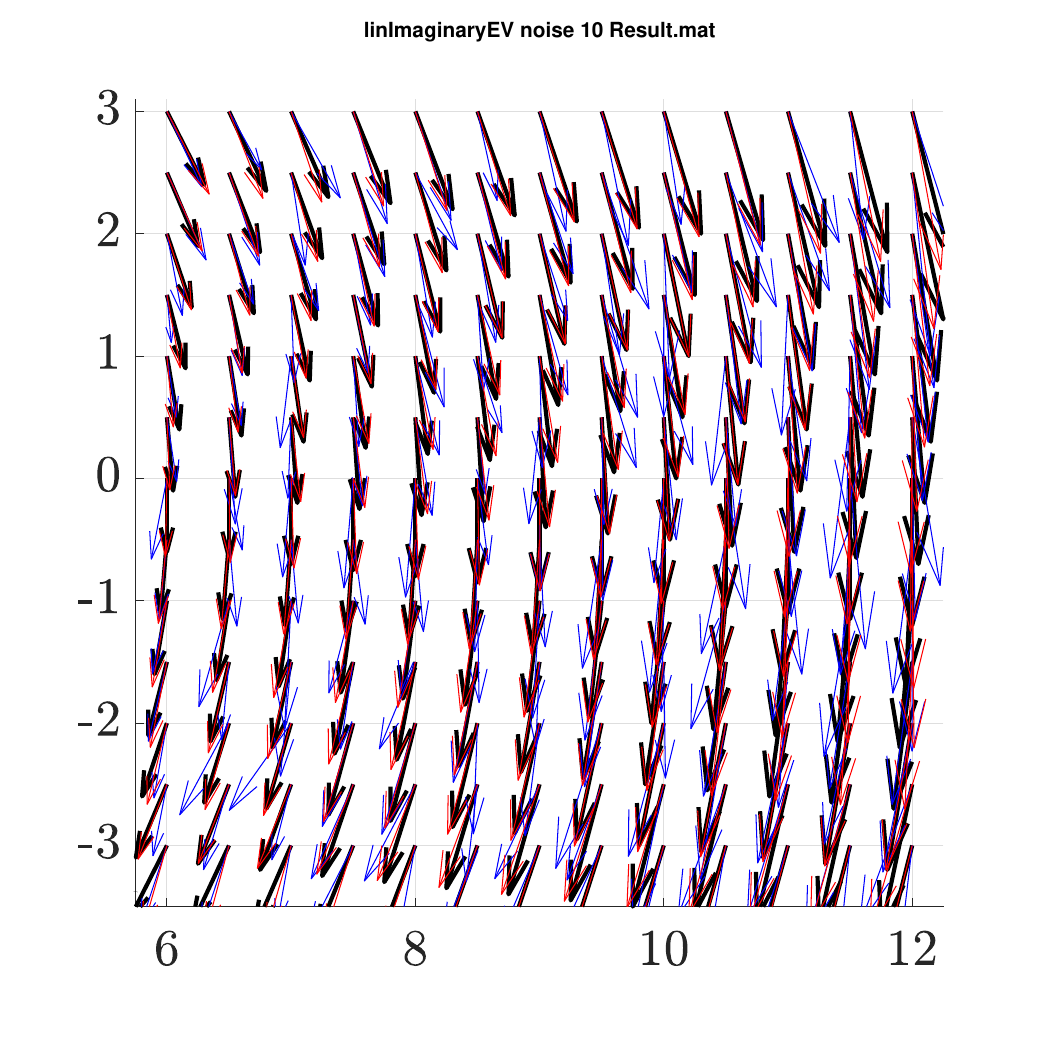}
        \caption{Noisy samples of vector field in $\mathcal{D}$ of Eq. \eqref{eq:linA} where
        $A = \dfrac{1}{10}\begin{bmatrix}
        0&1\\-1&0
    \end{bmatrix}$ and their restorations}
        \label{subfig:linImaginaryEV_noise_10_ResultVF}
    \end{subfigure}\\
    \begin{subfigure}[t]{0.495\textwidth} 
\includegraphics[trim=30 30 30 20, clip,width=1\textwidth,valign = t]{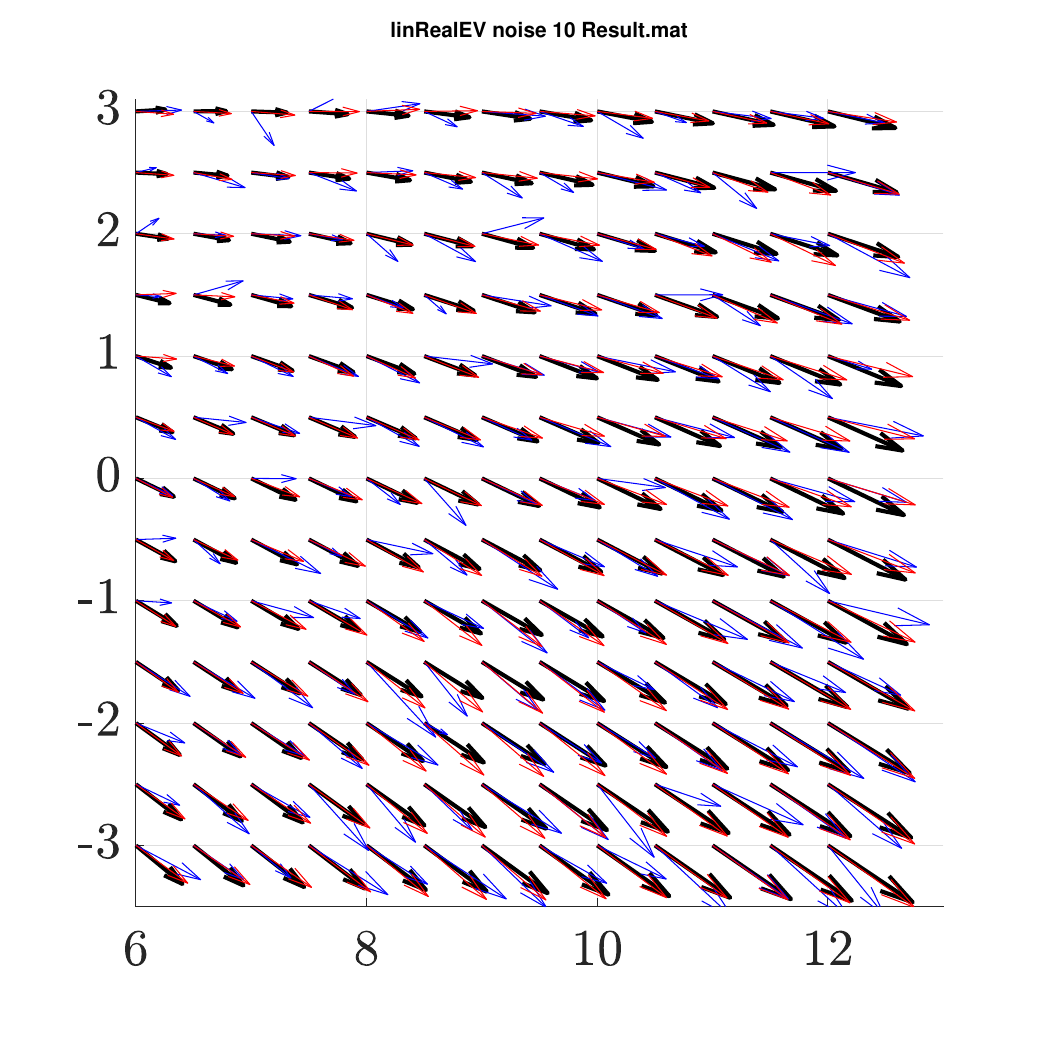}
        \caption{Noisy samples of vector field in $\mathcal{D}$ of Eq. \eqref{eq:linA} where
        $A = \dfrac{1}{200}\begin{bmatrix}
        11&-5\\-5&11
    \end{bmatrix}$ and their restorations}
        \label{subfig:linRealEV_noise_10_ResultVF}
    \end{subfigure}
    \begin{subfigure}[t]{0.495\textwidth} 
\includegraphics[trim=30 30 30 20, clip,width=1\textwidth,valign = t]{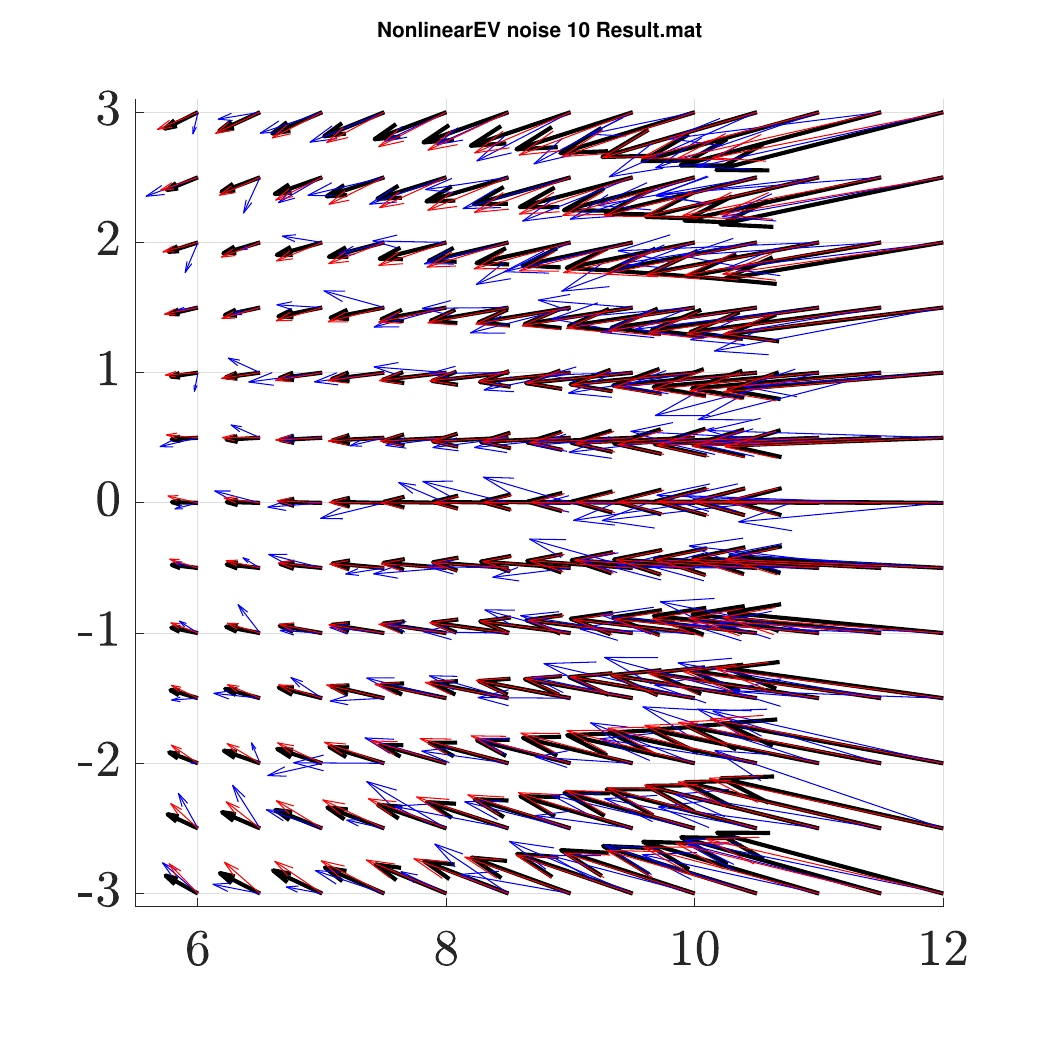}
        \caption{Noisy samples of vector field in $\mathcal{D}$ of Eq. \eqref{eq:nonlin} and their restorations}
        \label{subfig:NonlinearEV_noise_10_ResultVF}
    \end{subfigure}

    \caption{{\bf{Noise Reduction}} - noisy (blue), ground truth (black), and restored (red) vector fields are depicted.}
    \label{fig:Denoising} 
\end{figure}

Here, \emph{Koopman Regularization} is applied to Problem \ref{Pro:Main} with the setting of {\bf{noise reduction}} presented in \Cref{sec:settingNoiseRe}. The vector field sample set is measured from the (linear and nonlinear 2D) dynamical systems in Eqs. \eqref{eq:linA}-\eqref{eq:nonlin}. The domain $\mathcal{D}=[6,12]\times [-3,3]$. This domain is sampled every $dx=0.1$. The noise is normally distributed $\bm{n}_i\sim\mathcal{N}(0,0.1)$. The results are depicted in \Cref{fig:Denoising}. The blue arrows represent the noisy samples, the black arrows are the clean vector fields, and the red arrows are the restored vector fields. 

\begin{figure}[phtb!]
    \centering 
    \captionsetup[subfigure]{justification=centering}
    \begin{subfigure}[t]{0.495\textwidth} 
\includegraphics[trim=55 10 50 20, clip,width=1\textwidth,valign = t]{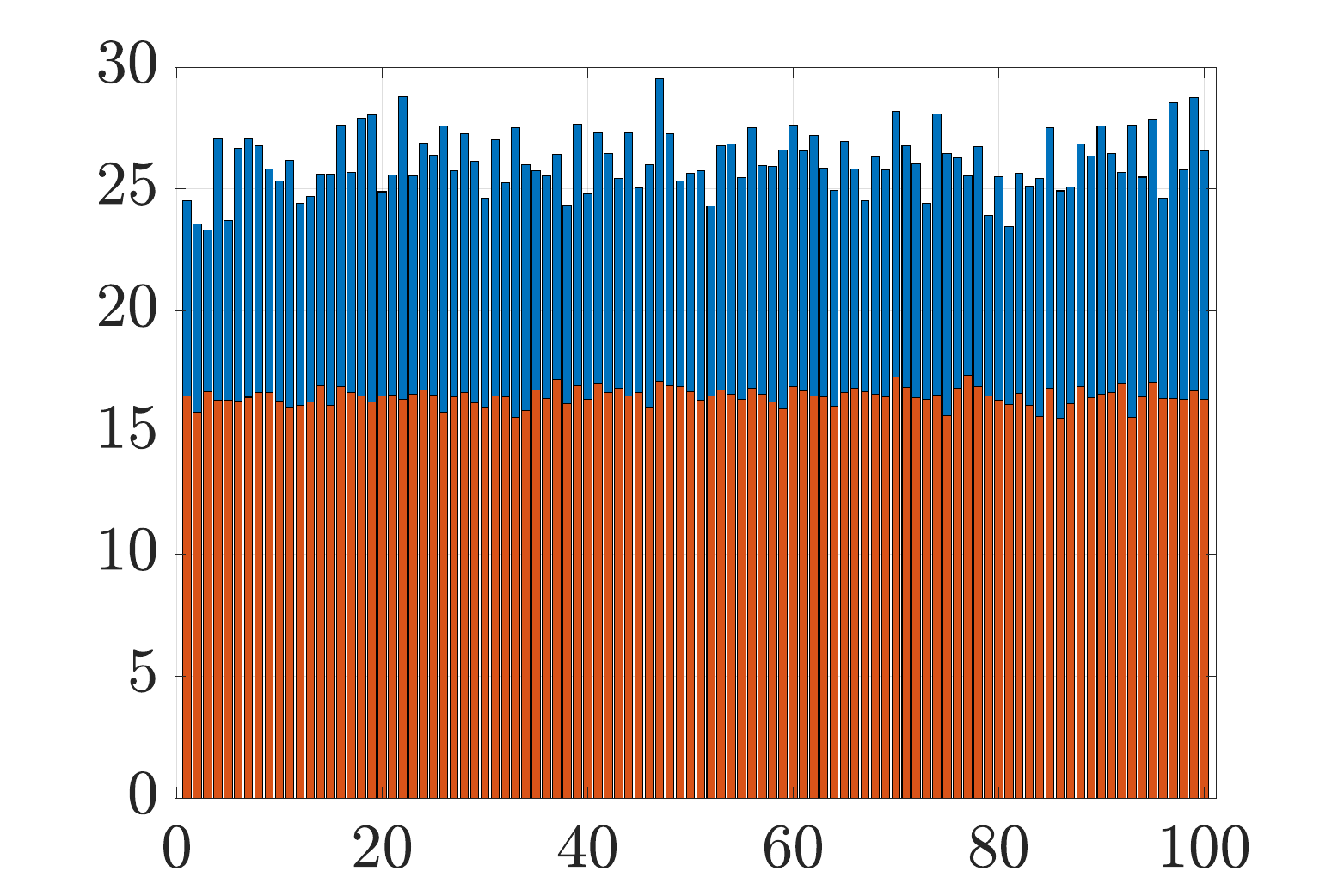}
\caption{{\bf{SNR Improvement.}} Results from experiments on linear 2D dynamic Eq. \eqref{eq:linA} with imaginary eigenvalues.}
        \label{subfig:ImDist}
    \end{subfigure}
    \begin{subfigure}[t]{0.495\textwidth} 
\includegraphics[trim=55 10 50 20, clip,width=1\textwidth,valign = t]{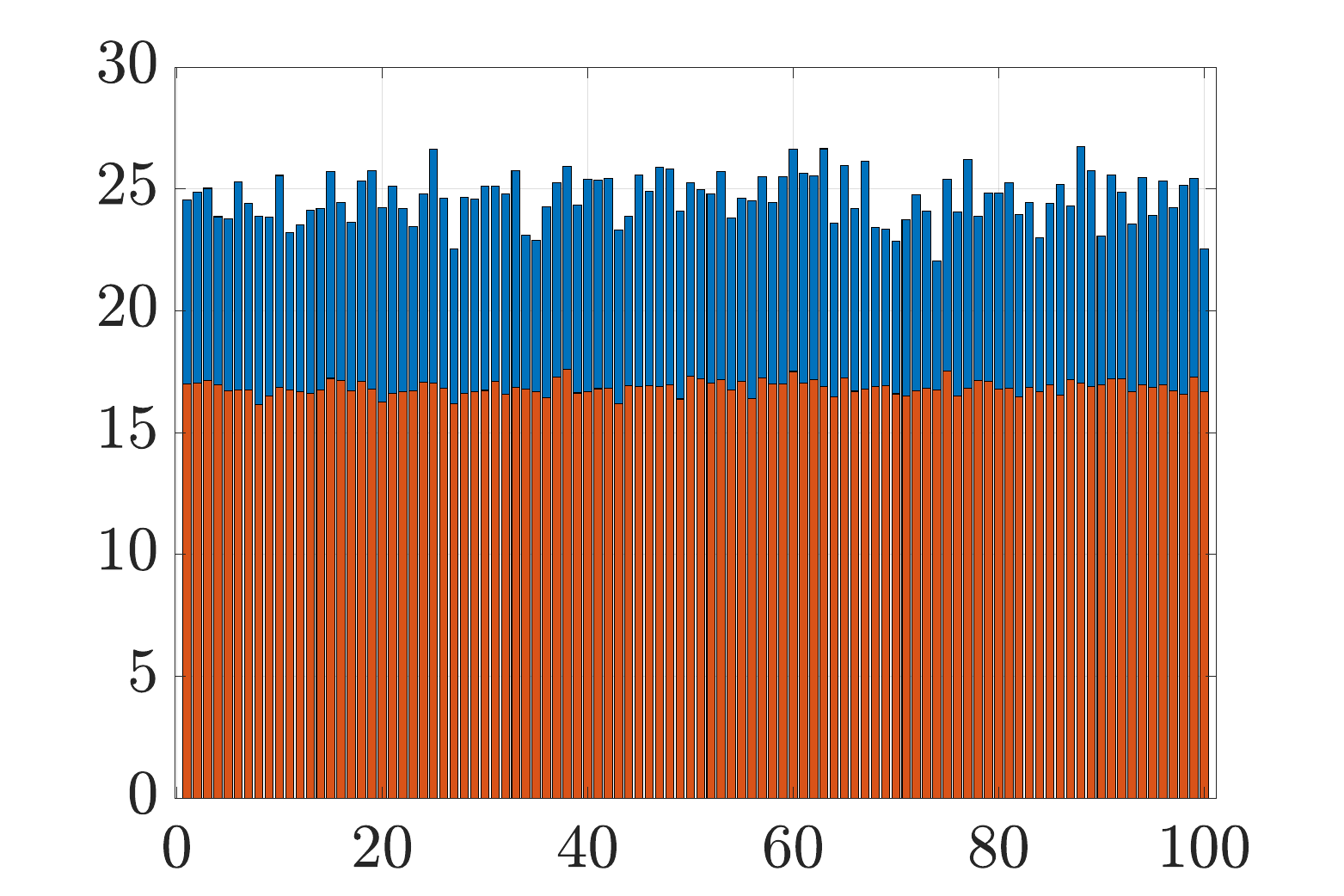}
        \caption{{\bf{SNR Improvement.}} Results from experiments on 2D nonlinear dynamic Eq. \eqref{eq:nonlin}.}
        \label{subfig:NonDist}
    \end{subfigure}
    \caption{{\bf{Noise Reduction Quality -- SNR Improvement.}} $Y$ axis -- SNR level in $[dB]$. $X$ axis -- number of experiment. Red Bars -- SNR level in $[dB]$ before \emph{Koopman Regularization}. Blue Bars -- SNR level in $[dB]$ after  \emph{Koopman Regularization}.}
    \label{Fig:qualityBatch} 
\end{figure}

To check the efficiency of noise reduction in Koopman Regularization, the aforementioned experiment was conducted $100$ times for linear Eq. \eqref{eq:linA} with $A = \dfrac{1}{10}\begin{bmatrix}
    0&1\\
    -1&0
\end{bmatrix}$ and $100$ times for the nonlinear dynamical systems in Eq. \eqref{eq:nonlin}. The results are summarized in \Cref{Fig:qualityBatch}. The results of the linear and nonlinear systems are summarized in \Cref{subfig:ImDist} and \Cref{subfig:NonDist}, respectively. Each experiment is represented by two bars (red and blue) in these figures. The red one is the SNR level of the noisy samples, and the blue one is the level of SNR after \emph{Koopman Regularization}. On average, the SNR improvement is $9.6dB$ in the linear system (\Cref{subfig:ImDist}) and $7.8dB$ in the nonlinear system  (\Cref{subfig:NonDist}).

\begin{figure}[phtb!]
    \centering 
    \captionsetup[subfigure]{justification=centering}
    \begin{subfigure}[t]{0.32\textwidth} 
\includegraphics[trim=0 0 0 20, clip,width=1\textwidth,valign = t]{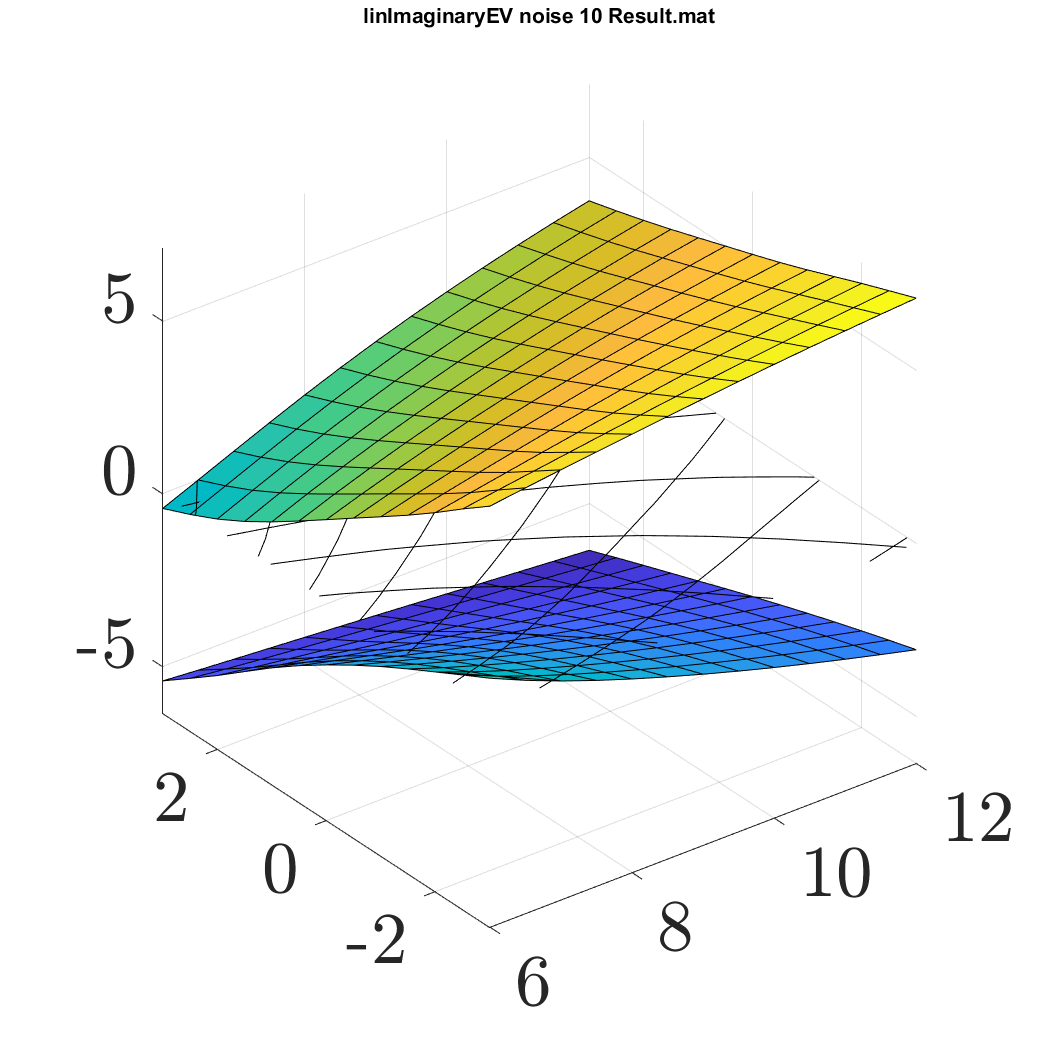}
        \label{subfig:linImaginaryEV_noise_10_Resultmanifold}
    \end{subfigure}
    \begin{subfigure}[t]{0.32\textwidth} 
\includegraphics[trim=0 0 0 20, clip,width=1\textwidth,valign = t]{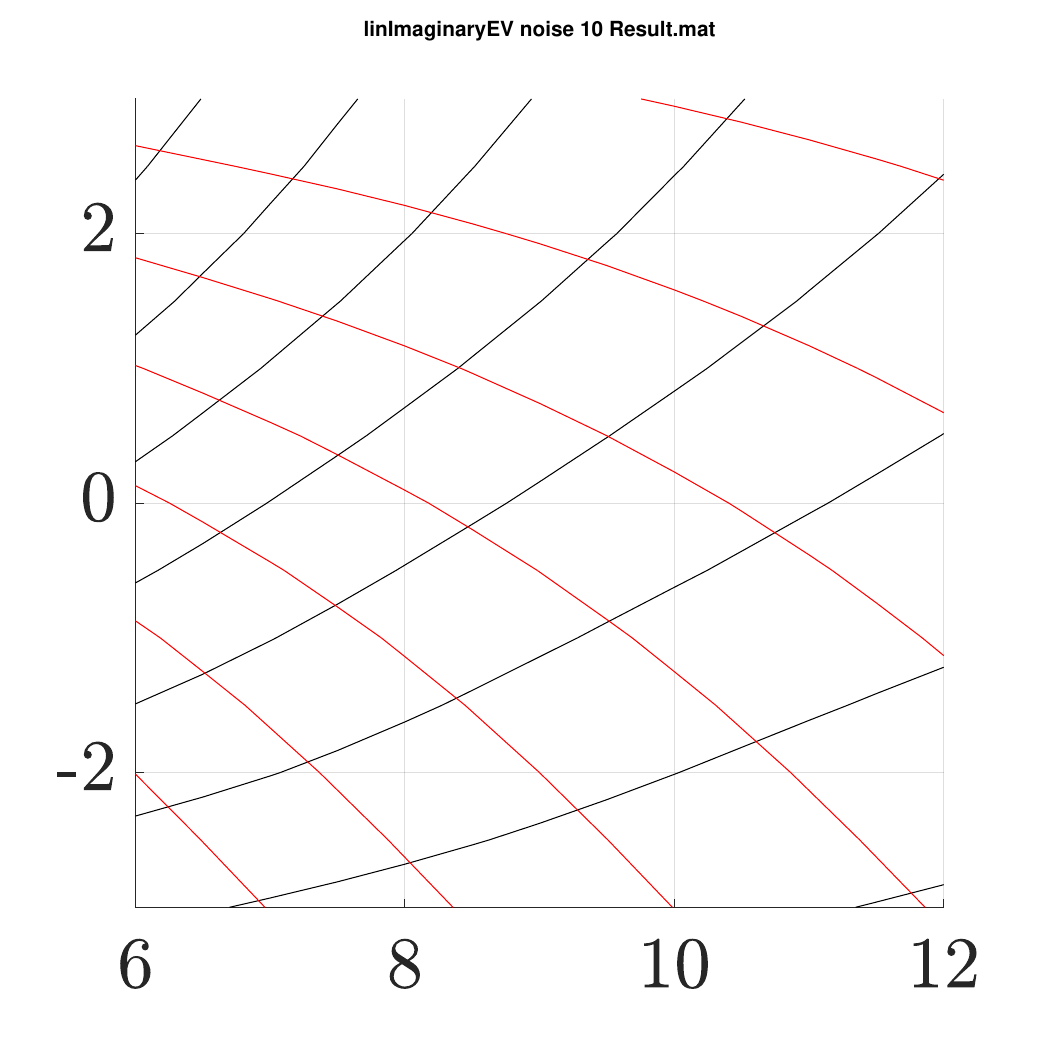}
        \label{subfig:linImaginaryEV_noise_10_Resultcontour}
    \end{subfigure}
    \begin{subfigure}[t]{0.32\textwidth} 
\includegraphics[trim=0 0 0 15, clip,width=1\textwidth,valign = t]{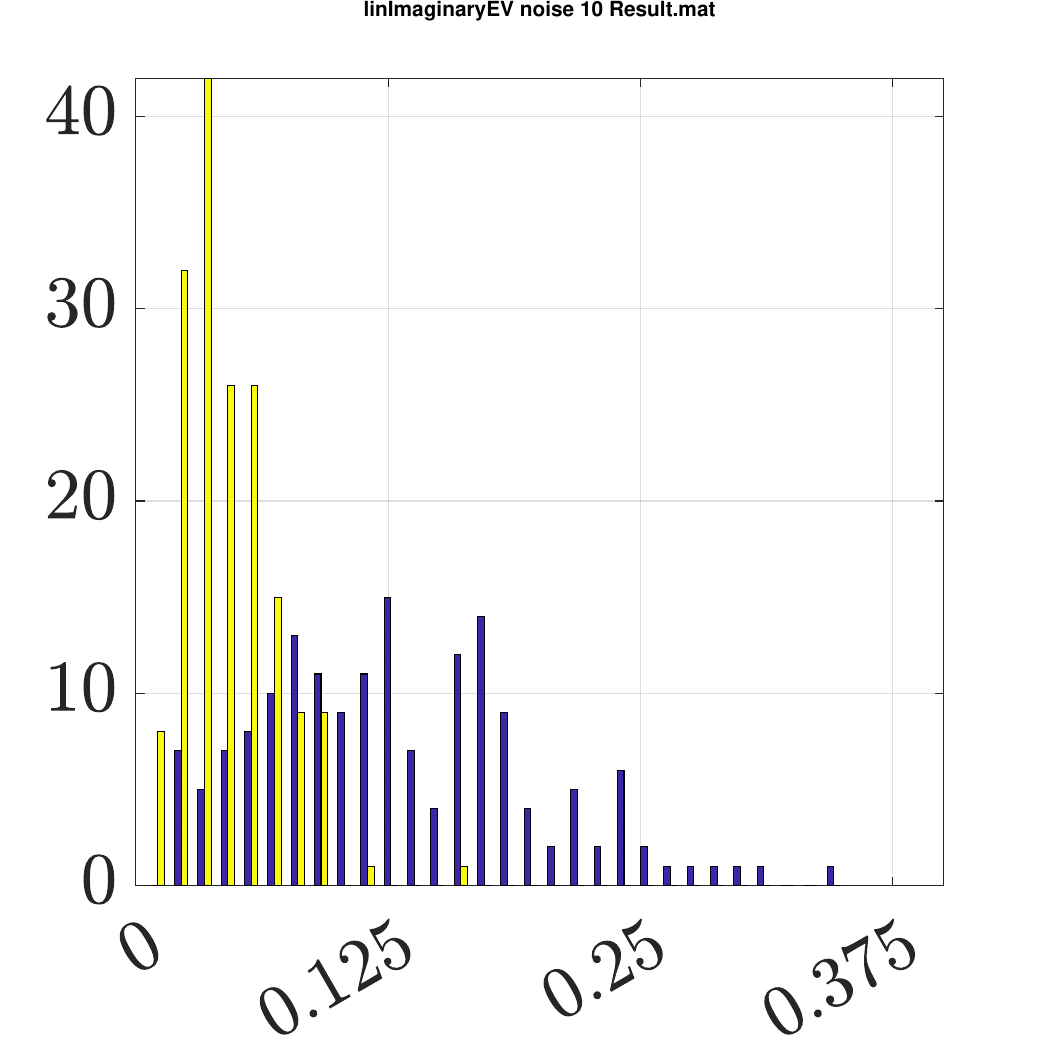}
        \label{subfig:linImaginaryEV_noise_10_Resultdist}
    \end{subfigure}
    \caption{{\bf{Noise Reduction Quality}} -- From left to right, graphs of $m_1(\bm{x})$ and $m_2(\bm{x})$ are given, the contours of these manifolds are depicted, and the histograms of the noise before and after the \emph{Koopman Regularization}.}
    \label{fig:quality} 
\end{figure}

The SNR improvement is $60\%$ or above. The SNR improvement in the linear system with imaginary eigenvalues is greater than $81\%$. The results of a specific case out of $100$ are shown in \Cref{fig:quality}. 

\subsubsection{Generalization}
\begin{figure}[phtb!]
    \centering 
    \captionsetup[subfigure]{justification=centering}
    \begin{subfigure}[t]{0.495\textwidth} 
\includegraphics[trim=30 30 30 20, clip,width=1\textwidth,valign = t]{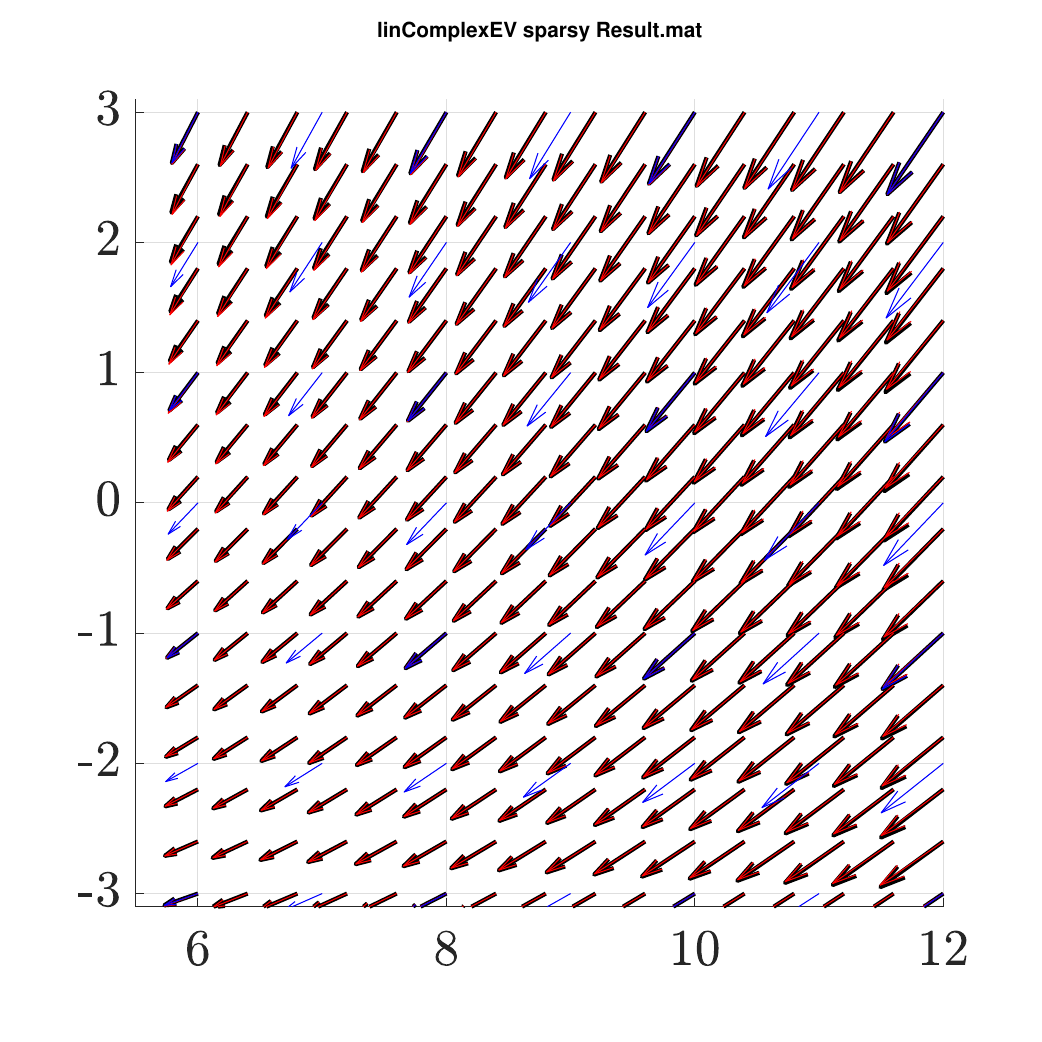}
\caption{Sparse samples of vector field in $\mathcal{D}$ of Eq. \eqref{eq:linA} where $A=\dfrac{1}{10}\begin{bmatrix}
            -0.4&0.1\\-0.4&-0.5
        \end{bmatrix}$ and their generalization}
        \label{subfig:linComplexEV_sparsy_ResultVF}
    \end{subfigure}
    \begin{subfigure}[t]{0.495\textwidth} 
\includegraphics[trim=30 30 30 20, clip,width=1\textwidth,valign = t]{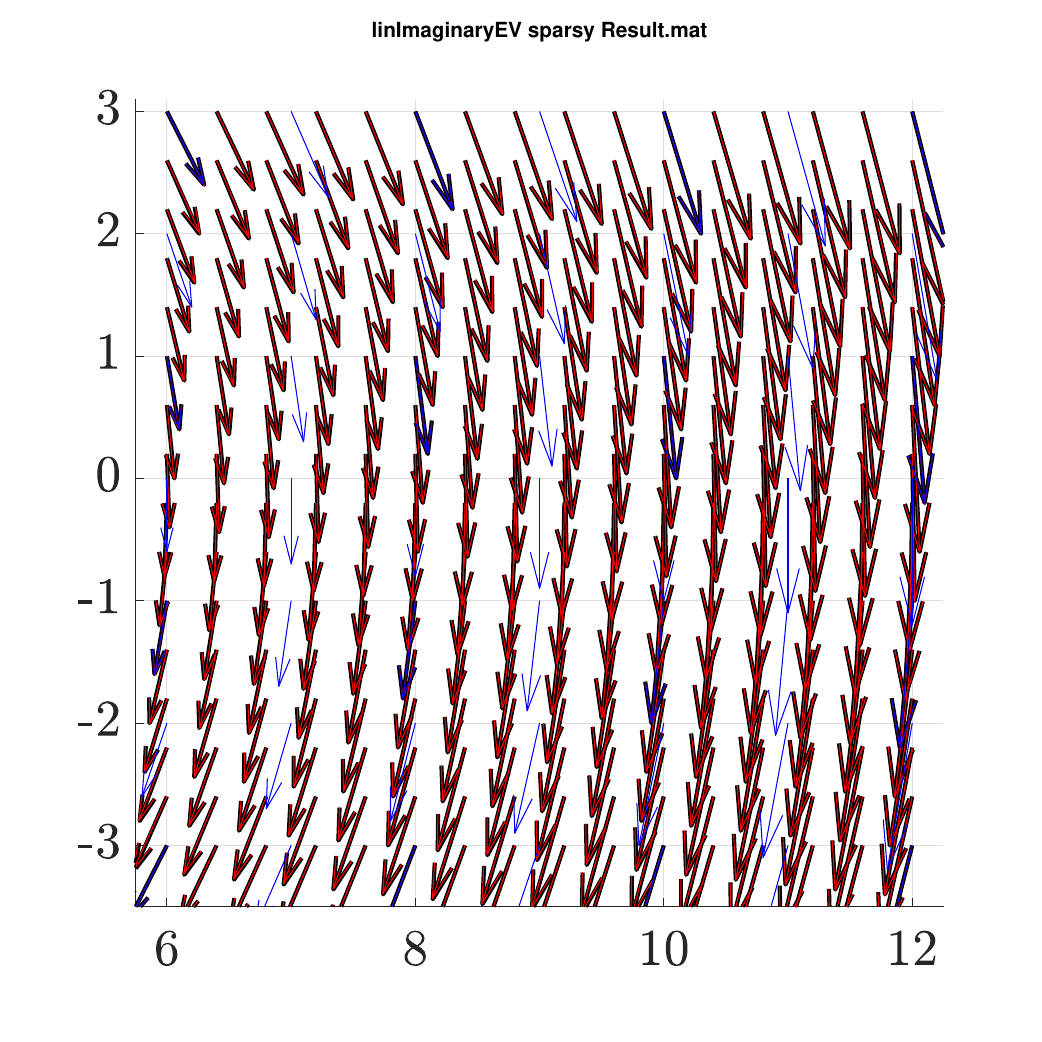}
\caption{Sparse samples of vector field in $\mathcal{D}$ of Eq. \eqref{eq:linA} where $A=\dfrac{1}{10}\begin{bmatrix}
            0&1\\-1&0
        \end{bmatrix}$ and their generalization}
        \label{subfig:linImaginaryEV_sparsy_ResultVF}
    \end{subfigure}\\
    \begin{subfigure}[t]{0.495\textwidth} 
\includegraphics[trim=30 30 30 20, clip,width=1\textwidth,valign = t]{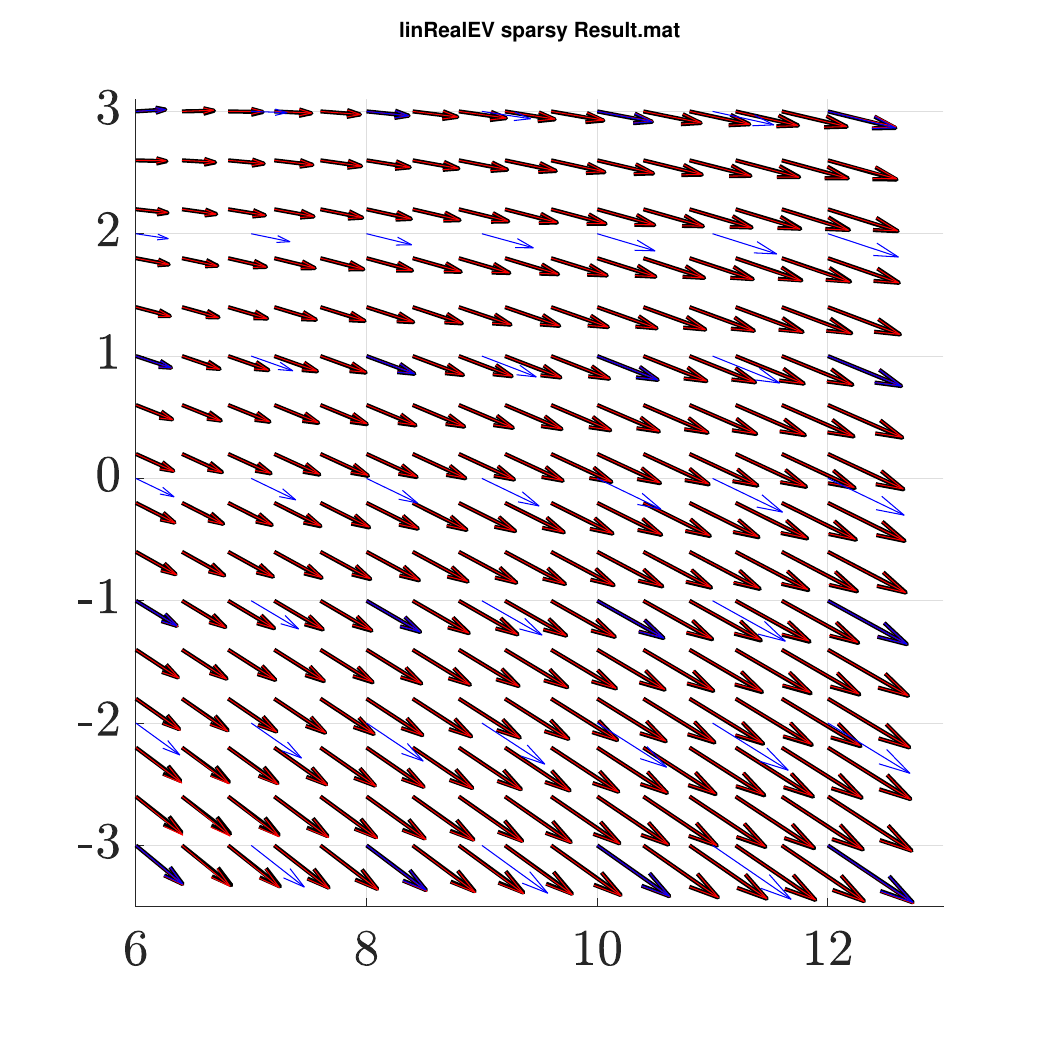}
\caption{Sparse samples of vector field in $\mathcal{D}$ of Eq. \eqref{eq:linA} where $A=\dfrac{1}{200}\begin{bmatrix}
            11&-5\\-5&11
        \end{bmatrix}$ and their generalization}
        \label{subfig:linRealEV_sparsy_ResultVF}
    \end{subfigure}
    \begin{subfigure}[t]{0.495\textwidth} 
\includegraphics[trim=30 30 30 20, clip,width=1\textwidth,valign = t]{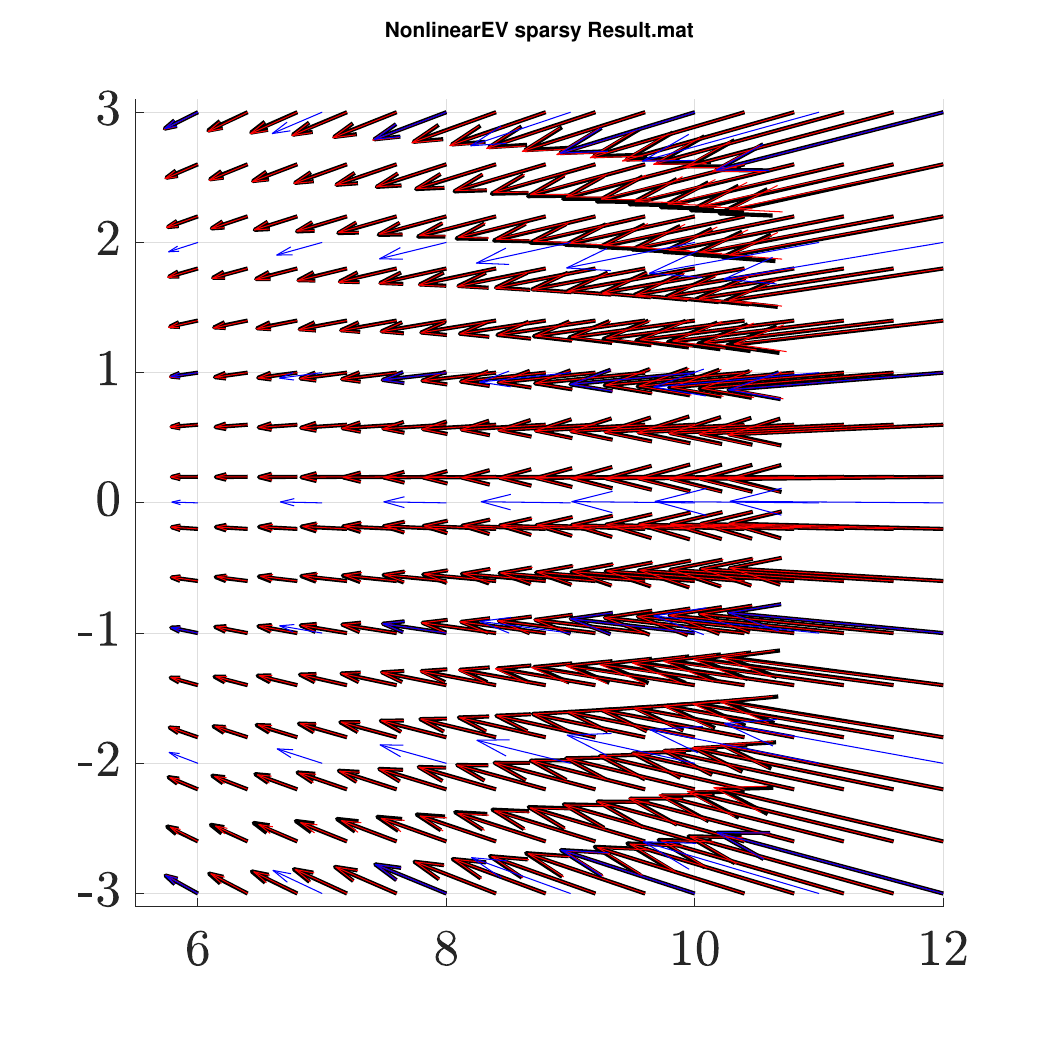}
\caption{Sparse samples of vector field in $\mathcal{D}$ of Eq. \eqref{eq:nonlin} and their generalization}
        \label{subfig:NonlinearEV_sparsy_ResultVF}
    \end{subfigure}

    \caption{{\bf{Generalization}} - sparse samples (blue), ground truth (black), generalized (red) of the vector fields are presented.}
    \label{fig:Generalization} 
\end{figure}

Here, the vector fields in Eqs. \eqref{eq:lin} and \eqref{eq:nonlin} are sampled every $dx=1$ in the domain $\mathcal{D}=[6,12]\times[-3,3]$. The generalization results from \Cref{alg:KR} are summarized in \Cref{fig:Generalization}, where the blue arrows indicate the sparse samples, the black arrows indicate the ground truth, and the red arrows indicate the generalized vector field. The generalization process yields accurate results. The error reconstruction (Eq. \eqref{eq:errorVF}) in the nonlinear dynamics is $3.01\%$, in the linear cases $8.45\%$, $0.25\%$, and $0.6\%$ for complex, imaginary, and real eigenvalues, respectively.

\subsubsection{Dimensionality Reduction}
\begin{figure}[phtb!]
    \centering 
    \captionsetup[subfigure]{justification=centering}
    \begin{subfigure}[t]{0.495\textwidth} 
\includegraphics[trim=0 0 0 20, clip,width=1\textwidth,valign = t]{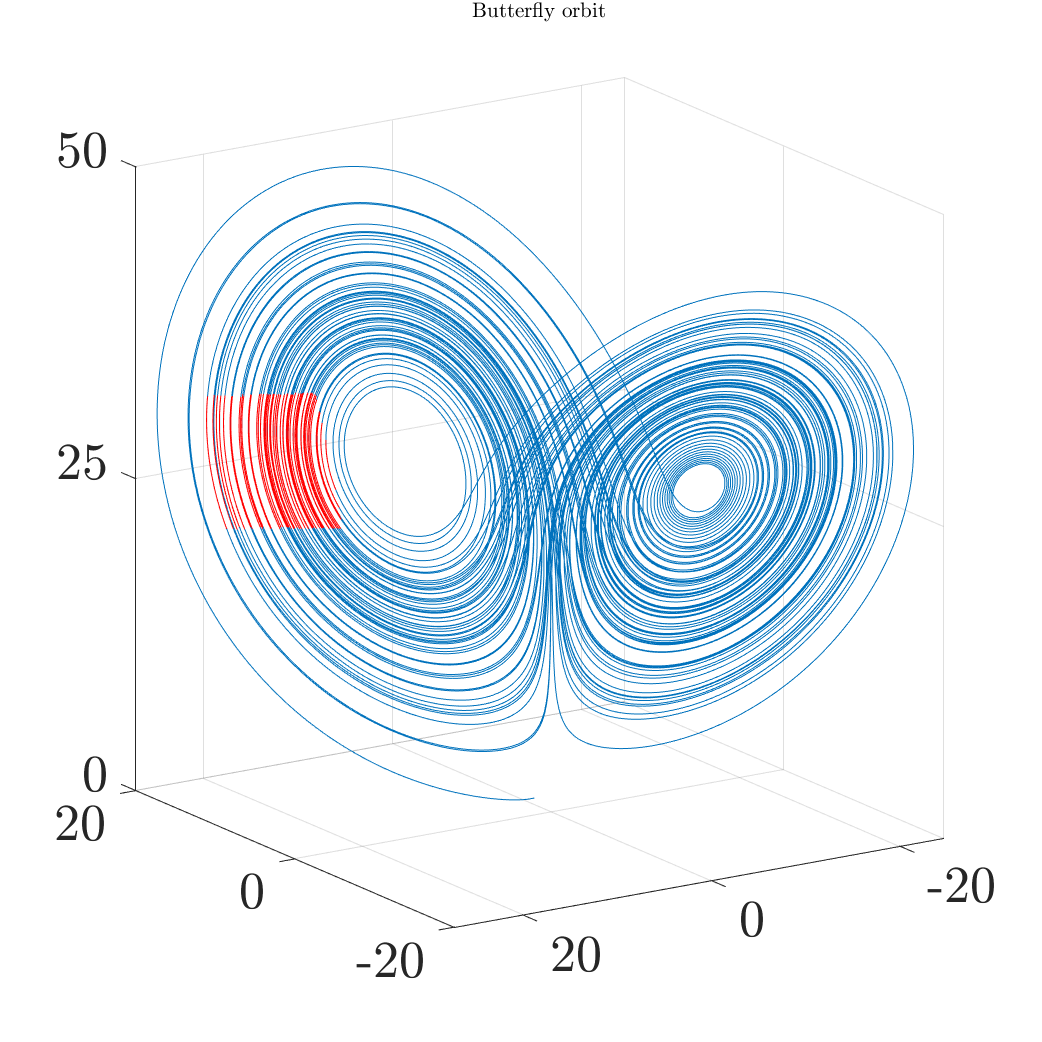}
        \label{subfig:ButterflyOrbit}
    \end{subfigure}
    \begin{subfigure}[t]{0.495\textwidth} 
\includegraphics[trim=0 0 0 20, clip,width=1\textwidth,valign = t]{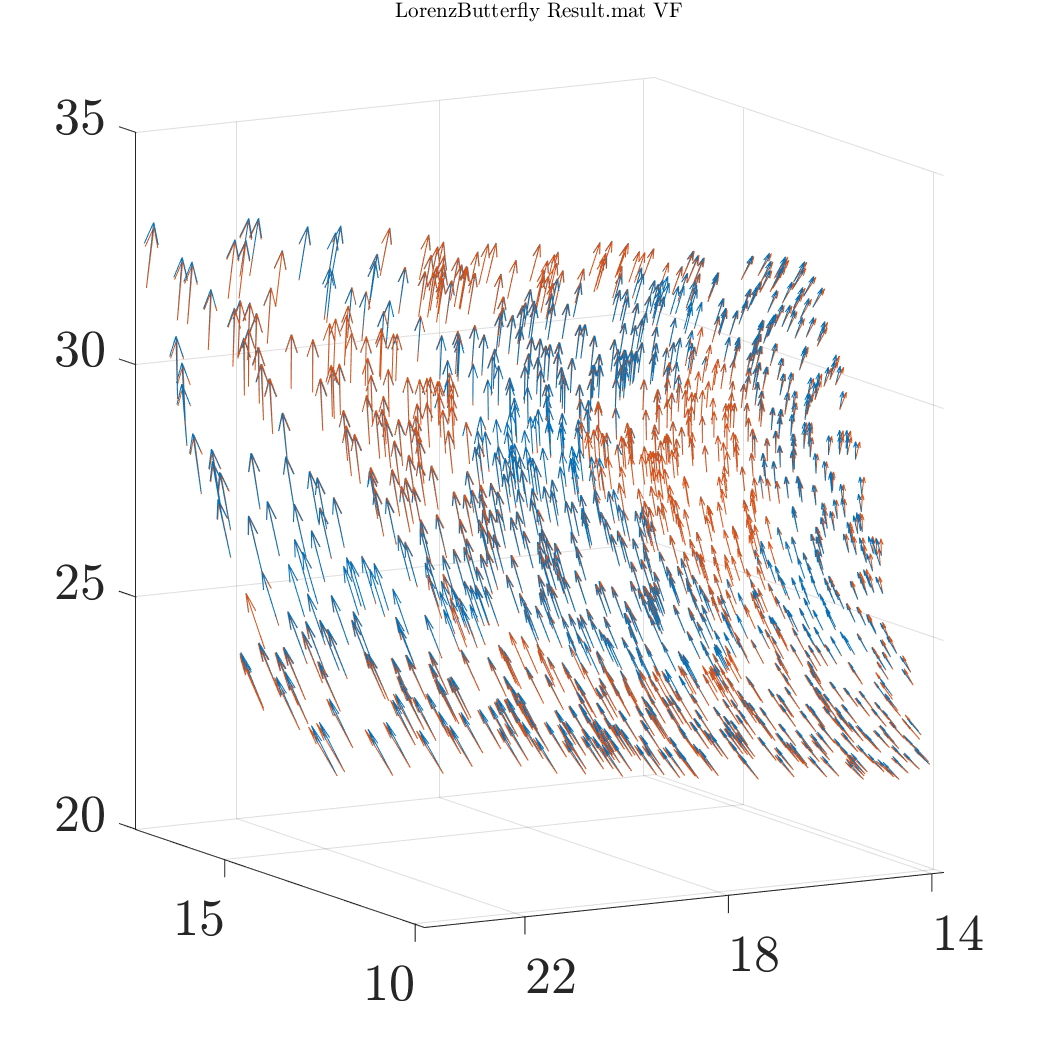}
        \label{subfig:LorenzButterfly_ResultVF}
    \end{subfigure}
    \caption{{\bf{Dimensionality Reduction}} -- left, Lorenz butterfly, the red part approximated as a two dimensional dynamical system; right, vector field restoration, the ground truth (blue) and the restored (red)}
    \label{fig:Butter} 
\end{figure}

For this use case, the Lorenz butterfly is studied as a toy example. Given the dynamical system
\begin{equation}\label{eq:butterfly}
    \begin{split}
        \dot{x}&=-\sigma x + \sigma y\\
        \dot{y}&=\rho x - y - xz\\
        \dot{z}&=-\beta z + xy
    \end{split}    
\end{equation}
where $\sigma = 10,\,\beta = 8/3,\,\rho = 28$ and $(x,y,z)^T(0) = (1,1,1)^T$, the dynamics follow the trajectory depicted in \Cref{fig:Butter} (top left). The study is limited to the region highlighted in red. One can see that the dynamic lies approximately on a plane in this part. Thus, in this neighborhood, the vector field can be restored with only two Koopman eigenfunctions; i.e., via the optimization problem in Eq. \eqref{eq:UVMDR} with $K=2$. \Cref{fig:Butter} depicts the vector field recovered using only two \acp{UVM}. The error from the vector field is $2.5\%$ (Eq. \eqref{eq:errorVF}). 

\subsubsection{Geometric Condition}
\Cref{subsec:cond} thoroughly discusses the geometric condition of a functionally independent set. This discussion yields the definition of the feasible region for a legitimate solution. This condition ensures that the Jacobian is a full-rank matrix, which is essential to system reconstructions (Eqs. \eqref{eq:dynReconUnit} and \eqref{eq:dynReconUnitK}) --- hence, its importance.

\begin{figure}[phtb!]
\includegraphics[trim=0 0 0 20, clip,width=0.85\textwidth,valign = t]{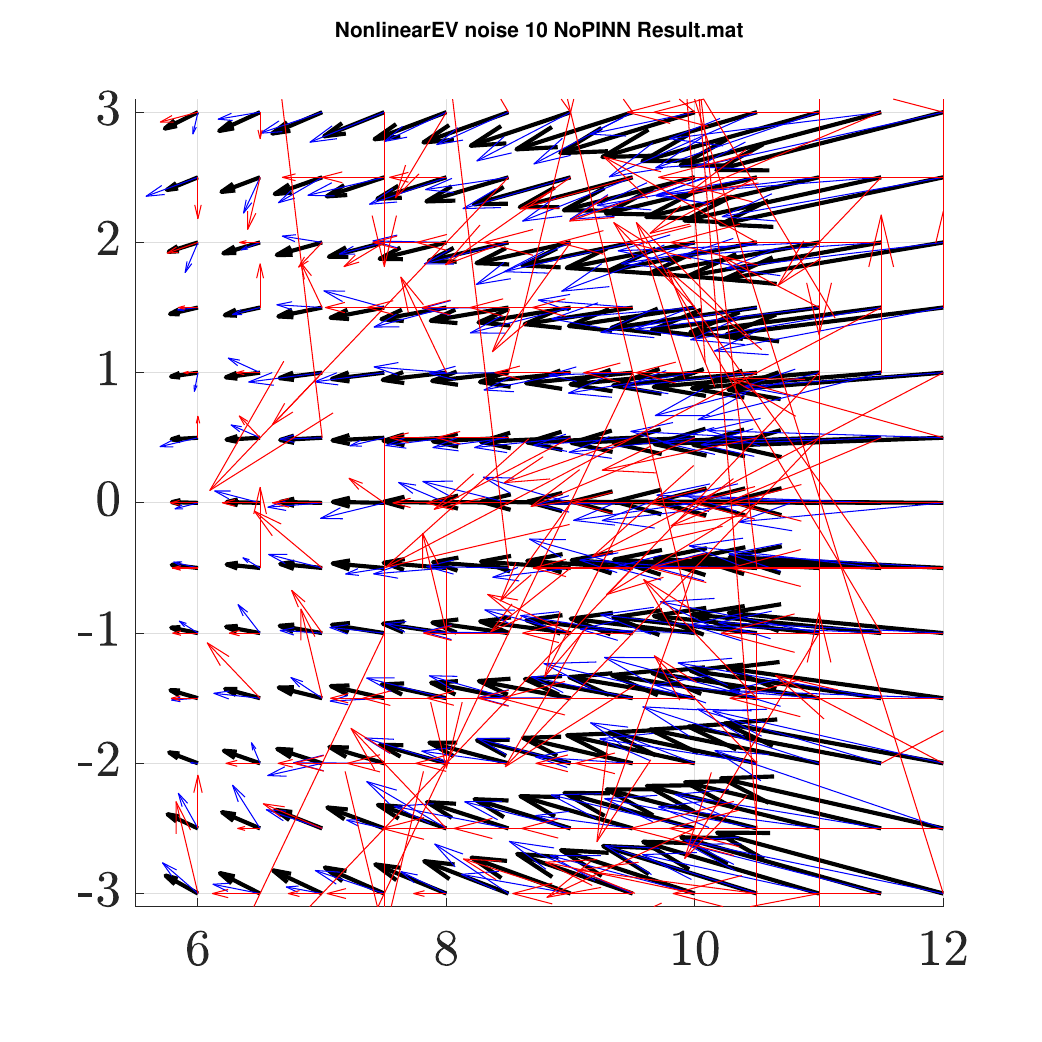}
    \caption{{\bf{PINN}} -- noisy (blue), ground truth (black), and restored (red) vector fields are depicted.}
    \label{fig:NonlinearEV_noise_10_NoPINN_ResultVF} 
\end{figure}

Omitting this constraint turns \ac{KR} into a simple minimization process of the functional $\mathcal{L}(\bm{m})$, meaning finding solutions to \eqref{eq:unitPDE}. Thus, \ac{KR} becomes a supervised learning algorithm as \ac{PINN} approach, since the \acp{KEF} includes the knowledge about the conservative laws of the system (see \Cref{theo:Gery}). This process results in system reconstructions struggling with singularity. We repeat the experiment from \Cref{subfig:NonlinearEV_noise_10_ResultVF} without the functional independence constraint and illustrate the result in \Cref{fig:NonlinearEV_noise_10_NoPINN_ResultVF}.
The result clearly shows inconsistent behavior of the reconstruction, typical of singularity.

\subsection{Independent Koopman Eigenfunction Approximation}
\subsubsection{Settings}
\begin{figure}[phtb!]
    \centering 
    \includegraphics[trim=0 30 0 20, clip,width=0.85\textwidth,valign = t]{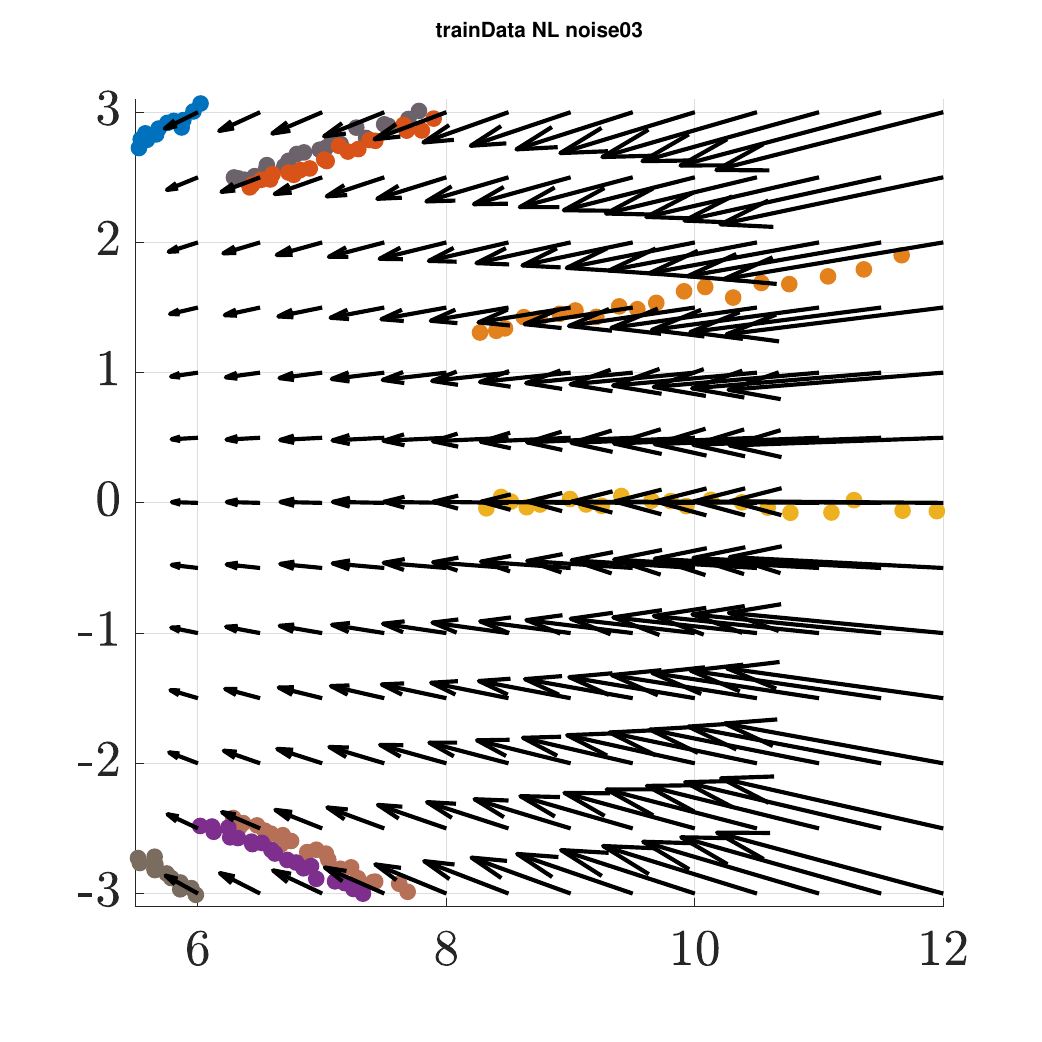}
    \caption{{\bf{Problem Settings}} -- Restore the vector field (black) using the eight sampled orbits.}
    \label{fig:OurProblems} 
\end{figure}

In this part of the experiments, we focus on the 2D nonlinear dynamical system (Eq. \eqref{eq:nonlin}) in the domain $\mathcal{D}$ as defined above with the noise-reduction setting (Eq. \eqref{eq:ikeaNoise}). We draw eight distinct initial conditions in $\mathcal{D}$ and sample each orbit $20$ times with a time step $dt = 0.2$ (see \Cref{fig:OurProblems}). 
We invoke this algorithm with clean data and with additional white Gaussian noise where $\textrm{SNR}=55.1dB,\,49.1dB,\,45.6dB$. With these levels of $\textrm{SNR}$, \ac{DMD} and its relative should work accurately enough. 


\subsubsection{\texorpdfstring{\acl{IKEA}}{}}

\begin{figure}[phtb!]
    \centering 
    \captionsetup[subfigure]{justification=centering}
    \begin{subfigure}[t]{0.495\textwidth} 
\includegraphics[trim=30 30 30 20, clip,width=1\textwidth,valign = t]{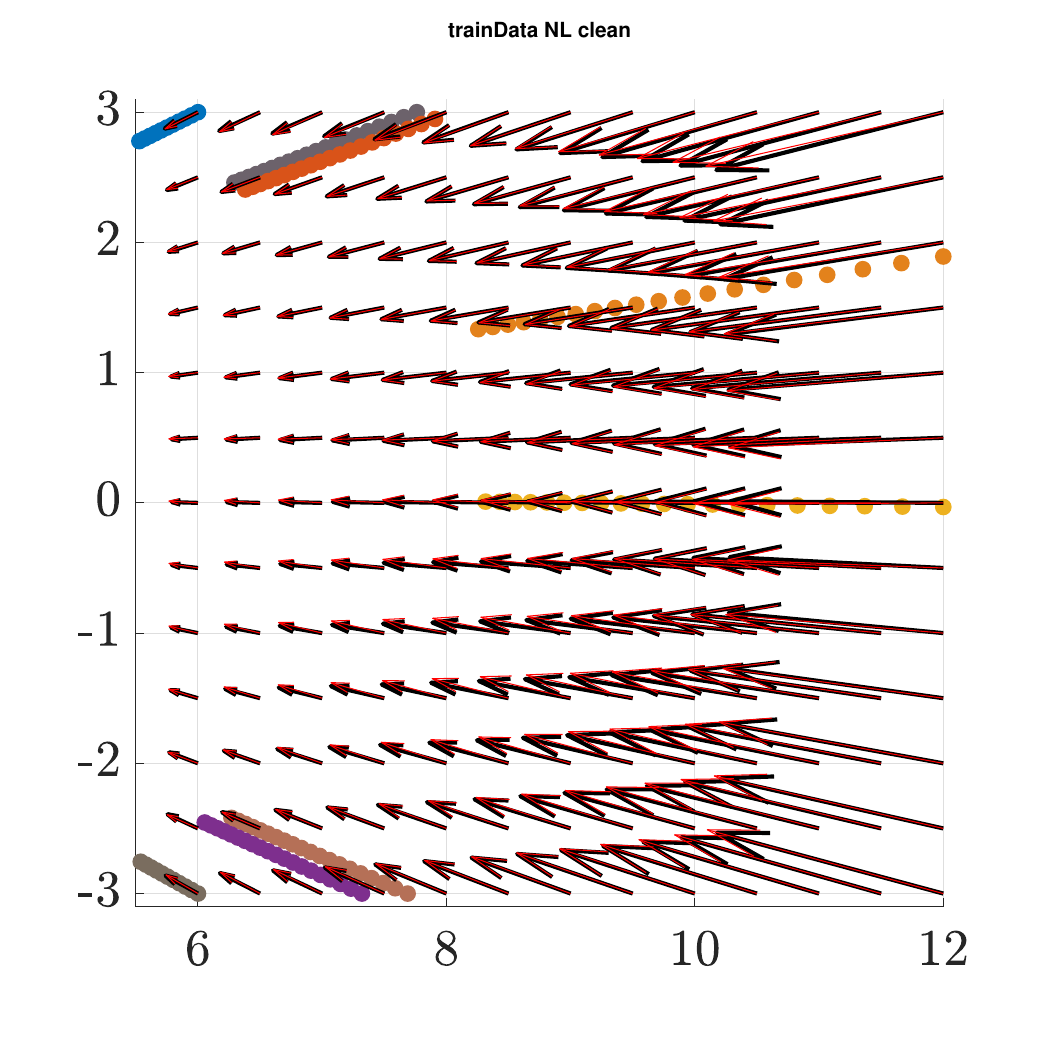}
\caption{{\bf IKEA} -- eight orbits without noise}
        \label{subfig:trainData_NL_clean_IKEA}
    \end{subfigure}
    \begin{subfigure}[t]{0.495\textwidth} 
\includegraphics[trim=30 30 30 20, clip,width=1\textwidth,valign = t]{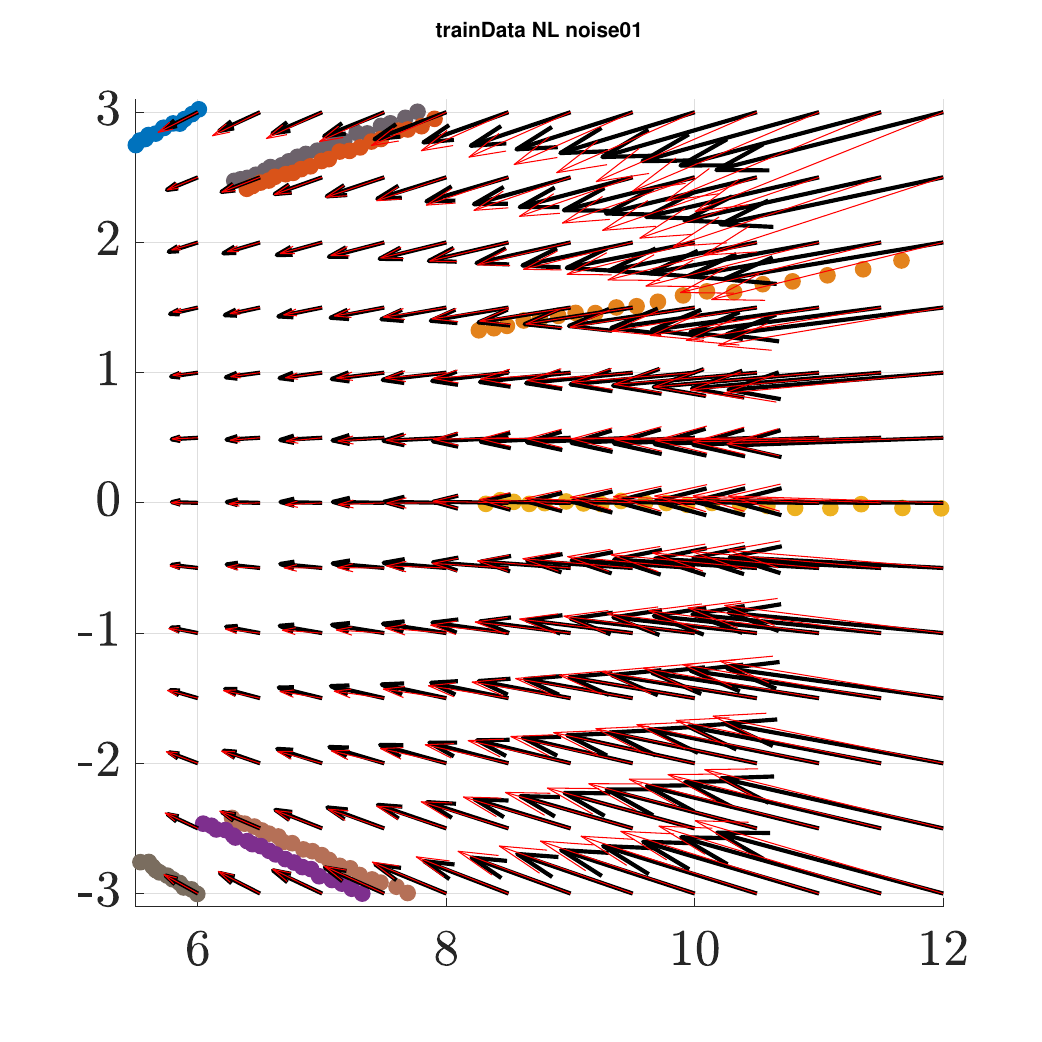}
\caption{{\bf IKEA} -- eight orbits with $\textrm{SNR}=55.1dB$}
        \label{subfig:trainData_NL_noise01_IKEA}
    \end{subfigure}\\
    \begin{subfigure}[t]{0.495\textwidth} 
\includegraphics[trim=30 30 30 20, clip,width=1\textwidth,valign = t]{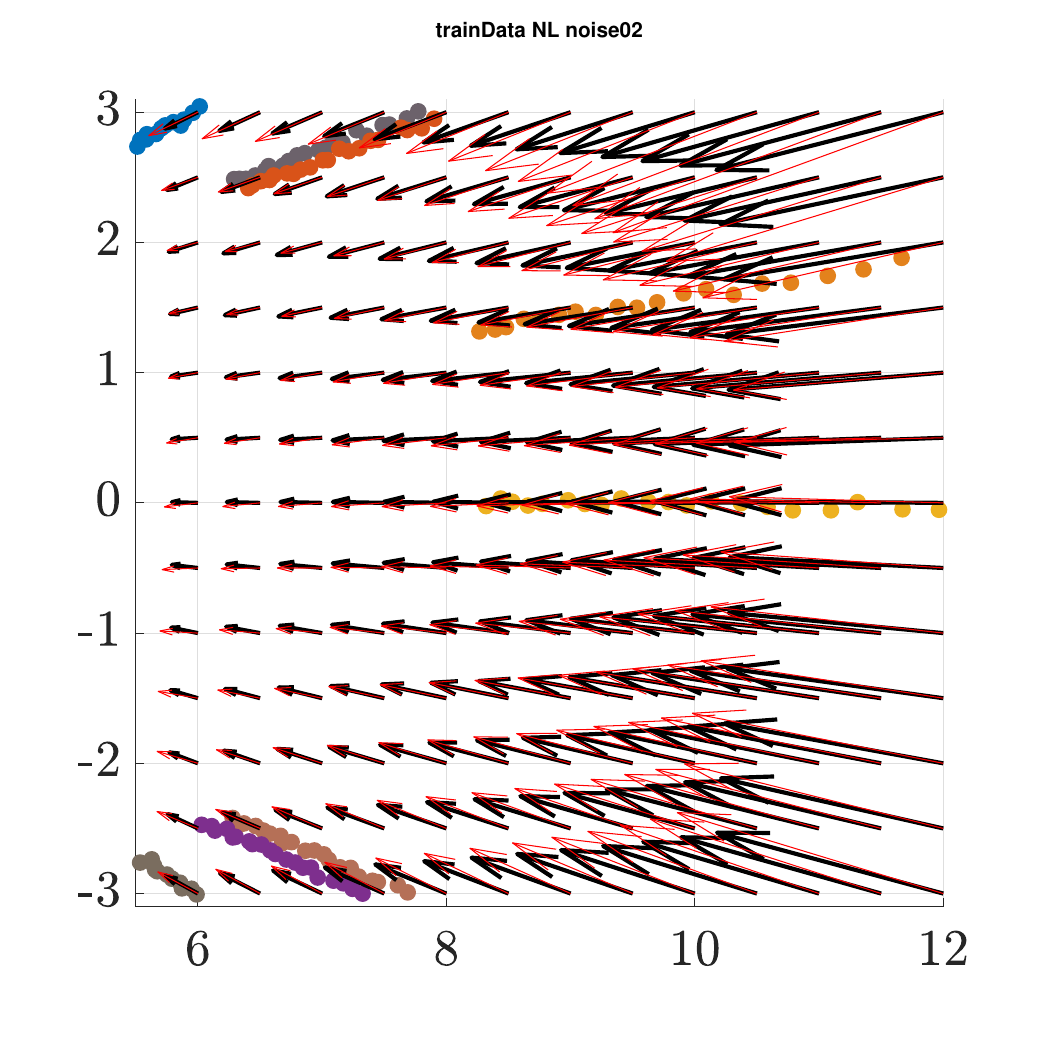}
\caption{{\bf IKEA} -- eight orbits with $\textrm{SNR}=49.1dB$}
        \label{subfig:trainData_NL_noise02_IKEA}
    \end{subfigure}
    \begin{subfigure}[t]{0.495\textwidth} 
\includegraphics[trim=30 30 30 20, clip,width=1\textwidth,valign = t]{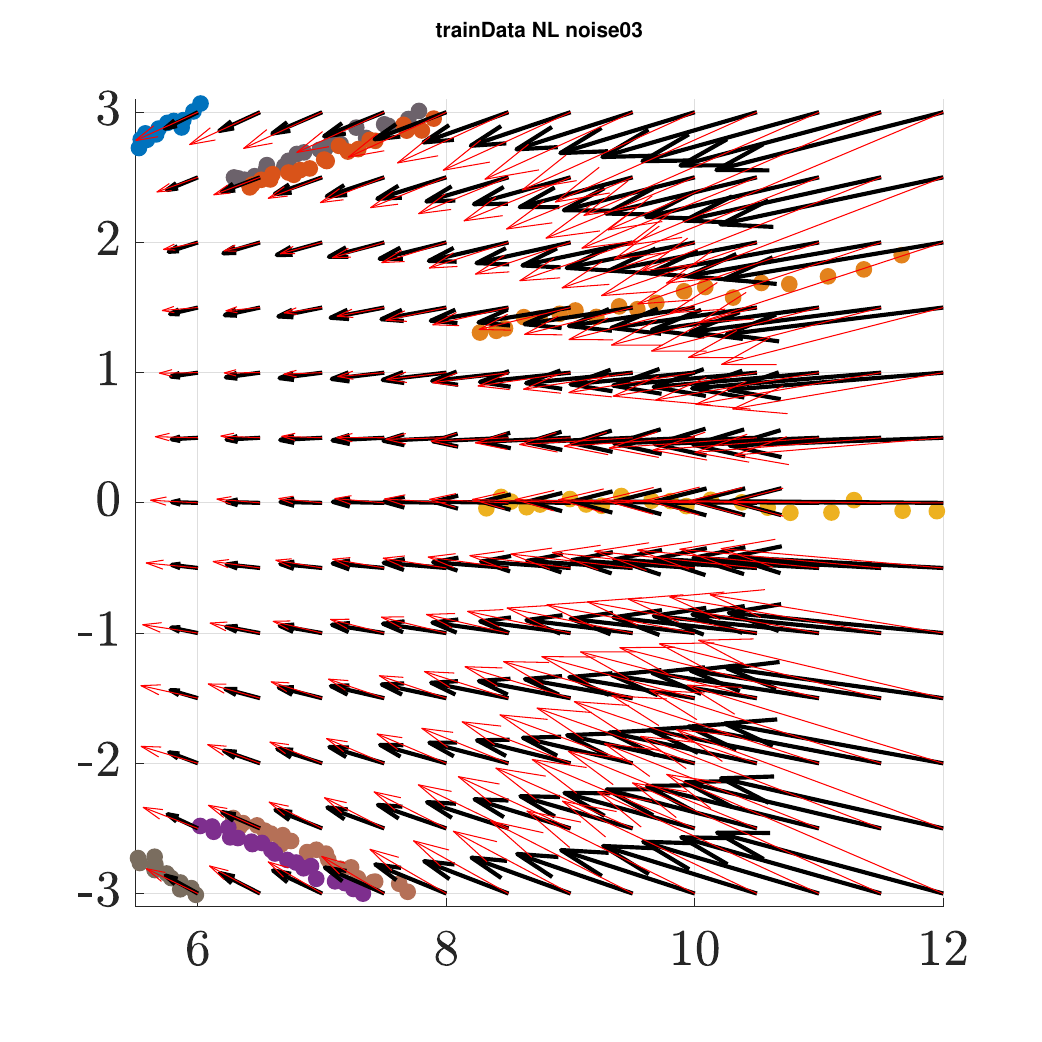}
\caption{{\bf IKEA} -- eight orbits with $\textrm{SNR}=45.6dB$}
        \label{subfig:trainData_NL_noise03_IKEA}
    \end{subfigure}
    \caption{{\bf{IKEA}} -- eight orbits with different levels of noise. Ground truth (black) and restoration (red) of the vector fields are presented.}
    \label{fig:IKEA} 
\end{figure}

We present the results in \Cref{fig:IKEA}. \ac{IKEA} reconstructs the vector field in $\mathcal{D}$ with precision level of $0.0229$, $0.0872$, $0.1674$, and $0.2444$, respectively. These results are presented in \Cref{subfig:trainData_NL_clean_IKEA,subfig:trainData_NL_noise01_IKEA,subfig:trainData_NL_noise02_IKEA,subfig:trainData_NL_noise03_IKEA}, respectively. The sampling rate is constant in all orbits, which lets us compare these results to \acs{SINDy}, \ac{DMD}, and the latter's relatives.

\begin{figure}[phtb!]
    \centering 
    \captionsetup[subfigure]{justification=centering}
    \begin{subfigure}[t]{0.495\textwidth} 
\includegraphics[trim=30 30 30 20, clip,width=1\textwidth,valign = t]{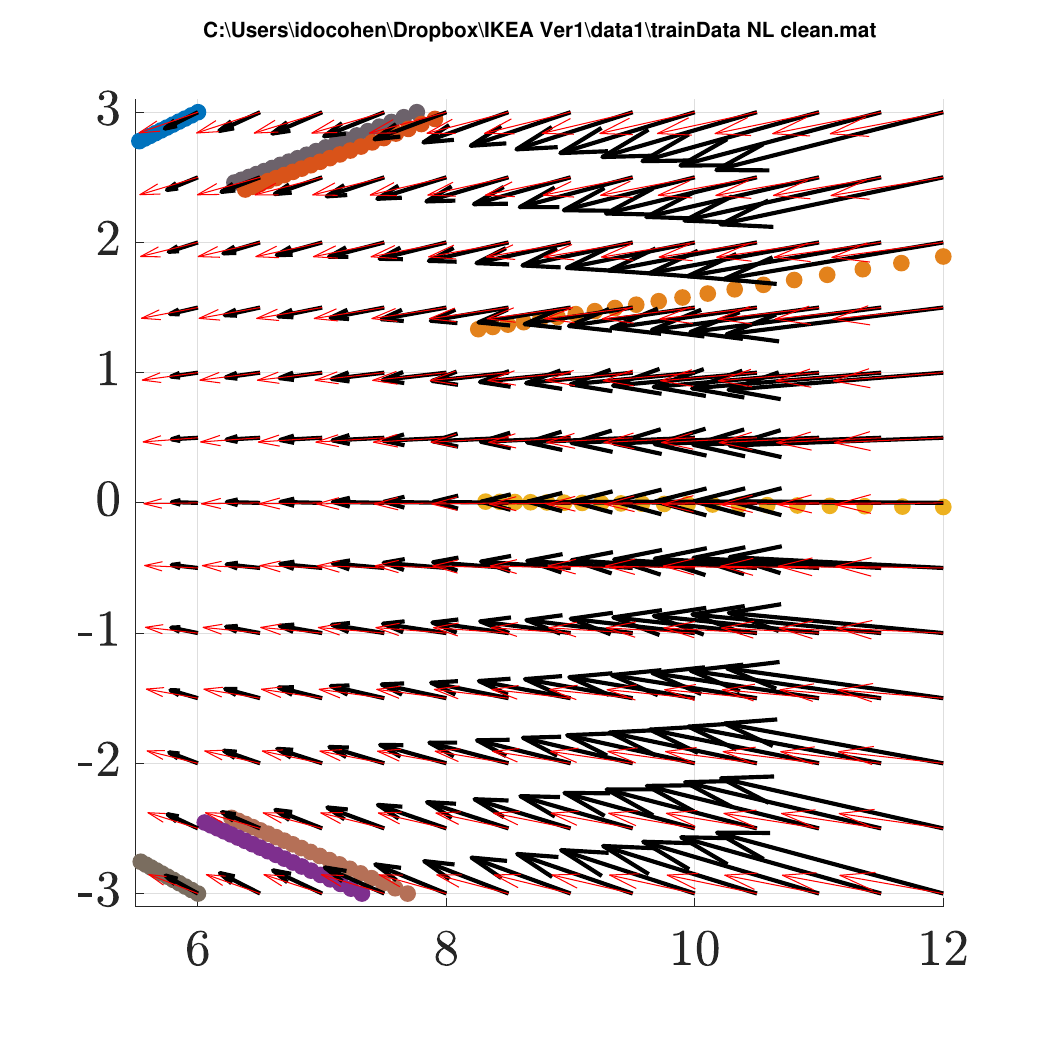}
\caption{{\bf DMD} -- eight orbits without noise}
        \label{subfig:trainData_NL_cleanDMD}
    \end{subfigure}
    \begin{subfigure}[t]{0.495\textwidth} 
\includegraphics[trim=30 30 30 20, clip,width=1\textwidth,valign = t]{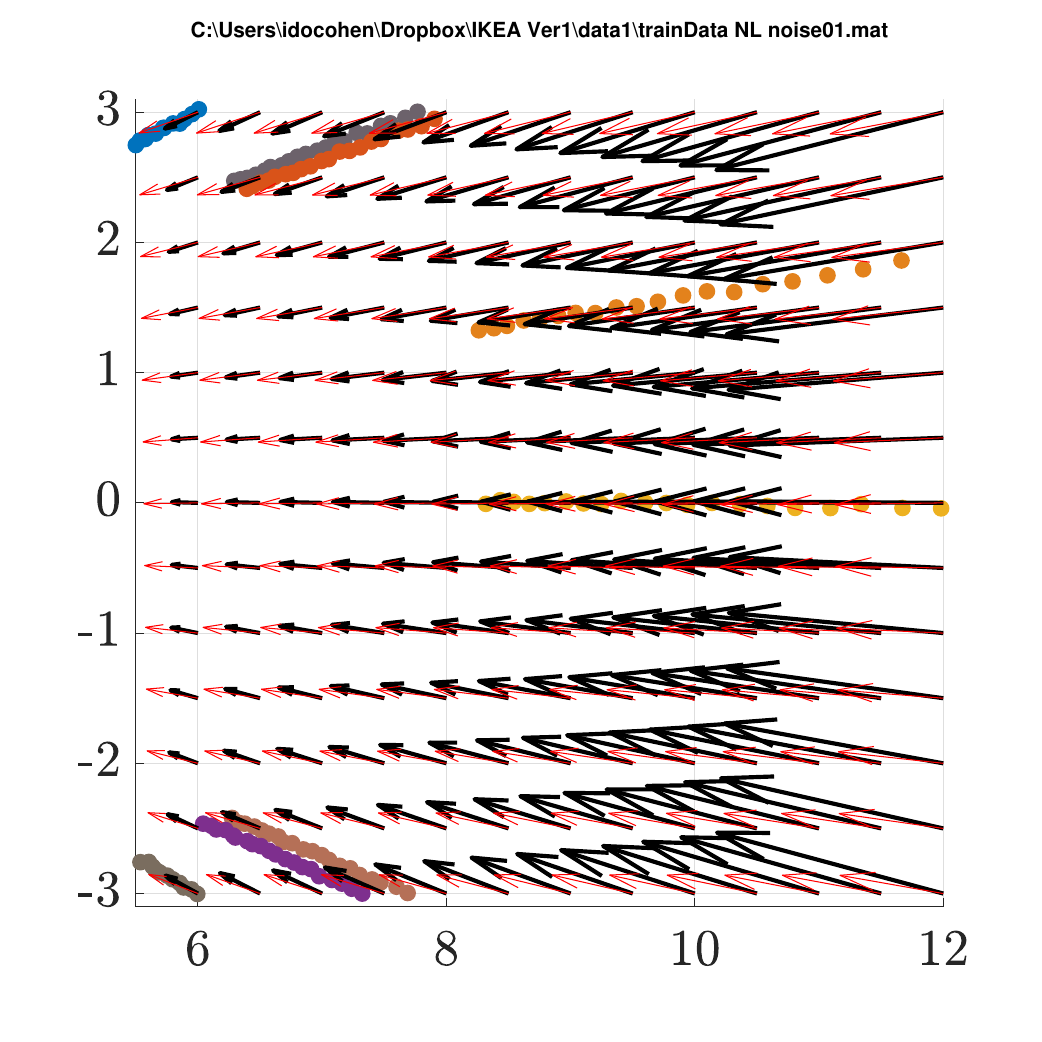}
\caption{{\bf DMD} -- eight orbits with noise $\textrm{SNR}=55.1dB$}
        \label{subfig:trainData_NL_noise01DMD}
    \end{subfigure}\\
    \begin{subfigure}[t]{0.495\textwidth} 
\includegraphics[trim=30 30 30 20, clip,width=1\textwidth,valign = t]{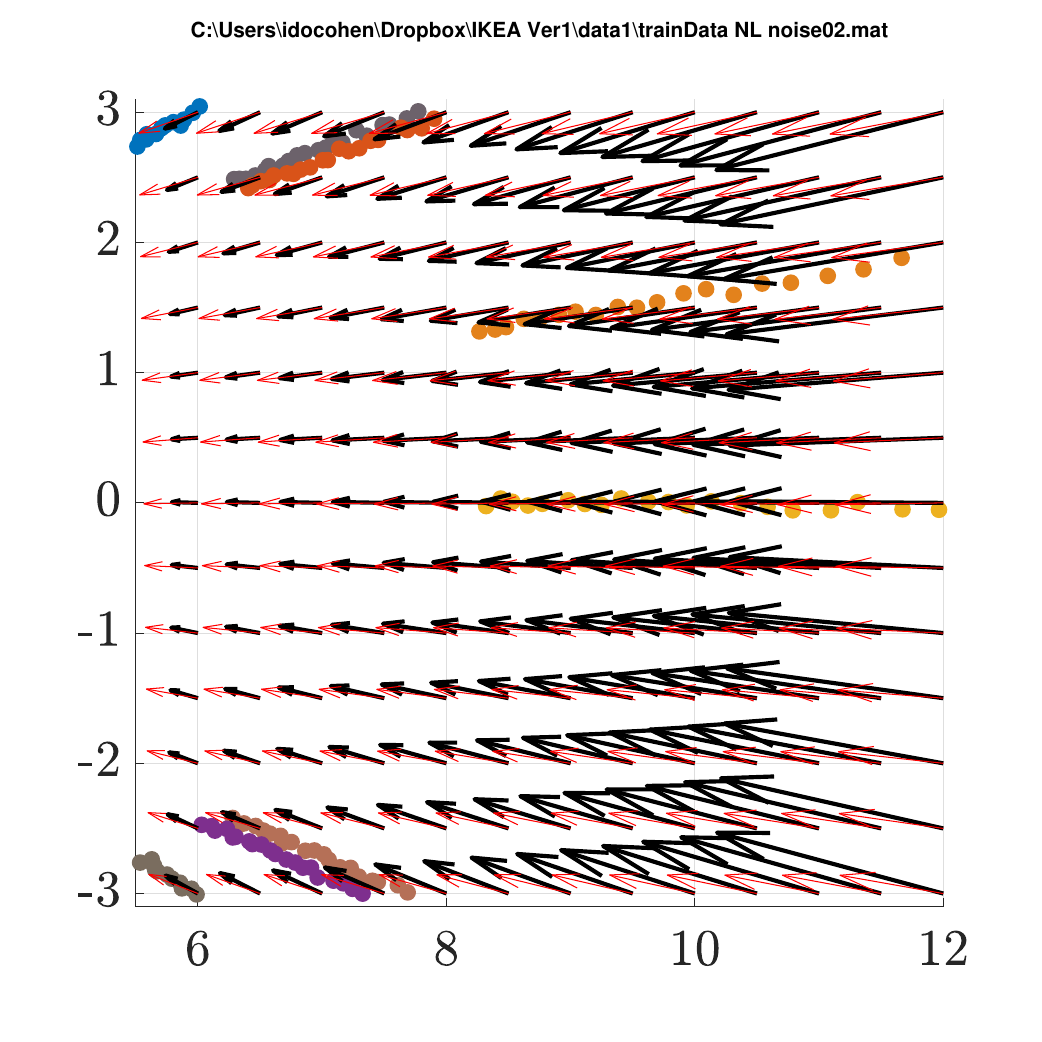}
\caption{{\bf DMD} -- eight orbits with noise $\textrm{SNR}=49.1dB$}
        \label{subfig:trainData_NL_noise02DMD}
    \end{subfigure}
    \begin{subfigure}[t]{0.495\textwidth} 
\includegraphics[trim=30 30 30 20, clip,width=1\textwidth,valign = t]{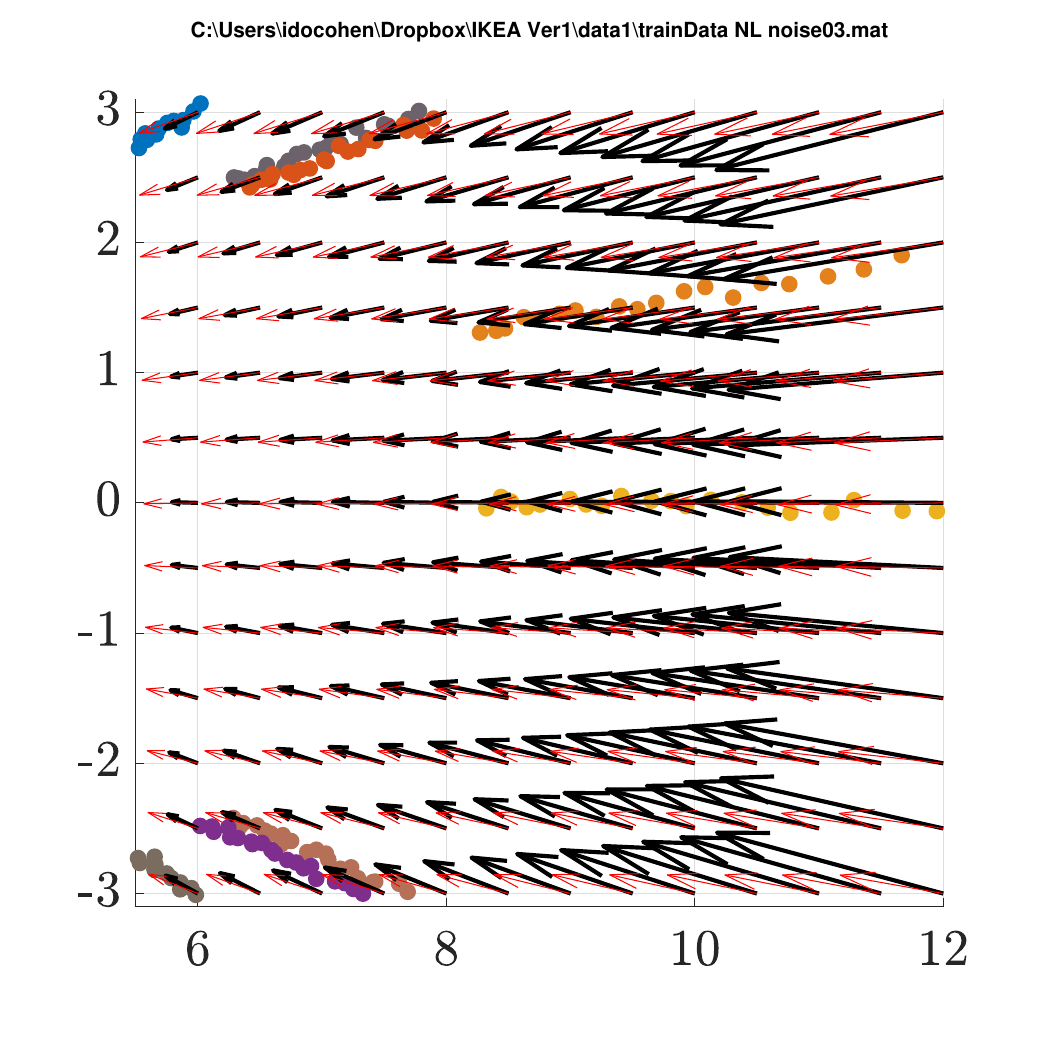}
\caption{{\bf DMD} -- eight orbits with noise $\textrm{SNR}=45.6dB$}
        \label{subfig:trainData_NL_noise03DMD}
    \end{subfigure}

    \caption{{\bf{DMD}} -- eight orbits with different levels of noise. Ground truth (black) and restoration (red) of the vector fields are presented.}
    \label{fig:DMD} 
\end{figure}

\paragraph{Comparison to \ac{DMD}}
This experiment contains eight different orbits from the same vector field. To compare \ac{IKEA} and \ac{DMD}'s family, some adaptations should be made. Usually, the \ac{DMD} algorithm acts on one orbit. In our case, the functional \ac{DMD} to be minimized is a sum of all the functionals of all orbits, correspondingly. The results from \ac{DMD}'s reconstruction are shown in \Cref{fig:DMD}. The errors from GT for these experiments are 
$0.4139$, $0.4139$, $0.4138$, and $0.4134$ and correspond to \Cref{subfig:trainData_NL_cleanDMD,subfig:trainData_NL_noise01DMD,subfig:trainData_NL_noise02DMD,subfig:trainData_NL_noise03DMD}. 

\begin{figure}[phtb!]
    \centering 
    \captionsetup[subfigure]{justification=centering}
    \begin{subfigure}[t]{0.495\textwidth} 
\includegraphics[trim=30 30 30 20, clip,width=1\textwidth,valign = t]{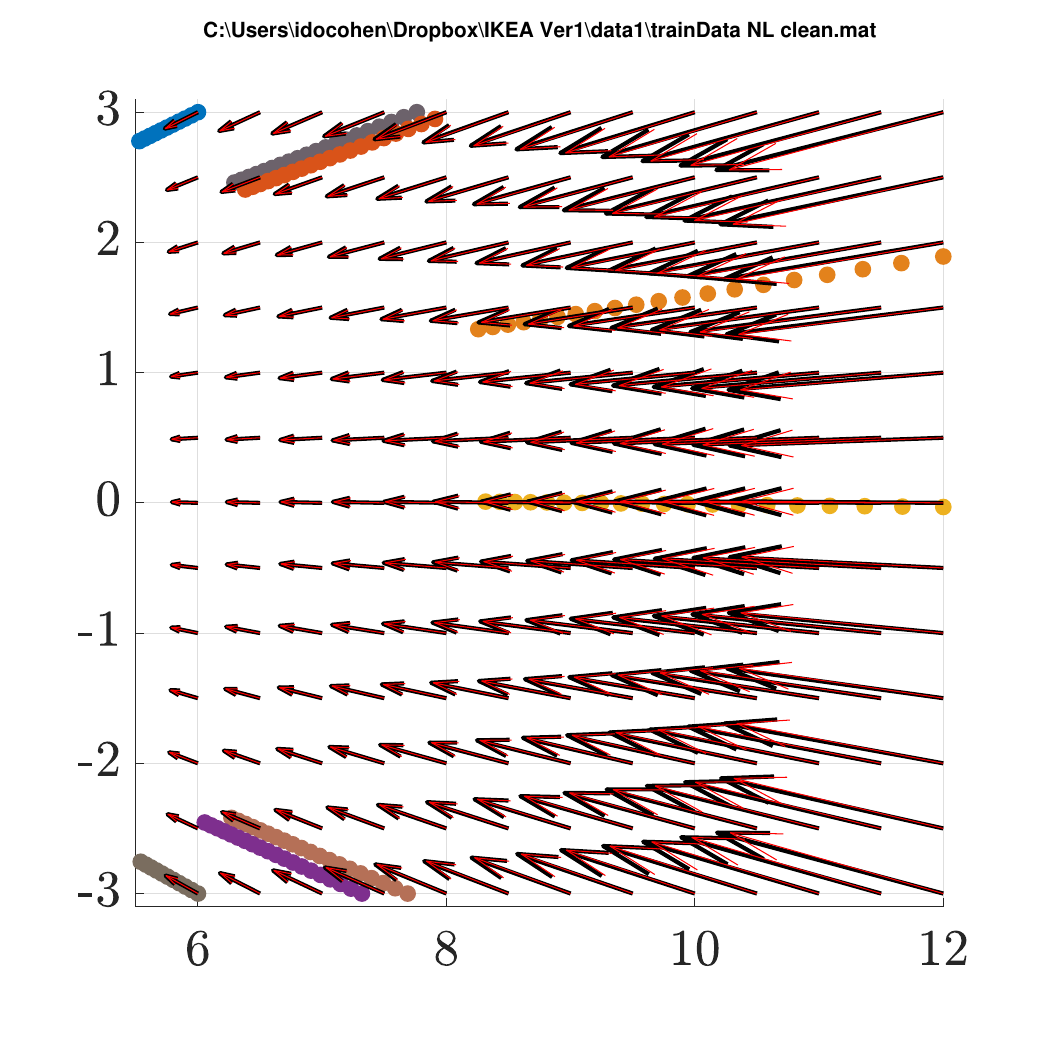}
\caption{{\bf EDMD} -- eight orbits without noise}
        \label{subfig:trainData_NL_cleanEDMD}
    \end{subfigure}
    \begin{subfigure}[t]{0.495\textwidth} 
\includegraphics[trim=30 30 30 20, clip,width=1\textwidth,valign = t]{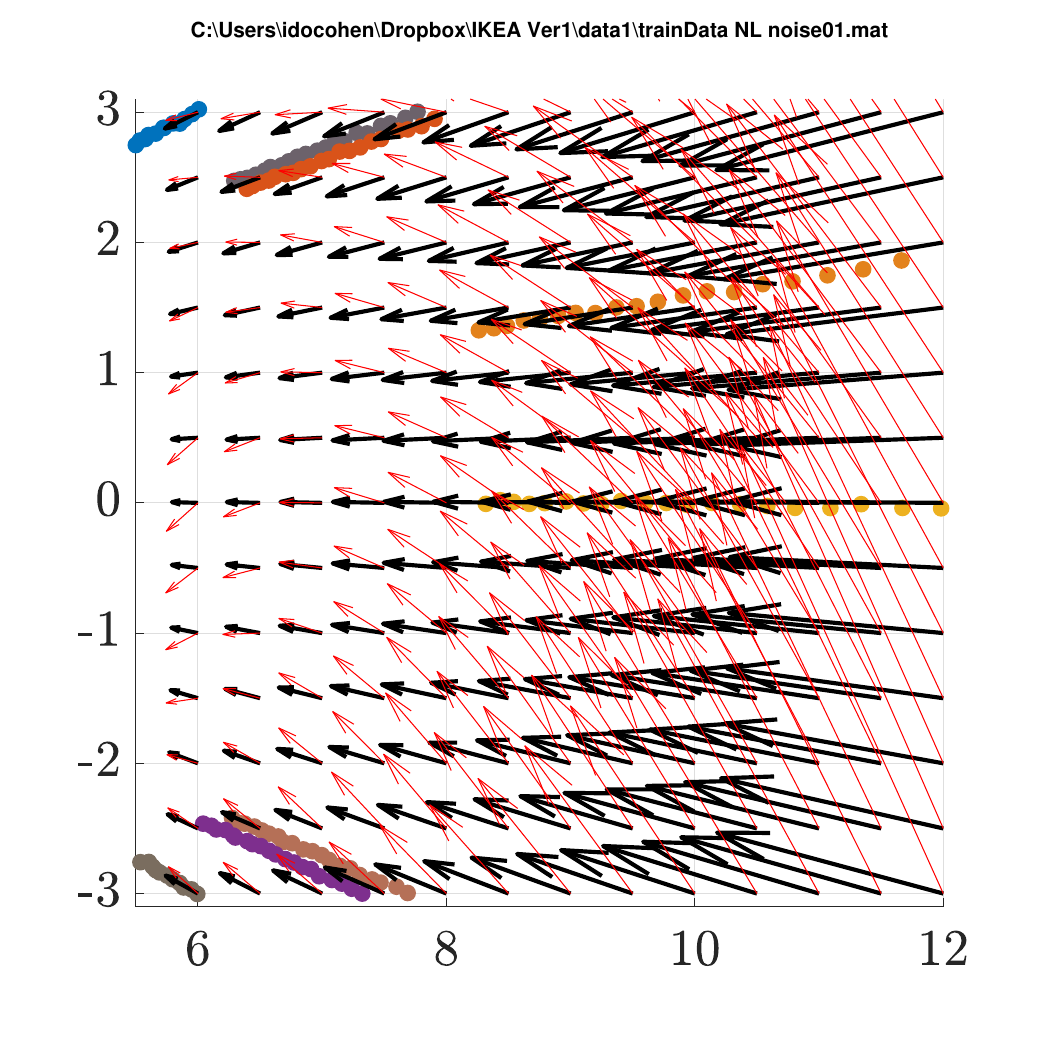}
\caption{{\bf EDMD} -- eight orbits with $\textrm{SNR}=55.1dB$}
        \label{subfig:trainData_NL_noise01EDMD}
    \end{subfigure}\\
    \begin{subfigure}[t]{0.495\textwidth} 
\includegraphics[trim=30 30 30 20, clip,width=1\textwidth,valign = t]{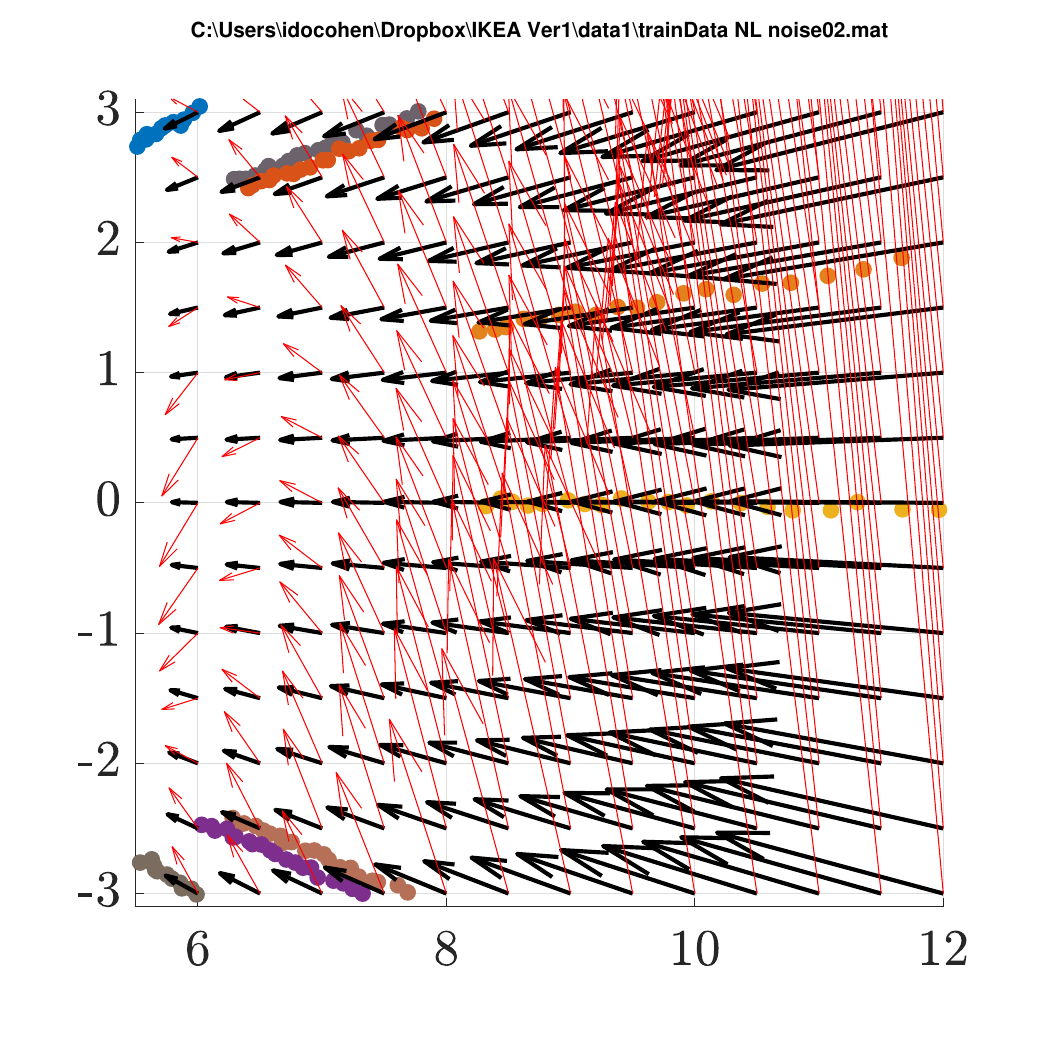}
\caption{{\bf EDMD} -- eight orbits with $\textrm{SNR}=49.1dB$}
        \label{subfig:trainData_NL_noise02EDMD}
    \end{subfigure}
    \begin{subfigure}[t]{0.495\textwidth} 
\includegraphics[trim=30 30 30 20, clip,width=1\textwidth,valign = t]{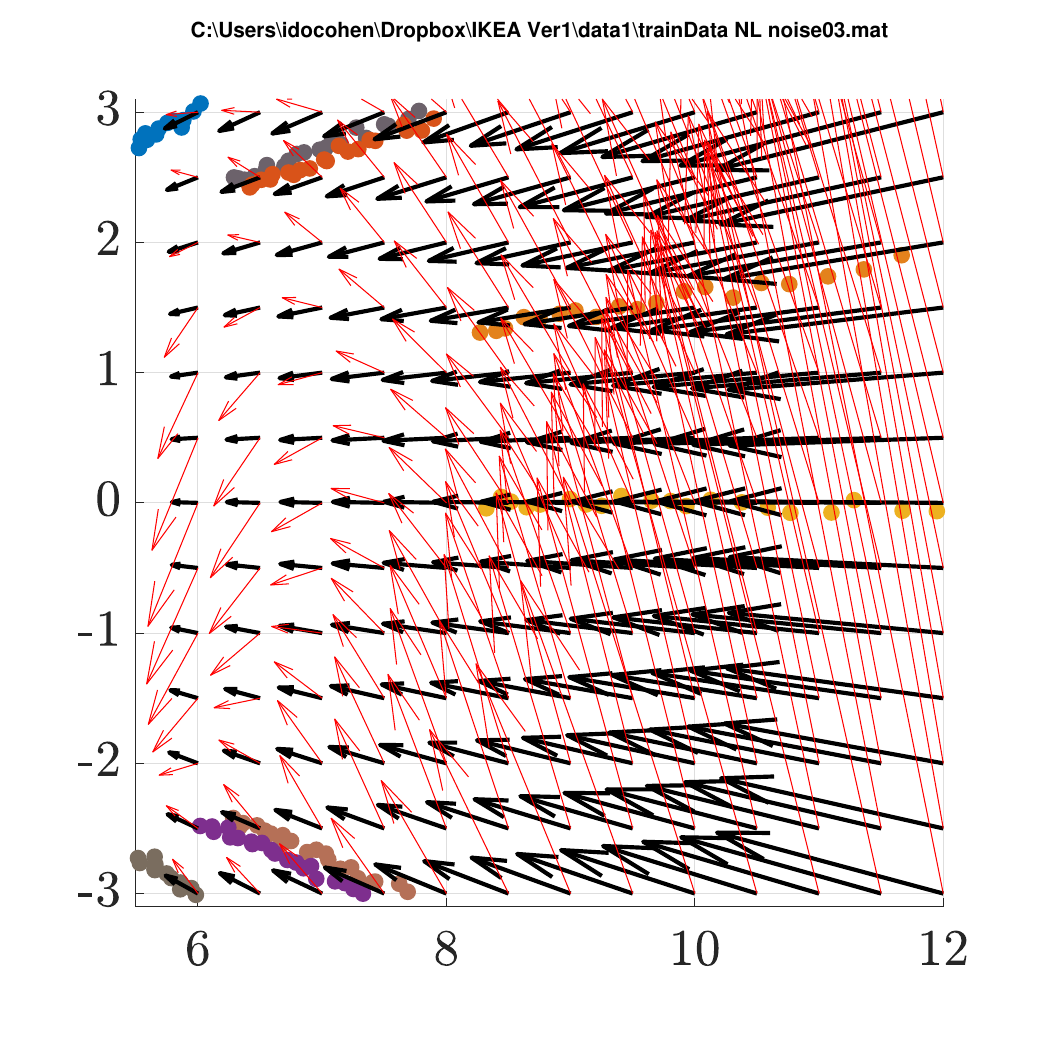}
\caption{{\bf EDMD} -- eight orbits with $\textrm{SNR}=45.6dB$}
        \label{subfig:trainData_NL_noise03EDMD}
    \end{subfigure}

    \caption{{\bf{EDMD}} -- eight orbits with different levels of noise. Ground truth (black) and restoration (red) of the vector fields are presented.}
    \label{fig:EDMD} 
\end{figure}

\paragraph{Comparison to \ac{EDMD}}
We apply these adaptations also to Extended \ac{DMD} and repeat these experiments, using fourth-order polynomial components. The results are depicted in \Cref{fig:EDMD}. The noise levels are the same as in the previous experiment, and the results are brought in \Cref{subfig:trainData_NL_cleanEDMD,subfig:trainData_NL_noise01EDMD,subfig:trainData_NL_noise02EDMD,subfig:trainData_NL_noise03EDMD}, respectively. The error for these results are $0.0463$, $1.1619$, $4.1182$, and $10.2040$, respectively.

\paragraph{Comparison to \texorpdfstring{\acs{SINDy}}{}}
The \acf{SINDy} is a framework for identifying the governing equations of a dynamical system directly from sampled orbits. Consider a dynamical system Eq. \eqref{eq:dynamics}. The central assumption of \ac{SINDy} is that the vector field $p$ can be represented using only a small number of functions from a sufficiently rich dictionary. Thus,
\begin{equation}\label{eq:SINDy}
    \dot{\bm{x}} = \bm{\Theta}(\bm{x})\bm{\Xi},
\end{equation}
where $\bm{\Theta}(\bm{x})$ is a dictionary of functions and $\bm{\Xi}$ is a sparse coefficient matrix (for more details, see \cite{brunton2016discovering}).

In our case, for a two-dimensional state $\bm{x}= \begin{bmatrix} x_1& x_2 \end{bmatrix}^T$, we choose 
a dictionary with fourth-order polynomials that contains

\begin{equation}\label{eq:SINDyDictionary}
    \bm{\Theta}(\bm{x}) = [ 1,\, x_1,\, x_2,\, x_1^2,\, \ldots ,\,x_1^4,\, x_1^3x_2,\,x_1^2x_2^2,\,x_1x_2^3,\,x_2^4].    
\end{equation}

\begin{figure}[phtb!]
    \centering 
    \captionsetup[subfigure]{justification=centering}
    \begin{subfigure}[t]{0.495\textwidth} 
\includegraphics[trim=30 30 30 20, clip,width=1\textwidth,valign = t]{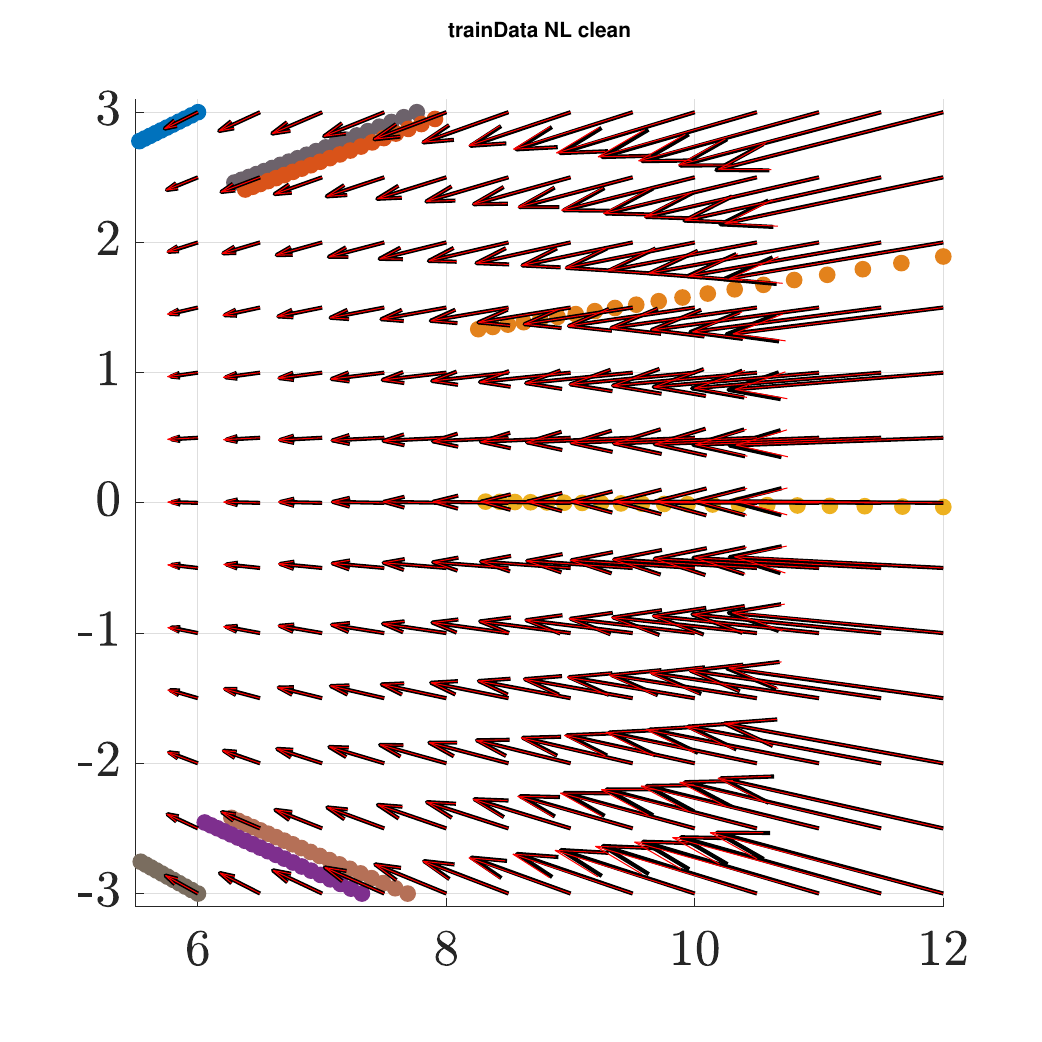}
\caption{{\bf SINDy} -- eight orbits without noise}
        \label{subfig:trainData_NL_cleanSidny}
    \end{subfigure}
    \begin{subfigure}[t]{0.495\textwidth} 
\includegraphics[trim=30 30 30 20, clip,width=1\textwidth,valign = t]{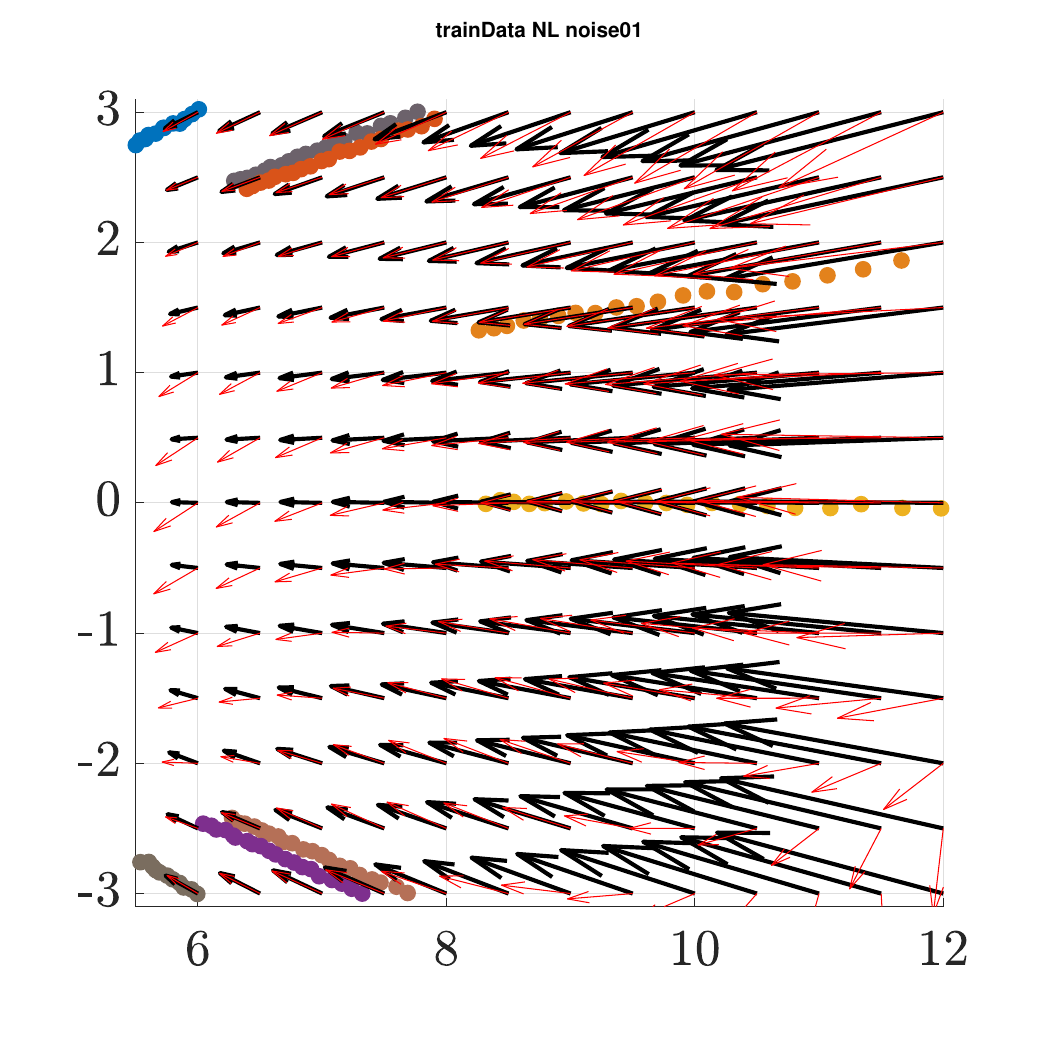}
\caption{{\bf SINDy} -- eight orbits with $\textrm{SNR}=55.1dB$}
        \label{subfig:trainData_NL_noise01Sidny}
    \end{subfigure}\\
    \begin{subfigure}[t]{0.495\textwidth} 
\includegraphics[trim=30 30 30 20, clip,width=1\textwidth,valign = t]{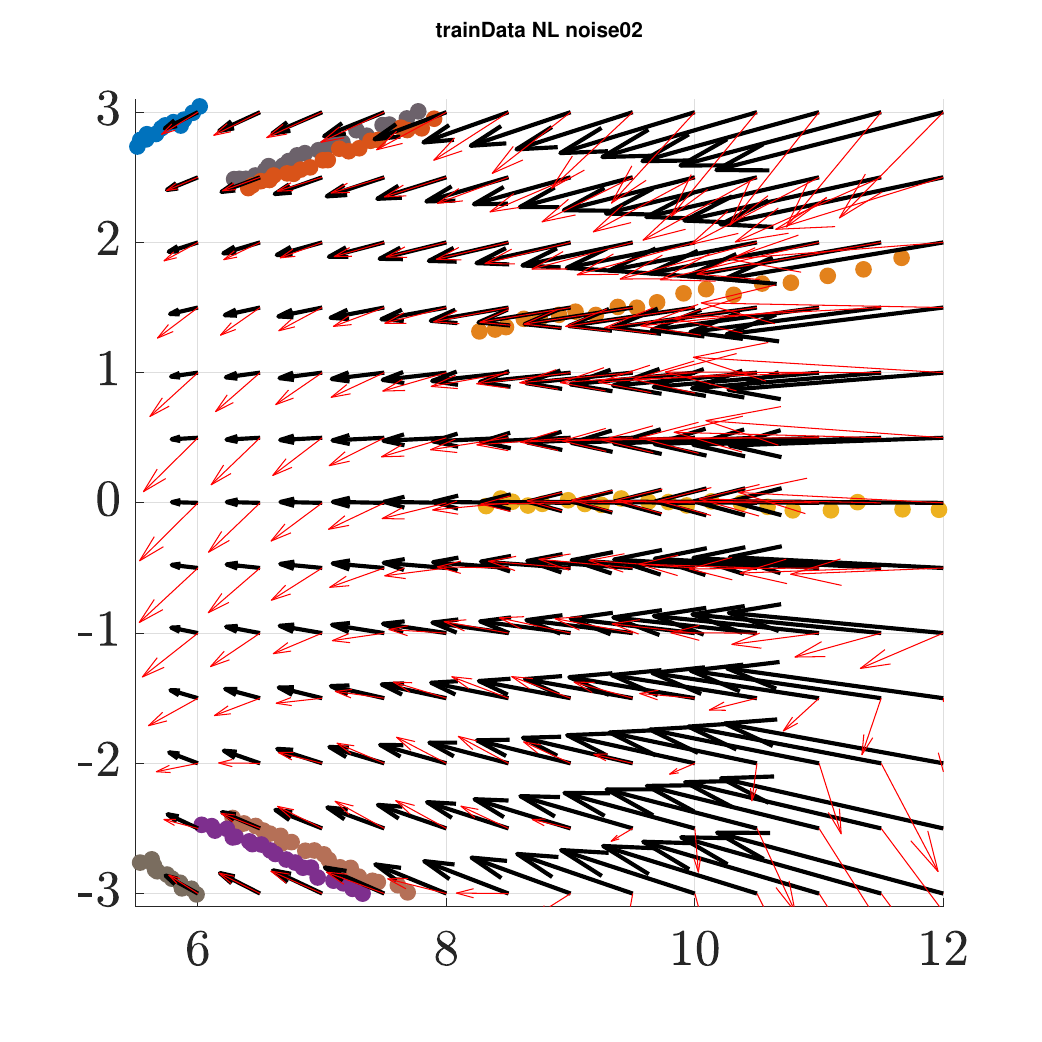}
\caption{{\bf SINDy} -- eight orbits with $\textrm{SNR}=49.1dB$}
        \label{subfig:trainData_NL_noise02Sidny}
    \end{subfigure}
    \begin{subfigure}[t]{0.495\textwidth} 
\includegraphics[trim=30 30 30 20, clip,width=1\textwidth,valign = t]{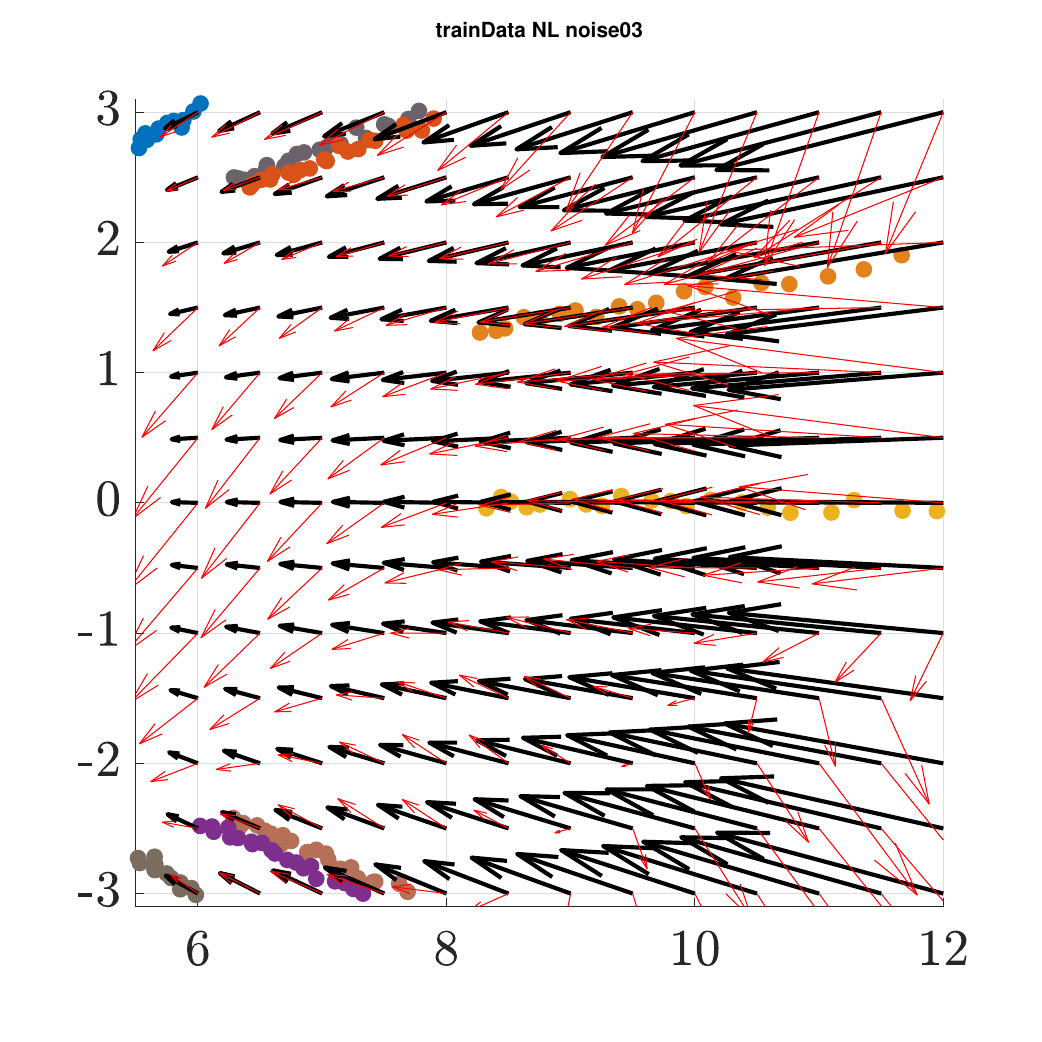}
\caption{{\bf SINDy} -- eight orbits with $\textrm{SNR}=45.6dB$}
        \label{subfig:trainData_NL_noise03Sidny}
    \end{subfigure}

    \caption{{\bf{SINDy}} -- eight orbits with different levels of noise. Ground truth (black) and restoration (red) of the vector fields are presented.}
    \label{fig:SIDNy} 
\end{figure}

The \ac{SINDy} procedure determines which of these candidate functions are actually required to describe the dynamics. We repeat the previous experiments with the same settings. The results are depicted in \Cref{subfig:trainData_NL_cleanSidny,subfig:trainData_NL_noise01Sidny,subfig:trainData_NL_noise02Sidny,subfig:trainData_NL_noise03Sidny}. The errors (Eq. \eqref{eq:errorVF}) are $0.0200$, $0.4611$, $0.9549$, and $1.3866$, respectively.

The following table (\Cref{table:results}) summarizes the results of the experiments on \ac{IKEA} compared to other methods.

\begin{table}[pthb!]
        \caption{Result Summary (Error \eqref{eq:errorVF})}
    \label{table:results}
\centering
\begin{tabular}{|c|c |c| c|c|}
 \hline
 SNR&\acs{IKEA}&\acs{DMD}& \acs{EDMD}&\acs{SINDy}\\
 \hline\hline
Without noise&$0.0229$& $0.4139$& $0.0463$&$\bm{0.0200}$\\
\hline
$55.1$&$\bm{0.0872}$& $0.4139$ &$1.1619$&$0.4611$\\
\hline
$49.1$&$\bm{0.1674}$& $0.4138$ &$4.1182$ &$0.9549$\\
\hline
$45.6$&$\bm{0.2444}$& $0.4134$ & $10.2040$  & $1.3866$\\
\hline

\end{tabular}
 
\end{table}

\begin{figure}[phtb!]
    \centering 
    \captionsetup[subfigure]{justification=centering}
    \begin{subfigure}[t]{0.495\textwidth} 
\includegraphics[trim=30 30 30 20, clip,width=1\textwidth,valign = t]{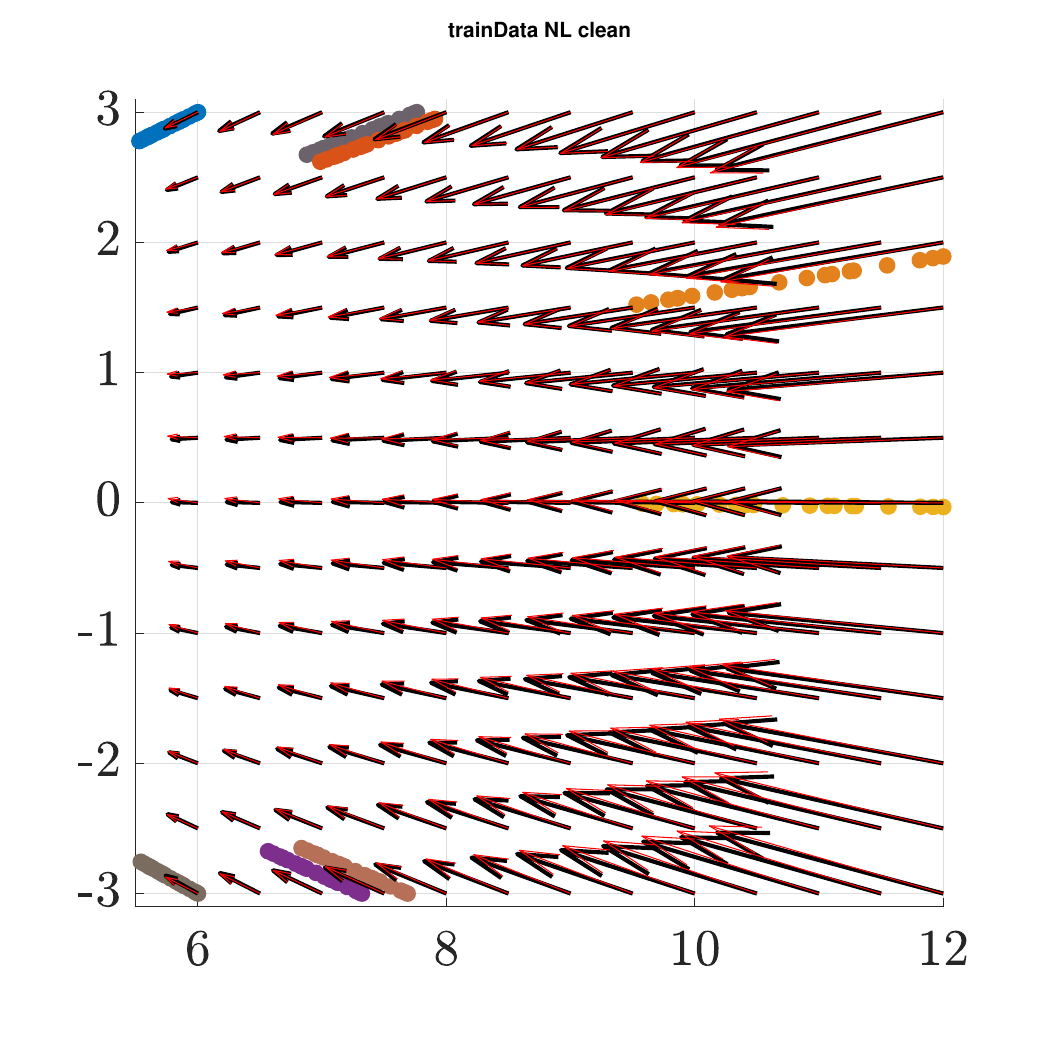}
\caption{{\bf IKEA} -- eight orbits without noise}
        \label{subfig:trainData_NL_clean_IKEAS}
    \end{subfigure}
    \begin{subfigure}[t]{0.495\textwidth} 
\includegraphics[trim=30 30 30 20, clip,width=1\textwidth,valign = t]{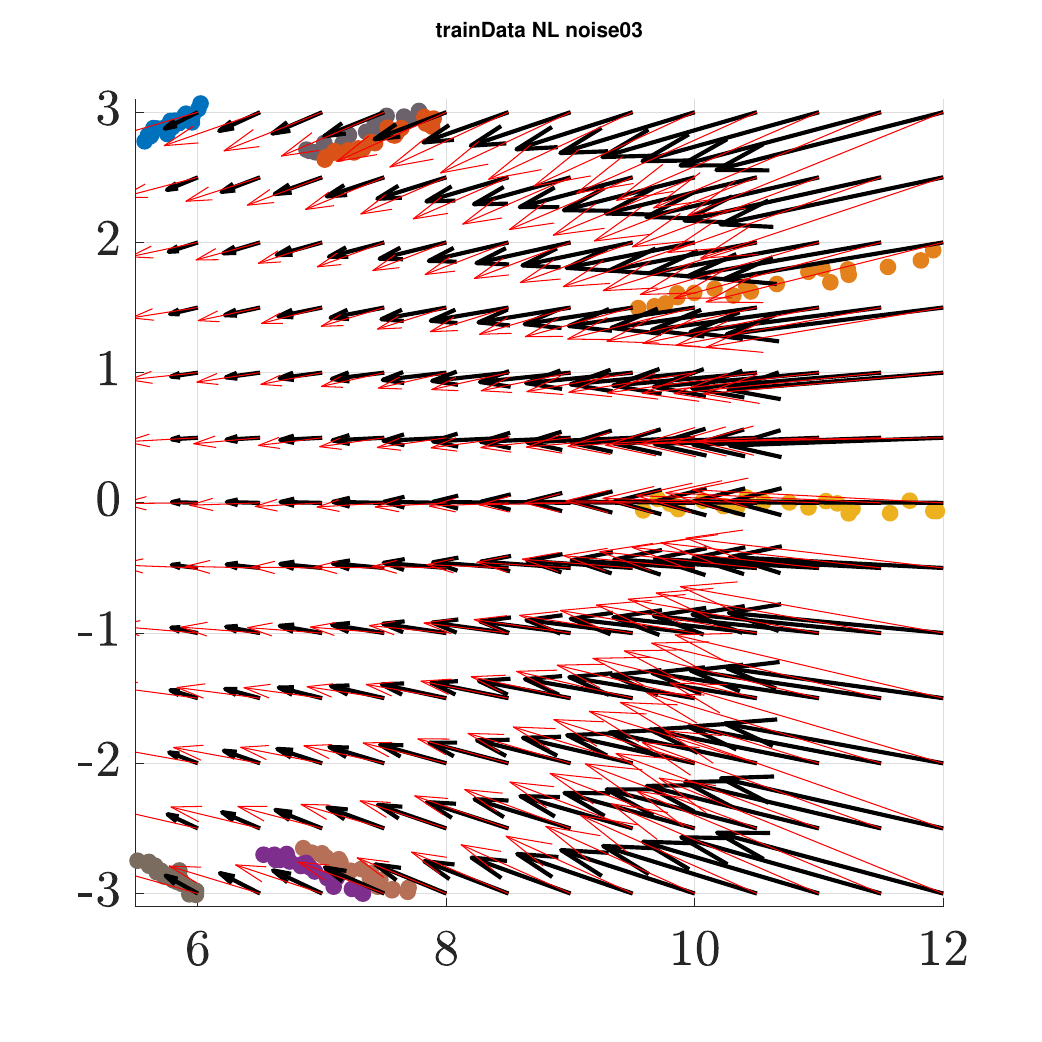}
\caption{{\bf IKEA} -- eight orbits with $\textrm{SNR}=45.6dB$}
        \label{subfig:trainData_NL_noise03_IKEAS}
    \end{subfigure}
    \caption{{\bf{IKEA}} -- eight orbits with different levels of noise. Ground truth (black) and restoration (red) of the vector fields are presented.}
    \label{fig:IKEAS} 
\end{figure}

\niceparagraph{Variability in Temporal Sampling} Now, we repeat the experiments from \Cref{subfig:trainData_NL_clean_IKEA} and \Cref{subfig:trainData_NL_noise03_IKEA} when the time step size is chosen randomly $dt\sim U[0,0.2]$. We show the results in \Cref{fig:IKEAS}. The errors (Eq. \eqref{eq:errorVF}) are $0.0248$ and $0.2653$, respectively. The settings for this experiment are different from the previous ones; therefore, these results are excluded from \Cref{table:results}.

\subsection{Discussion}
One of the advantages of \ac{DMD} and its family is that it was designed to find the dynamics that minimize the RMSE from the actual samples. Although \ac{DMD} achieves inaccurate results for one sampled orbit \cite{cohen2021modes}, applying it to a few orbits can significantly improve the results. However, the inherent problems of this method are revealed when the system is not homogeneous, or when it is not in $C^1$ (c.f. \cite{cohen2021latent,wu2021challenges,brokman2024enhancing}). \ac{EDMD} may circumvent some of these problems and achieve good results (as in \Cref{subfig:trainData_NL_cleanEDMD}); however, it dramatically fails in the presence of noise. \ac{EDMD} overfits the data; hence, its sensitivity to noise. In addition, the almost perfect result in \Cref{subfig:trainData_NL_cleanEDMD} is due to the optional functions containing the correct ones. Therefore, the functions in \ac{EDMD} must be versatile enough to restore the correct dynamical system, which may cause overfitting in the presence of noise. The same arguments hold for \ac{SINDy}, since it is also a dictionary-based approach.

\ac{KR} and \ac{IKEA} present a new point of view in the field of inverse problems. ROF and Tikhonov models, and even \ac{PINN},  combine two expressions: the fidelity term and a regularization term. That is, these models distinguish between the observables and their expected interaction. On the other hand, \ac{KR} and \ac{IKEA} force the results to address a geometric condition that the dynamical system must satisfy. This constraint acts as a regularization term. It does not mitigate the objective functional, but it prevents singularity and hence smooths the results. 

\section{Conclusions}\label{sec:conclusions}

This work presented two algorithms for data-driven dynamical-system reconstruction, \ac{KR}, which uses samples of the vector field, and \ac{IKEA}, which uses samples of orbits. Both methods formulate the reconstruction as a constrained optimization problem and seek a functionally independent set of \acp{UVM}.

The main contribution of the work is the formulation of conditions for functional independence based on the Jacobian and the Gershgorin Circle Theorem. These conditions are incorporated into the optimization problem as constraints.


The numerical experiments demonstrate the application of \ac{KR} and \ac{IKEA} to vector-field reconstruction, including noise reduction, generalization from sparse samples, and dimensionality reduction. The results also show robustness to noisy data and variability in temporal sampling. Overall, the proposed methods provide a geometry-based approach to dynamical-system reconstruction from sampled vector fields and trajectories.

This paper presents a new methodology and demonstrate its performance in low-dimensional systems. The next step in the research is to check the scalability of this approach and examine ways to circumvent singularities and close orbits. In addition, these algorithms, breaking data down into motifs in the sense of \ac{KEF}, are applicable in signal processing and data mining. The connection between vector field representation and data points cannot be ignored, and this work can be a cornerstone for forthcoming research.

\section*{Acronym List}
\begin{acronym}
\acro{DMD}[DMD] {\emph{Dynamic Mode Decomposition}}
\acro{EDMD}[EDMD] {\emph{Extended Dynamic Mode Decomposition}}
\acro{IKEA}[IKEA] {\emph{Independent Koopman Eigenfunction Approximation}}
\acro{KAE}[KAE]{\emph{Koopman Autoencoder}}
\acro{KEF}[KEF]{\emph{Koopman Eigenfunction}}
\acro{KPDE}[KPDE]{\emph{Koopman Partial Differential Equation}}
\acro{KR}[KR] {\emph{Koopman Regularization}}
\acro{PDE}[PDE]{\emph{Partial Differential Equation}}
\acro{PINN}[PINN] {\emph{Physics-Informed Neural Networks}}
\acro{SINDy}[SINDy] {\emph{Sparse Identification of Nonlinear Dynamics}}
\acro{UVM}[UVM] {\emph{Unit Velocity Measurement}}
\end{acronym}

\begin{appendices}

\section*{Acknowledgments}
``He who walks with the wise grows wise'' Proverbs 13:20. The author acknowledges Prof. Guy Gilboa, Prof. Gershon Wolansky, Dr. Eli Appleboim, Mr. (soon Dr.) Meir Yossef Levi, and Mr. Gilad Feldman for their good advice and support. In addition, the author thanks the reviewers for suggesting that we emphasize the importance of functional independence for system reconstruction.

\end{appendices}
\bibliographystyle{amsplain}
\bibliography{smartPeople}

@Inbook{Arnold1989-1,
author="Arnold, V. I.",
title="Lagrangian mechanics on manifolds",
bookTitle="Mathematical Methods of Classical Mechanics",
year="1989",
publisher="Springer New York",
address="New York, NY",
pages="75--97",
isbn="978-1-4757-2063-1",
doi="10.1007/978-1-4757-2063-1_4",
url="https://doi.org/10.1007/978-1-4757-2063-1_4"
}

@article{raissi2017physics,
  title={Physics informed deep learning (part i): Data-driven solutions of nonlinear partial differential equations},
  author={Raissi, Maziar and Perdikaris, Paris and Karniadakis, George Em},
  journal={arXiv preprint arXiv:1711.10561},
  pages = {1-22},
  year={2017}
}

@article{gershgorin1931uber,
  title        = {Über die Abgrenzung der Eigenwerte einer Matrix},
  author       = {Gershgorin, Semen Aronovich},
  journal      = {Izvestiya Rossiiskoi Akademii Nauk. Seriya Matematicheskaya},
  volume       = {7},
  number       = {6},
  pages        = {749--754},
  year         = {1931},
  publisher    = {Russian Academy of Sciences, Steklov Mathematical Institute}
}

@article{wu2021challenges,
  title={Challenges in dynamic mode decomposition},
  author={Wu, Ziyou and Brunton, Steven L and Revzen, Shai},
  journal={Journal of the Royal Society Interface},
  volume={18},
  number={185},
  pages={20210686},
  year={2021},
  publisher={The Royal Society}
}

@misc{gilboa2024system,
  title={System and method of training a neural network model},
  author={Gilboa, Guy and Turjeman, Rotem and Berkov, Tom and Cohen, Ido},
  year={2024},
  month=feb # "~15",
  publisher={Google Patents},
  note={US Patent App. 18/231,968}
}

@misc{Folland1995Introduction,
publisher = {Princeton University Press},
title = {Introduction to Partial Differential Equations},
year = {1995},
author = {Folland, G. B.},
address = {Princeton, N.J},
booktitle = {Introduction to Partial Differential Equations},
isbn = {1-5231-4842-X},
language = {eng},
}

@book{zorich2016mathematical,
  title={Mathematical analysis II},
  author={Zorich, Vladimir Antonovich and Paniagua, Octavio},
  volume={220},
  year={2016},
  publisher={Springer}
}

@misc{cohen2023minimal,
      title={A Minimal Set of Koopman Eigenfunctions -- Analysis and Numerics}, 
      author={Ido Cohen and Eli Appleboim},
      year={2023},
      eprint={2303.05837},
      archivePrefix={arXiv},
      primaryClass={math.DS}
}

@article{koopman1931hamiltonian,
  title={Hamiltonian systems and transformation in Hilbert space},
  author={Koopman, Bernard O},
  journal={Proceedings of the national academy of sciences of the united states of america},
  volume={17},
  number={5},
  pages={315},
  year={1931},
  publisher={National Academy of Sciences}
}

@misc{lawrence1991differential,
  title={Differential equations and dynamical systems},
  author={Lawrence, Perko},
  year={1991},
  publisher={Springer-Verlag, New York},
  pages={9}
}

@ARTICLE{9454411,
  author={Haseli, Masih and Cortés, Jorge},
  journal={IEEE Transactions on Control of Network Systems}, 
  title={Parallel Learning of Koopman Eigenfunctions and Invariant Subspaces for Accurate Long-Term Prediction}, 
  year={2021},
  volume={8},
  number={4},
  pages={1833-1845},
  doi={10.1109/TCNS.2021.3088791}
}

@article{cohen2021modes,
  title={Modes of homogeneous gradient flows},
  author={Cohen, Ido and Azencot, Omri and Lifshits, Pavel and Gilboa, Guy},
  journal={SIAM Journal on Imaging Sciences},
  volume={14},
  number={3},
  pages={913--945},
  year={2021},
  publisher={SIAM}
}

@article{cohen2021latent,
  title={Latent Modes of Nonlinear Flows--a Koopman Theory Analysis},
  author={Cohen, Ido and Gilboa, Guy},
  journal={arXiv preprint arXiv:2107.07456},
pages={26},
  year={2021}
}

@article{askham2018variable,
  title={Variable projection methods for an optimized dynamic mode decomposition},
  author={Askham, Travis and Kutz, J Nathan},
  journal={SIAM Journal on Applied Dynamical Systems},
  volume={17},
  number={1},
  pages={380--416},
  year={2018},
  publisher={SIAM}
}

@article{brunton2016discovering,
  title={Discovering governing equations from data by sparse identification of nonlinear dynamical systems},
  author={Brunton, Steven L and Proctor, Joshua L and Kutz, J Nathan},
  journal={Proceedings of the national academy of sciences},
  volume={113},
  number={15},
  pages={3932--3937},
  year={2016},
  publisher={National Acad Sciences}
}

@article{li2017extended,
  title={Extended dynamic mode decomposition with dictionary learning: A data-driven adaptive spectral decomposition of the Koopman operator},
  author={Li, Qianxiao and Dietrich, Felix and Bollt, Erik M and Kevrekidis, Ioannis G},
  journal={Chaos: An Interdisciplinary Journal of Nonlinear Science},
  volume={27},
  pages={--},
  number={10},
  year={2017},
  publisher={AIP Publishing}
}

@book{cohen_gilboa_2023, place={Cambridge}, series={Elements in Non-local Data Interactions: Foundations and Applications}, title={Latent Modes of Nonlinear Flows: A Koopman Theory Analysis}, publisher={Cambridge University Press}, author={Cohen, Ido and Gilboa, Guy}, year={2023}, collection={Elements in Non-local Data Interactions: Foundations and Applications}}

@article{schmid2010dynamic,
  title={Dynamic mode decomposition of numerical and experimental data},
  author={Schmid, Peter J},
  journal={Journal of fluid mechanics},
  volume={656},
  pages={5--28},
  year={2010},
  publisher={Cambridge University Press}
}

@article{schmid2022dynamic,
  title={Dynamic mode decomposition and its variants},
  author={Schmid, Peter J},
  journal={Annual Review of Fluid Mechanics},
  volume={54},
  pages={225--254},
  year={2022},
  publisher={Annual Reviews}
}

@article{OPTbao2004computing,
  title={Computing the ground state solution of Bose--Einstein condensates by a normalized gradient flow},
  author={Bao, Weizhu and Du, Qiang},
  journal={SIAM Journal on Scientific Computing},
  volume={25},
  number={5},
  pages={1674--1697},
  year={2004},
  publisher={SIAM}
}

@article{OPTgarcia2001optimizing,
  title={Optimizing Schr{\"o}dinger functionals using Sobolev gradients: Applications to quantum mechanics and nonlinear optics},
  author={Garc{\'\i}a-Ripoll, Juan Jos{\'e} and P{\'e}rez-Garc{\'\i}a, V{\'\i}ctor M},
  journal={SIAM Journal on Scientific Computing},
  volume={23},
  number={4},
  pages={1316--1334},
  year={2001},
  publisher={SIAM}
}

@article{OPTbao2003ground,
  title={Ground-state solution of Bose--Einstein condensate by directly minimizing the energy functional},
  author={Bao, Weizhu and Tang, Weijun},
  journal={Journal of Computational Physics},
  volume={187},
  number={1},
  pages={230--254},
  year={2003},
  publisher={Elsevier}
}

@article{OPTcaliari2009minimisation,
  title={A minimisation approach for computing the ground state of Gross--Pitaevskii systems},
  author={Caliari, Marco and Ostermann, Alexander and Rainer, Stefan and Thalhammer, Mechthild},
  journal={Journal of Computational Physics},
  volume={228},
  number={2},
  pages={349--360},
  year={2009},
  publisher={Elsevier}
}

@book{OPTekeland1999convex,
  title={Convex analysis and variational problems},
  author={Ekeland, Ivar and Temam, Roger},
  volume={28},
  year={1999},
  publisher={Siam}
}

@article{OPTdem2004exact,
  title={Exact penalty functions and problems of variation calculus},
  author={Dem'yanov, Vladimir Fedorovich},
  journal={Automation and Remote Control},
  volume={65},
  number={2},
  pages={280--290},
  year={2004},
  publisher={Springer}
}

@article{OPTdemyanov2011exact,
  title={Exact penalty functions in isoperimetric problems},
  author={Demyanov, VF and Tamasyan, G Sh},
  journal={Optimization},
  volume={60},
  number={1-2},
  pages={153--177},
  year={2011},
  publisher={Taylor \& Francis}
}

@article{OPTCOHEN20181138,
title = {Energy dissipating flows for solving nonlinear eigenpair problems},
journal = {Journal of Computational Physics},
volume = {375},
pages = {1138-1158},
year = {2018},
issn = {0021-9991},
doi = {https://doi.org/10.1016/j.jcp.2018.09.012},
url = {https://www.sciencedirect.com/science/article/pii/S0021999118306065},
author = {Ido Cohen and Guy Gilboa}
}

@article{calder2018game,
  title={The game theoretic p-Laplacian and semi-supervised learning with few labels},
  author={Calder, Jeff},
  journal={Nonlinearity},
  volume={32},
  number={1},
  pages={301},
  year={2018},
  publisher={IOP Publishing}
}

@article{shi2016weighted,
  title={Weighted graph laplacian and image inpainting},
  author={Shi, Zuoqiang and Osher, Stanley and Zhu, W},
  journal={Journal of Scientific Computing},
  volume={577},
pages={1-13},
  year={2016}
}

@Inbook{Debnath2012,
author="Debnath, Lokenath",
title="First-Order, Quasi-linear Equations and Method of Characteristics",
bookTitle="Nonlinear Partial Differential Equations for Scientists and Engineers",
year="2012",
publisher="Birkh{\"a}user Boston",
address="Boston",
pages="201--226",
isbn="978-0-8176-8265-1",
doi="10.1007/978-0-8176-8265-1_3",
url="https://doi.org/10.1007/978-0-8176-8265-1_3"
}

@article{cohen2024functional,
title = {Functional Dimensionality of Koopman Eigenfunction Space},
journal = {Results in Applied Mathematics},
volume = {26},
pages = {100585},
year = {2025},
issn = {2590-0374},
doi = {https://doi.org/10.1016/j.rinam.2025.100585},
url = {https://www.sciencedirect.com/science/article/pii/S2590037425000494},
author = {Ido Cohen and Eli Appleboim and Gershon Wolansky}
}

@book{kutz2016dynamic,
  title={Dynamic mode decomposition: data-driven modeling of complex systems},
  author={Kutz, J Nathan and Brunton, Steven L and Brunton, Bingni W and Proctor, Joshua L},
  year={2016},
  publisher={SIAM}
}

@book{tu2013dynamic,
  title={Dynamic mode decomposition: Theory and applications},
  author={Tu, Jonathan H},
  year={2013},
  publisher={Princeton University}
}

@article{hemati2017biasing,
  title={De-biasing the dynamic mode decomposition for applied Koopman spectral analysis of noisy datasets},
  author={Hemati, Maziar S and Rowley, Clarence W and Deem, Eric A and Cattafesta, Louis N},
  journal={Theoretical and Computational Fluid Dynamics},
  volume={31},
  number={4},
  pages={349--368},
  year={2017},
  publisher={Springer}
}

@article{dawson2016characterizing,
  title={Characterizing and correcting for the effect of sensor noise in the dynamic mode decomposition},
  author={Dawson, Scott TM and Hemati, Maziar S and Williams, Matthew O and Rowley, Clarence W},
  journal={Experiments in Fluids},
  volume={57},
  number={3},
  pages={42},
  year={2016},
  publisher={Springer}
}

@inproceedings{brokman2024enhancing,
  title={Enhancing neural training via a correlated dynamics model},
  author={Brokman, Jonathan and Betser, Roy and Turjeman, Rotem and Berkov, Tom and Cohen, Ido and Gilboa, Guy},
  booktitle={International Conference on Learning Representations},
  volume={2024},
  pages={48186--48220},
  year={2024}
}

@article{azencot2019consistent,
  title={Consistent dynamic mode decomposition},
  author={Azencot, Omri and Yin, Wotao and Bertozzi, Andrea},
  journal={SIAM Journal on Applied Dynamical Systems},
  volume={18},
  number={3},
  pages={1565--1585},
  year={2019},
  publisher={SIAM}
}

@article{lusch2018deep,
  title={Deep learning for universal linear embeddings of nonlinear dynamics},
  author={Lusch, Bethany and Kutz, J Nathan and Brunton, Steven L},
  journal={Nature communications},
  volume={9},
  number={1},
  pages={4950},
  year={2018},
  publisher={Nature Publishing Group UK, London}
}

@misc{miller2023eigenvalue,
  title  = {Eigenvalue Initialisation and Regularisation for {K}oopman Autoencoders},
  author = {Miller, Jack W. and O'Neill, Charles and Constantinou, Navid C. and Azencot, Omri},
  year   = {2023},
  note   = {ICLR 2023},
  howpublished = {OpenReview},
  url    = {https://openreview.net/forum?id=6TugHflAGRU}
}

@inproceedings{azencot2020forecasting,
  title={Forecasting sequential data using consistent Koopman autoencoders},
  author={Azencot, Omri and Erichson, N Benjamin and Lin, Vanessa and Mahoney, Michael},
  booktitle={International conference on machine learning},
  pages={475--485},
  year={2020},
  organization={PMLR}
}

@article{nayak2025temporally,
  title={Temporally consistent Koopman autoencoders for forecasting dynamical systems},
  author={Nayak, Indranil and Chakrabarti, Ananda and Kumar, Mrinal and Teixeira, Fernando L and Goswami, Debdipta},
  journal={Scientific Reports},
  volume={15},
  number={1},
  pages={22127},
  year={2025},
  publisher={Nature Publishing Group UK, London}
}

@article{rowley2009spectral,
  title={Spectral analysis of nonlinear flows},
  author={Rowley, Clarence W and Mezi{\'c}, Igor and Bagheri, Shervin and Schlatter, Philipp and Henningson, Dan S},
  journal={Journal of fluid mechanics},
  volume={641},
  pages={115--127},
  year={2009},
  publisher={Cambridge University Press}
}
\end{document}